\documentclass{amsart}
\usepackage[nohug]{diagrams}
\usepackage{graphicx}
\usepackage[mathscr]{eucal}
\usepackage{amssymb, amsmath, amsthm}
\usepackage{amsfonts}
\usepackage{latexsym}
\usepackage{tabularx,array}
\usepackage{enumerate}
\usepackage[all,arc]{xy}
\usepackage{latexsym}
\usepackage{mathrsfs}
\usepackage{color}
\usepackage{mathtools}
\usepackage{leftidx}
\usepackage{bm}
\usepackage{url}
\newfont{\msam}{msam10}

\newtheorem{theorem}[]{Theorem}
\newtheorem{proposition}[]{Proposition}
\newtheorem{corollary}[]{Corollary}
\newtheorem{lemma}[]{Lemma}

\theoremstyle{definition}
\newtheorem{definition}[]{Definition}
\newtheorem{defn}[theorem]{Definition}
\newtheorem{remark}[theorem]{Remark}
\newtheorem{example}[]{Example}
\newtheorem{conj}[]{Conjecture}

\def\remark{\noindent\textbf{Remark.}}
\let\nc\newcommand

\def\bthm{\begin{theorem}}
\def\ethm{\end{theorem}}
\def\blemma{\begin{lemma}}
\def\elemma{\end{lemma}}
\def\bproof{\begin{proof}}
\def\eproof{\end{proof}}
\def\bprop{\begin{proposition}}
\def\eprop{\end{proposition}}
\def\bcor{\begin{corollary}}
\def\ecor{\end{corollary}}
\def\bconj{\begin{conj}}
\def\econj{\end{conj}}
\nc{\la}{\label}
\def\O{\mathcal{O}}
\def\Z{\mathbb{Z}}

\def\Q{\mathbb{Q}}
\def\c{\mathbb{C}}

\def\L {\boldsymbol{L}}
\def\R{\mathbb{R}}

\def\Com{\mathtt{Com}}
\def\Vect{\mathtt{Vect}}

\def\Top{\mathtt{Top}}

\def\Mod{\mathtt{Mod}}

\def\cAlg{\mathtt{Comm\,Alg}}
\def\scAlg{\mathtt{sComm\,Alg}}

\def\DGCA{\mathtt{dgComm\,Alg}}
\def\cDGA{\mathtt{dgComm\,Alg}}

\def\Ho{{\mathtt{Ho}}}

\nc{\Ob}{{\rm Ob}}
\nc{\Hom}{{\rm{Hom}}}
\nc{\Homcont}{{\mathcal{H}om}}
\nc{\HOM}{\underline{\rm{Hom}}}
\nc{\DER}{\underline{\rm{Der}}}
\nc{\END}{\underline{\rm{End}}}
\nc{\bSym}{\mathbf{Sym}}
\nc{\Ext}{{\rm{Ext}}}
\nc{\Rep}{{\rm{Rep}}}
\nc{\DRep}{{\rm{DRep}}}
\nc{\NCRep}{\widetilde{\rm{Rep}}}
\nc{\RAct}{{\rm{RAct}}}
\nc{\bs}{\backslash}
\nc{\ob}{{\tt{Obs}}}
\nc{\CE}{\mathcal{C}}
\nc{\TP}{{T\!P}}
\nc{\un}{\underline{n}}
\nc{\um}{\underline{m}}
\nc{\rn}{\langle n \rangle}
\nc{\nn}{{{\natural} {\natural}}}
\nc{\n}{{{\natural}}}
\nc{\A}{\mathbb A}
\nc{\B}{{\mathrm{B}}}
\nc{\Ba}{\overline{\mathrm{B}}}
\nc{\bC}{\overline{C}}
\nc{\bOmega}{\boldsymbol{\Omega}}
\nc{\bB}{\boldsymbol{B}}
\nc{\EXT}{\underline{\rm{Ext}}}
\nc{\TOR}{\underline{\rm{Tor}}}
\def\H{\mathrm H}

\def\HR{\mathrm{HR}}

\def\T{\mathrm T}

\nc{\End}{{\rm{End}}}
\nc{\GL}{{\rm{GL}}}
\nc{\gl}{{\mathfrak{gl}}}
\nc{\rgl}{\overline{{\mathfrak{gl}}}}
\nc{\g}{{\mathfrak{g}}}
\nc{\h}{{\mathfrak{h}}}
\nc{\PGL}{{\rm{PGL}}}
\nc{\SL}{{\rm{SL}}}
\nc{\sll}{\mathfrak{sl}}
\nc{\cn}{ \mbox{\rm c\^{o}ne} }
\nc{\PSL}{{\rm{PSL}}}
\nc{\ad}{{\rm{ad}}}
\nc{\Ad}{{\rm{Ad}}}
\nc{\dlim}{\varinjlim}
\nc{\plim}{\varprojlim}
\nc{\colim}{{{\rm colim}}}

\newcommand{\Loc}{{\sl L}{\rm oc}}

\newcommand{\HH}{{\rm{HH}}}

\newcommand{\Tor}{{\rm{Tor}}}
\newcommand{\Spec}{{\rm{Spec}}}

\newcommand{\Sym}{{\rm{Sym}}}

\newcommand{\id}{{\rm{Id}}}
\newcommand{\Der}{{\rm{Der}}}

\newcommand{\Ker}{{\rm{Ker}}}

\newcommand{\diag}{{\rm{diag}}}
\newcommand{\Coker}{{\rm{Coker}}}

\newcommand{\im}{{\rm{Im}}}

\newcommand{\sonto}{\,\stackrel{\sim}{\twoheadrightarrow}\,}
\def\bs{\backslash}

\def\ve{\mathtt{Vect}}

\def\sGr{\mathtt{sGr}}
\def\Gr{\mathtt{Gr}}

\def\set{\mathtt{Set}}
\def\sset{\mathtt{sSet}}
\def\lgr{\mathbb{G}}

\def\Char{\mathrm{Char}}
\def\DChar{\mathrm{DChar}}

\newcommand{\rar}{\xrightarrow{}}

\nc{\env}{\mathrm{End}(V)}
\nc{\FT}{\mathcal{C}}

\numberwithin{equation}{section}
\numberwithin{theorem}{section}
\numberwithin{lemma}{section}
\numberwithin{proposition}{section}
\numberwithin{definition}{section}
\numberwithin{corollary}{section}
\numberwithin{example}{section}
\newcommand{\ab}{\mathrm{ab}}
\def\cC{\mathrm{C}}

\def\cO{\mathcal O}
\newcommand{\sMod}{{\tt sMod}}

\newcommand{\op}{{\rm op}}

\def\bdf{\begin{defn}}
\def\edf{\end{defn}}

\def\brm{\begin{remark}}
\def\erm{\end{remark}}

\newcommand{\sCommAlg}{\mathtt{sComm\,Alg}}
\theoremstyle{definition}
\def\bdf{\begin{definition}}
\def\edf{\end{definition}}
\def\cQ{\mathcal{Q}}

\newcommand{\sAff}{\mathtt{sAff}}
\newcommand{\Aff}{\mathtt{Aff}}
\newcommand{\CommAlg}{\mathtt{Comm\,Alg}}

\newcommand{\dAff}{\mathtt{dAff}}
\newcommand{\st}{{\rm st}}
\newcommand{\bS}{{\mathbb S}}

\newcommand{\m}{{\mathfrak m}}

\def\arbreBA{\vcenter{\xymatrix@R=2pt@C=2pt{
&&&&\\
&&&*{}\ar@{-}[ul] & \\
&&*{}\ar@{-}[uurr] \ar@{-}[uull] \ar@{-}[d]     &&\\
&&&&
}}}

\def\arbreAB{\vcenter{\xymatrix@R=2pt@C=2pt{
&&&&\\
&*{}\ar@{-}[ur] &&& \\
&&*{}\ar@{-}[uurr] \ar@{-}[uull] \ar@{-}[d]     &&\\
&&&&
}}}

\def\arbreABC{\vcenter{\xymatrix@R=1pt@C=1pt{
&&&&&&\\
&*{}\ar@{-}[ur] &&&&& \\
&&*{}\ar@{-}[uurr] &&&&\\
&&&*{}\ar@{-}[uuurrr] \ar@{-}[uuulll] \ar@{-}[d] &&&\\
&&&&&&
}}}

\def\arbreBAC{\vcenter{\xymatrix@R=1pt@C=1pt{
&&&&&&\\
&&&*{}\ar@{-}[ul] &&& \\
&&*{}\ar@{-}[uurr] &&&&\\
&&&*{}\ar@{-}[uuurrr] \ar@{-}[uuulll] \ar@{-}[d] &&&\\
&&&&&&
}}}

\def\arbreACB{\vcenter{\xymatrix@R=1pt@C=1pt{
&&&&&&\\
&*{}\ar@{-}[ur] &&&&& \\
&&&&*{}\ar@{-}[uull] &&\\
&&&*{}\ar@{-}[uuurrr] \ar@{-}[uuulll] \ar@{-}[d] &&&\\
&&&&&&
}}}

\def\arbreBCA{\vcenter{\xymatrix@R=1pt@C=1pt{
&&&&&&\\
&&&&&*{}\ar@{-}[ul] & \\
&&*{}\ar@{-}[uurr] &&&&\\
&&&*{}\ar@{-}[uuurrr] \ar@{-}[uuulll] \ar@{-}[d] &&&\\
&&&&&&
}}}

\def\arbreCAB{\vcenter{\xymatrix@R=1pt@C=1pt{
&&&&&&\\
&&&*{}\ar@{-}[ur] &&& \\
&&&&*{}\ar@{-}[uull] &&\\
&&&*{}\ar@{-}[uuurrr] \ar@{-}[uuulll] \ar@{-}[d] &&&\\
&&&&&&
}}}

\def\arbreCBA{\vcenter{\xymatrix@R=1pt@C=1pt{
&&&&&&\\
&&&&&*{}\ar@{-}[ul] & \\
&&&&*{}\ar@{-}[uull] &&\\
&&&*{}\ar@{-}[uuurrr] \ar@{-}[uuulll] \ar@{-}[d] &&&\\
&&&&&&
}}}

\def\arbreACA{\vcenter{\xymatrix@R=1pt@C=1pt{
&&&&&&\\
&*{}\ar@{-}[ur] &&&&*{}\ar@{-}[ul] & \\
&&&&&&\\
&&&*{}\ar@{-}[uuurrr] \ar@{-}[uuulll] \ar@{-}[d] &&&\\
&&&&&&
}}}

\begin{document}

\title{Vanishing theorems for representation homology\\ and the derived
cotangent complex}
\author{Yuri Berest}
\address{Department of Mathematics,
Cornell University, Ithaca, NY 14853-4201, USA}
\email{berest@math.cornell.edu}
\author{Ajay C. Ramadoss}
\address{Department of Mathematics,
Indiana University,
Bloomington, IN 47405, USA}
\email{ajcramad@indiana.edu}
\author{Wai-kit Yeung}
\address{Department of Mathematics,
Indiana University,
Bloomington, IN 47405, USA}
\email{yeungw@iu.edu}
\begin{abstract}
Let $G$ be a reductive affine algebraic group defined over a field $k$ of characteristic zero.  In this paper, we study the  cotangent complex of the derived $G$-representation scheme $ \DRep_G(X)$ of a pointed connected topological space $X$. We use an (algebraic version of) unstable Adams spectral sequence 
relating the cotangent homology of $ \DRep_G(X) $ to the representation homology $ \HR_\ast(X,G) := 
\pi_\ast {\mathcal O}[\DRep_G(X)] $ to prove some vanishing theorems for groups and geometrically interesting spaces. Our examples include virtually free groups, Riemann surfaces, link complements in $ \R^3 $ and generalized lens spaces. In particular, for any f.g. virtually free group $ \Gamma $, we show
that $\, \HR_i(\B\Gamma, G) = 0 \,$ for all $ i > 0 $. For a closed Riemann surface $\Sigma_g $ of genus $ g \ge 1 $, we have $\, \HR_i(\Sigma_g, G) = 0 \,$ for all $ i > \dim G $. The {\it sharp} vanishing bounds for $ \Sigma_g $ depend actually on the genus: we conjecture that if $ g = 1 $,  then
$\, \HR_i(\Sigma_g, G) = 0 \,$ for $ i > {\rm rank}\,G $, and if $ g \ge 2 $, then
$\, \HR_i(\Sigma_g, G) = 0 \,$ for $ i > \dim\,{\mathcal Z}(G) \,$, where $ {\mathcal Z}(G)  $
is the center of $G$. We prove these bounds locally on the smooth locus of the representation 
scheme $ \Rep_G[\pi_1(\Sigma_g)]\,$ in the case of complex connected reductive groups. One important consequence of our results is the existence of a well-defined $K$-theoretic virtual fundamental class for $ \DRep_G(X)$ in the sense of Ciocan-Fontanine and Kapranov \cite{CK}. We give a new ``Tor formula'' for this class in terms of functor homology.

\end{abstract}
\maketitle
\section{Introduction and main results}
Topologically, it is natural to think of discrete groups as homotopy types of aspherical spaces\footnote{Recall that a (pointed connected) topological space $X$ is {\it aspherical} if all its higher homotopy groups vanish: i.e., $\pi_i(X) = 0 $ for $ i >1 $.
Assigning to a discrete group $ \Gamma $ the homotopy type of its classifying space $ {\rm B} \Gamma $ gives a fully faithful functor $ \Gr \to \Ho(\Top_{0,*}) $ from groups to the homotopy category of pointed connected spaces. The essential image of this functor comprises precisely the  aspherical spaces.}. Given some algebraic or combinatorial property of a group, or a group-theoretic construction, it is then natural to ask whether one can extend it to more general spaces. A most familiar example is the ordinary homology of spaces: by a classical theorem of Kan \cite{Kan0}, it is an extension of abelianization of groups. Another less known but equally important example, which plays a fundamental role in homotopy theory, is the Bousfield-Kan completion of a space: it extends the operation of pronilpotent completion on the category of groups (see \cite{BK}).

In this paper, we will study {\it representation homology of spaces}, which is a topological extension of another classical construction from group theory: a representation variety of a group. Given a discrete group $ \Gamma $ and an affine algebraic group $G$ defined over a field $k$, one can consider the space $ \Rep_G(\Gamma) $ of all representations of $ \Gamma $ in $G$. This space carries a natural structure of an algebraic variety (more precisely, an affine $k$-scheme) called the $G$-{\it representation variety} of $\Gamma $. The geometry of $ \Rep_G(\Gamma) $ reflects the properties of the group $ \Gamma $, the geometry of $G$ as well as the algebraic structure of representations of $ \Gamma $ in $G$. The group $G$ acts naturally on $ \Rep_G(\Gamma) $ by conjugation, with orbits corresponding to the equivalence classes of representations. The categorical quotient $  \Char_G(\Gamma) := \Rep_G(\Gamma)/\!/G $ is called the $G$-{\it character variety} of $\Gamma$; when $k$ is algebraically closed and $ \Gamma $ is finitely generated, $ \Char_G(\Gamma) $ parametrizes the closed $G$-orbits in $ \Rep_G(\Gamma) $ which correspond to the semisimple representations of $\Gamma $ in $G$ (see, e.g., \cite{LuM}).

Representation varieties and associated character varieties play an important role in many areas of mathematics, most notably in geometry and topology\footnote{See, for example, the recent survey papers \cite{Coh}, \cite{Sh} and references therein.}. It is therefore natural to ask how to describe these algebro-geometric objects in homotopy-theoretic terms which would make sense for  arbitrary topological spaces. Put concisely, the question is: {\it What is a ``representation variety of a space''}\,?

There is an obvious answer to this question: given a (pointed)  space $X$, just take $ \Rep_G $ of its fundamental group $ \pi_1(X) $. If $ \pi_1(X) $ is nontrivial, this gives a geometric invariant of $X$ that often encodes an interesting topological information\footnote{Representation varieties of fundamental groups are particularly rich and geometrically interesting objects in the case of surfaces and 3-dimensional manifolds, see, e.g., \cite{Go} and \cite{CS, Sh}.}. By its very definition, however, the variety $ \Rep_G[\pi_1(X)] $
depends only on the fundamental group and contains no information about the higher homotopy groups of $X$; in particular, it is trivial when $X$ is simply connected.

The first ``non-obvious'' answer to the above question was proposed -- to the best of our knowledge -- by Kapranov \cite{K}.
He observed that if $X$ is a finite pointed connected CW-complex, the representation space  $ \Rep_G[\pi_1(X)] $ can be naturally identified with the moduli space $ \Loc_G(X) $ of $G$-local systems on $X$ (trivialized at the basepoint). The homotopy type of $X$ can be represented simplicially by its \v{C}ech nerve $ N_X $, and the moduli space $  \Loc_G(X) $ can then be  expressed as $ [N_X, BG] $, the set of homotopy classes of simplicial maps from $N_X$ to the
simplicial classifying space of $G$. Replacing $BG$ in this construction
by a simplicial DG scheme $ \boldsymbol{R}BG $, which plays the role of injective resolution of $BG$ in the category of simplicial DG schemes, Kapranov introduced an affine DG scheme
$  \boldsymbol{R} \Loc_G(X) $, which he called the derived moduli space of $G$-local systems on $X$.
The homotopy type of $  \boldsymbol{R} \Loc_G(X) $ in the category of DG schemes depends only on the homotopy type of $X$  and, unlike $ \Loc_G(X) $, it is sensitive to the higher homotopy structure of $X$.
Since $\, \pi_0[\boldsymbol{R} \Loc_G(X)] \cong
\Loc_G(X) \cong  \Rep_G[\pi_1(X)] $, the DG scheme $  \boldsymbol{R} \Loc_G(X) $ can be regarded as a higher homotopical extension of the representation variety $ \Rep_G[\pi_1(X)] $. Clearly, by replacing  \v{C}ech nerves with arbitrary (pointed connected) simplicial sets, one can extend this construction to  arbitrary (pointed connected) spaces.
The problem, however, is that the `injective resolution' $ \boldsymbol{R}BG $ is `too big' and complicated, and it is hard to understand $  \boldsymbol{R} \Loc_G(X) $ (in particular, to compute $ \pi_i[\boldsymbol{R} \Loc_G(X)] $ for $ i> 0$) even in simplest examples. Kapranov's definition of 
the derived moduli space of $G$-local systems was refined and generalized within the framework of
derived algebraic geometry by To\"{e}n, Vezzosi, Pridham and others (see, e.g., \cite{To1, TV05, TV08,  PTVV13, Pri08, Pri13, Pri132}
and \cite{To} for a general overview)\footnote{We will briefly review the To\"{e}n-Vezzosi construction and compare it to the original Kapranov's one in the Appendix.}; however, explicit computations seem still to be missing.

A different approach, closer in spirit to classical homotopy theory, was developed by the authors in 
\cite{BRY}. 
Instead of looking for a generalization of representation variety for a fixed space $X$
by resolving the classifying scheme of the algebraic group $G$, we fix the group $G$ and construct the derived representation scheme $ \DRep_G(X) $ by `resolving' the space $X$. More precisely, our starting point is the simple observation that the functor of points $\, G: \cAlg_k \to \Gr \,$ defining the algebraic group $G$ has a left adjoint
\begin{equation}
\la{lgrsc}
(\,\mbox{--}\,)_G\,:\,\mathtt{Gr} \to \cAlg_k\ 
\end{equation}
which \mbox{---} when applied to a given group $ \Gamma $ \mbox{---} represents the affine scheme $ \Rep_G(\Gamma) $: i.e. $\Gamma_G = \O[\Rep_G(\Gamma)]\,$. Therefore we can replace \mbox{---} in an equivalent way \mbox{---} the representation scheme functor
$ \Rep_G $ by  \eqref{lgrsc} which we call the {\it representation functor in} $G$.

The  functor   \eqref{lgrsc} extends naturally to
the category $ \sGr $ of simplicial groups, taking values in the category
$ \mathtt{s}\cAlg_k $ of simplicial commutative algebras. Both $ \sGr $
and $ \mathtt{s}\cAlg_k $ are (simplicial) model categories, with weak equivalences
being the weak homotopy equivalences of underlying simplicial sets. Although the representation functor $\,(\,\mbox{--}\,)_{G}:\, \sGr \to \mathtt{s}\cAlg_k\,$ is not homotopic \mbox{---} it does not preserve weak equivalences and hence, does not induce a functor between  homotopy categories \mbox{---} it has a well-behaved (total) left derived functor
\begin{equation}
\la{Lrep}
\L(\,\mbox{--}\,)_{G} :\ \Ho(\mathtt{sGr}) \rar \Ho(\mathtt{s}\cAlg_k)\,\text{.}
\end{equation}
The derived functor  \eqref{Lrep} is uniquely determined by \eqref{lgrsc} and provides, in a sense, the `closest' 
approximation to the representation functor at the level of homotopy categories.

Now, for a fixed simplicial group $\Gamma \in \sGr $, we set\footnote{Abusing notation, we will
also write $ \mathcal{A}_{G}(\Gamma) $ for  a specific representative (model) of  $ \L(\Gamma)_{G} $ in  $\mathtt{s}\cAlg_k$.}
$\,\mathcal{A}_{G}(\Gamma) := \L(\Gamma)_{G}\,$ and formally define the {\it derived representation
scheme} $\, \DRep_{G}(\Gamma) $ as $ \Spec[\mathcal{A}_{G}(\Gamma)]\, $, i.e. the simplicial
algebra $\mathcal{A}_{G}(\Gamma)$ viewed as an object of the opposite category $\Ho(\scAlg_k)^{\mathrm{op}}$.
Following \cite{BRY}, we call the  homotopy groups of  $ \DRep_{G}(\Gamma) $ the {\it representation homology of $\Gamma$ in $G\,$} and write
$$
\HR_\ast(\Gamma,G) := \L_\ast(\mathbf{\Gamma})_{G}\ ,
$$
where $ \L_\ast(\,\mbox{--}\,)_{G} $ denotes the composition of functors $ \pi_\ast \L(\,\mbox{--}\,)_{G} \cong 
\H_\ast[N\L(\,\mbox{--}\,)_{G}]\,$ (see Section~\ref{S2.1}).

By comparing the universal mapping properties, it is easy to check that, when extended to simplicial groups, the  functor $ (\,\mbox{--}\,)_G $  commutes with $ \pi_0 $;
hence, for any $\Gamma \in \sGr $, there is a natural isomorphism in $ \cAlg_k $:
\begin{equation}
\la{commpi}
\HR_0(\Gamma,G)\,\cong\, [\pi_0(\Gamma)]_{G}
\end{equation}
In particular, if $\Gamma \in \Gr $ is a constant simplicial group, we have
$\,\HR_0(\Gamma,G) \cong \Gamma_{G} \,$, which
justifies our notation and terminology for $ \DRep_{G}(\Gamma) $.

Now, the reason why we choose to work with simplicial groups is the fundamental theorem
of D. Kan \cite{Kan1} that identifies the homotopy types of simplicial groups with those of pointed connected spaces. To be precise, the Kan Theorem asserts that the category of simplicial groups is Quillen equivalent to the category $ \sset_0 $ of reduced simplicial sets, which is, in turn, Quillen equivalent to the category $ \Top_{0,*} $ of pointed connected (CGWH)  spaces. As a result, we have natural equivalences of homotopy categories
\begin{equation}\la{tsg}
\Ho(\mathtt{Top}_{0,\ast}) \,\cong\, \Ho(\sset_0)\,\cong\, \Ho(\sGr)\ .
\end{equation}
This leads us to the main definition.
\begin{definition}
\la{DRepX}
For a space $X \in \Top_{0,*} $, we define the {\it derived representation scheme}
$\,\DRep_G(X)\,$ to be $\,\DRep_G(\mathbf{\Gamma}X)\,$, where $ \mathbf{\Gamma}X$  is a(ny) simplicial group model\footnote{that is, a simplicial group that corresponds to $X$ under the Kan equivalence \eqref{tsg}.} of $X$. The {\it representation homology of $X$ in $G$} is then defined by
\begin{equation}
\la{hrX}
\HR_\ast(X,G) := \L_\ast (\mathbf{\Gamma}X)_{G}  := \pi_\ast \L(\mathbf{\Gamma}X)_{G} \ .
\end{equation}
\end{definition}

The advantage of our definition (compared to Kapranov's)
lies in its flexibility in the choice of simplicial group model $\mathbf{\Gamma}X$. In fact, there are two natural choices (both due to Kan) that play a role in practical computations: the Kan loop group
$ \lgr{X} $, which is a `large' semi-free simplicial group depending functorially on $X$,
and for finite CW-complexes, a `small' geometric model $ \mathbf{\Gamma}X$,
the structure of which reflects the cell structure of $X$ (see Section~\ref{simpgr} below).

For a detailed discussion of representation homology we refer the reader
to \cite{BRY}. Here, we only recall a few basic properties which we will need in the present paper:

\vspace*{1ex}

\begin{enumerate}
\item[(1)] $ \HR_*(X,G) $ is a graded commutative algebra, with
$ \HR_0(X,G) $  naturally isomorphic to $ [\pi_1(X)]_G = \cO[\Rep_G(\pi_1(X))]$,
the coordinate ring of the representation scheme $ \Rep_G[\pi_1(X)] $.
The last isomorphism is the composition of \eqref{commpi} with the natural isomorphism $ \pi_0(\mathbf{\Gamma}X) \cong
\pi_1(X) $.

\vspace{1ex}

\item[(2)]
If $ X $ is a $ K(\Gamma, 1)$-space (e.g., $ X = {\rm B} \Gamma $ for
a discrete group $\Gamma$), then $ \HR_\ast(X,G) \cong \HR_\ast(\Gamma,G) $
(see \cite[Corollary~4.1]{BRY}); in general, however,  $\,\HR_\ast(X,G) \not\cong
\HR_\ast(\pi_1(X), G)\,$.

\vspace{1ex}

\item[(3)] If $X$ is  simply connected and $ k=\Q $, then
$ \HR_\ast(X,G) $ is a {\it rational} homotopy invariant of $X$. In fact,
when $X$ has finite rational type, the main theorem of
\cite{BRY} describes $ \HR_\ast(X,G) $ explicitly in terms of
 Quillen's Lie algebra model of $X$ (see \cite[Theorem~6.1]{BRY}).

\vspace{1ex}

\item[(4)] The derived representation functor \eqref{Lrep} preserves homotopy pushouts (\cite[Lemma~4.2]{BRY}). This implies, in particular, that
$$
\HR_\ast(X \vee Y, G) \cong \HR_\ast(X, G) \otimes \HR_\ast(Y, G)
$$
for any spaces  $ X, Y \in \Top_{0,\ast} $.

\vspace{1ex}

\item[(5)] For any space $X$ (not necessarily pointed), there is a natural isomorphism of  algebras
$$
\HR_\ast(\Sigma(X_{+}), G) \cong \HH_\ast(X, \cO(G))\ ,
$$
where the right-hand side is the higher Hochschild homology of $X$ with coefficients in the commutative algebra $ \cO(G) $ as defined in \cite{P1} (see \cite[Theorem 5.2]{BRY}).
\end{enumerate}

\vspace{1ex}

The aim of the present paper is to study the linearization of representation
homology with a view to proving some vanishing theorems for groups and  geometrically interesting spaces. If $X$ is a pointed connected space, any representation of its fundamental group $ \varrho: \pi_1(X) \to G(k)$ determines  an augmentation $\, \mathcal{A}_{G}(X) \to k \,$ of the derived representation algebra $ \mathcal{A}_{G}(X) $. The derived linearization functor
$ \L{\mathcal Q}: \Ho(\scAlg_{k/k}) \to \Ho(\sMod_k) $ applied to
$ (\mathcal{A}_{G}(X), \varrho) $ gives a simplicial $k$-module
$ \T_{\varrho}^*\DRep_G(X) $, which is called \mbox{---} as suggested by its notation
\mbox{---} the derived cotangent complex
of $ \DRep_G(X) $ at $\,\varrho \in \Rep_G[\pi_1(X)](k) $. A fairly straightforward
calculation (see Theorem~\ref{LinRep}) shows 
$$
\H_i[\mathrm{T}^{\ast}_{\varrho}\DRep_G(X)]\,\cong\, \left\{
                                                 \begin{array}{ll}
                                                    Z^1(\pi_1(X), \mathrm{Ad}\,\varrho)^{\ast} & \text{if}\ \, i=0\\*[1ex]

                                                    \H_{i+1}(X,\mathrm{Ad}^{\ast}\!\varrho) & \text{if}\ \,i>0 \ \\
                                                 \end{array}
                                                 \right.
$$
where $ \H_\ast(X, \mathrm{Ad}^{\ast}\!\varrho)  $ is the homology of
$X$ with coefficients in the local system associated to the coadjoint representation $  \mathrm{Ad}^{\ast} \varrho:\, \pi_1(X) \to G(k) \to \GL_k(\g^*)\,$. Perhaps more interesting is the following simple formula for the Euler characteristic of $ \T_{\varrho}^*\DRep_G(X) $ for a finite reduced CW-complex $X$ (see Proposition~\ref{Euler}):
$$
\chi[\T^{\ast}_{\varrho}\DRep_G(X)]\,=\,[1-\chi_{{\scriptsize \rm top}}(X)] \dim G\ .
$$
This formula shows that $\,\chi[\T^{\ast}_{\varrho}\DRep_G(X)]\,$ is independent of $ \varrho $,
which may be interpreted geometrically as the `absence of singularities' of the derived scheme
$ \DRep_G(X) $ in agreement with Kontsevich's `hidden smoothness' philosophy in derived algebraic geometry ({\it cf.} \cite{To}).

Now, assume that $G$ is a reductive group over $k$. Then, it is natural to consider the $G$-invariant part of the representation functor \eqref{lgrsc}: indeed,
for a given space $X$, the $G$-invariant subalgebra $\,[\pi_1(X)]^G_G $ of $[\pi_1(X)]_G $ represents the classical character scheme $ \Char_G[\pi_1(X)]$.
Just as  the representation functor itself, its invariant subfunctor $\, (\,\mbox{--}\,)^{G}_{G} \,$ has a left derived functor $\,\L (\,\mbox{--}\,)^{G}_{G}\,$, which one can use to define the derived character scheme $ \DChar_G(X) $. Computing the cotangent homology
for  $ \DChar_G(X) $ is, however, a much more delicate problem than that for $ \DRep_G(X) $. Our first main result in this paper provides such a computation for certain `good'  representations of $\pi_1(X)$.
To be precise, following \cite{JM} (see also \cite{S}), we call
a representation $ \varrho: \Gamma \to G $ {\it good}\, if its $G$-orbit
in $ \Rep_G(\Gamma) $ is closed and its stabilizer coincides with
the center of $G$. One can show that, at least when $k$ is algebraically
closed, every good representation in $G$ is irreducible but the converse is
not always true.
\bthm[see Theorem~\ref{lincharhom}]
\la{T0}
Let $G$ be a connected reductive affine algebraic group over $ \c $, and
let $X$ be a  space such that $\pi_1(X)$ is finitely generated. Then

\vspace{1ex}

$(a)$ for every completely irreducible representation $\varrho :\,\pi_1(X) \rar G $, there are natural  maps
 \begin{equation*}
 \la{ecothom1}
 \H_i[\T^{\ast}_{\varrho}\DChar_G(X)] \rar \H_{i+1}(X, \mathrm{Ad}^{\ast}\varrho)\,,\,\,\, \forall\, i \geq 0 \, ,
 \end{equation*}

 $(b)$ if a representation $ \varrho :\,\pi_1(X) \rar G $ is good, then the above maps are all isomorphisms.
\ethm
We remark that Theorem~\ref{T0} is a natural generalization of the main result of \cite{S} on tangent spaces of the classical character scheme $\, \Char_G(\Gamma)\, $ ({\it cf.} \cite[Theorem~53]{S}); however, for the derived scheme $\DChar_G(X)$,  the proof  is rather more intricate.

Next, we explain the relation between representation homology and the
derived cotangent complex. This relation is given
in the form of the spectral sequence
\begin{equation}
\la{CKsp}
E_{p,q}^2 = \pi_{p+q}[\Sym^p(\T_{\varrho}^*\DRep_G(X))]\ \ \Longrightarrow\ \ \widehat{\HR}_{p+q}(X,G)_{\varrho}
\end{equation}
which converges to the local $\m$-adic
completion of $\, \HR_\ast(X,G) $ at the maximal ideal corresponding to 
$ \varrho \in \Rep_G[\pi_1(X)](k) $.
In characteristic zero  ($ \Q \subseteq k $), it is known
that the derived cotangent complex of an augmented simplicial commutative algebra $A$ is a {\it stable} homotopy invariant, which is represented (in the stable homotopy category of $ \scAlg_{k/k}$) by the suspension spectrum $ \Sigma^{\infty}\! A $ of the unreduced suspension of $A$ (see \cite{Sch}). From this point of view, one can regard  \eqref{CKsp} as an algebraic analogue of an unstable Adams spectral sequence in classical topology ({\it cf.} \cite{BK1}). In the above form, the spectral sequence \eqref{CKsp} first appeared in \cite{CK} but it essentially
goes back to Quillen \cite{Q70}
(see Remark in  Section~\ref{S2} after Proposition~\ref{Quillen_specseq}).

We will use the spectral sequence \eqref{CKsp} to prove some
vanishing theorems for representation homology. Our examples include virtually
free groups, Riemann surfaces $ \Sigma_g $ for $ g \ge 1 $ as well as some non-aspherical spaces such as link complements in $ \R^3 $ and lens spaces. We summarize our main vanishing results in the following  four theorems.

\bthm[see Theorem~\ref{vfgroups}]
\la{T1}
Let $ \Gamma $ be a f.g. virtually free group. Then, for any affine algebraic group $G$
over a field $k$ of characteristic zero,
\begin{equation}
\la{vfgr}
\HR_i(\Gamma, G) = 0 \ ,\ \forall\, i > 0\ .
\end{equation}
In particular,  for any finite group $ \Gamma $, the higher representation homology vanishes.
\ethm

For the next two theorems, we assume that $ k = \c $ and $G$ is a  complex reductive group. We will say that a representation $ \varrho: \Gamma \to G $ is {\it smooth} if the corresponding (closed) point of the representation scheme
$ \Rep_G(\Gamma) $ is simple\footnote{that is, $\varrho $ belongs to a unique irreducible component of
$ \Rep_G(\Gamma) $ and the dimension of that component coincides with
$ \dim_k [ \T_\varrho \Rep_G(\Gamma)] \,$.}. If $\, \varrho \in  \Rep_G[\pi_1(X)] \,$, we write
$ \HR_\ast(X,G)_{\varrho} $ for the localization of $ \HR_\ast(X,G) $ at the maximal ideal
of $ [\pi_1(X)]_G $ corresponding to $ \varrho $.
\bthm[see Corollary~\ref{vrephomtorus} and Corollary~\ref{vrephomhighergenus}]
\la{T2}
Let $ \Sigma_g $ be a closed connected orientable surface of genus $ g \ge 1 $.
Then

\vspace*{1ex}

$(1)$ $\,\HR_i(\Sigma_g, G) = 0 \ $ for all $\ i > \dim G \,$

\vspace*{1ex}

$(2)$ if $\,g=1\,$, i.e. $ \Sigma_g = T^2 $ is a $2$-torus, then for every smooth  representation $ \varrho: \pi_1(T^2) \to G $,
$$
\HR_i(T^2, G)_\varrho = 0\quad  \mbox{for all}\ \, i > {\rm rank}\,G
$$

\vspace*{1ex}

$(3)$ if $\,g \ge 2 \,$, then for every smooth representation $ \varrho: \pi_1(\Sigma_g) \to G $,
$$
\HR_i(\Sigma_g, G)_{\varrho} = 0\quad  \mbox{for all}\ \,  i > \dim \mathcal{Z}(G)\ ,
$$
where $ \mathcal{Z}(G) $ is the center of  $G$.
\ethm
We remark that Part $(2)$ and Part $(3)$ of Theorem~\ref{T2}  give
\mbox{---} although locally \mbox{---} much sharper vanishing bounds for representation homology
than Part $(1)$. We expect that these local bounds are actually {\it global}, i.e. they hold for
{\it all}\, representations of $ \pi_1(\Sigma_g) $. To be precise, we propose
\bconj
\la{Conj1} {\it
For any complex connected reductive affine algebraic group $G$,

\vspace*{1ex}

$(i)$ $\,\HR_i(T^2, G) = 0\,$ for all $\,i > {\rm rank}\,G\,$;

\vspace*{1ex}

$(ii)$ if $\,g\ge 2\,$, then $\,\HR_i(\Sigma_g, G) = 0\,$ for all $\, i > \dim \mathcal{Z}(G)\,$.
}
\econj
Note that Part $(ii)$ of Conjecture~\ref{Conj1} implies that, for any complex 
{\it semisimple}\, group $G$,  $\,\HR_i(\Sigma_g, G) = 0 \,$ for all $\, i > 0 \,$. 
This, in turn, implies that the corresponding $K$-theoretic class
$$ 
[\DRep_G(\Sigma_g)]_K^{\rm vir} := \sum_{i \ge 0}\,(-1)^{i} \,[\HR_i(\Sigma_g, G)]
$$ 
called the virtual fundamental class of $\, \DRep_G(\Sigma_g) \,$, does not vanish in the
rational Grothendieck group $\, K_0[\Rep_G(\Sigma_g)]_{\Q} $ of 
$\,\Rep_G(\Sigma_g)$. We prove that this is indeed the case (see Corollary~\ref{vfclass}), 
which can be viewed as an indirect evidence for Conjecture~\ref{Conj1}.
As for other evidence \mbox{---} besides Theorem~\ref{T2} \mbox{---} we mention that Part $(i)$ is known to be true for $ \GL_n $ for all $ n\ge 1 $ as a consequence of \cite[Theorem~27]{BFR};  Part $(ii)$ has been verified for the surface groups of genus $2$ and $3$ in the case of
$ \GL_2 $ and $ \SL_2$, using the {\tt Macaulay2} computer software,
and it also holds for $ G = (\GL_1)^r $ for all $\,r \ge 1 \,$, which is an easy calculation (see
Example~\ref{conj1tori} in Section~\ref{S5}).

Next, we consider the link complements in $ \R^3 $. By a link $ L $ in $ \R^3 $ we mean a smooth  embedding of a finite disjoint union of copies of $ \bS^1 $ into $ \R^3 $. The link complement
$ \R^3\!\setminus\! L $ is then defined to be the complement of an (open) tubular neighborhood of the image of $L$.

\bthm[see Theorem \ref{linkcompl}]
\la{T3}
Let $ X := \R^3\!\setminus\! L $ be a link complement in $ \R^3 $.
Then

\vspace*{1ex}

$(1)$ $\,\HR_i(X, G)=0 \ $ for all $\ i> M := \max \{\dim_k \T_\varrho \Rep_{G} [\pi_1(X)]\,:\,\varrho \in \Rep_G[\pi_1(X)]\}\,$.

\vspace*{1ex}

$(2)$ For every smooth representation $\,\varrho:\,\pi_1(X) \rar G\,$,
$$
\HR_i(X, G)_{\varrho}=0 \ \, \mbox{for all}\ \,i>\dim_{\varrho}\Rep_G[\pi_1(X)]\,,
$$
where $\,\dim_{\varrho}\,$ denotes the local dimension of the irreducible component of 
$\,\Rep_G[\pi_1(X)] \,$ containing $\varrho$.

\vspace*{1ex}

$(3)$ For every smooth $\varrho$ in the irreducible component of $\,\Rep_G[\pi_1(X)]\,$ of the trivial representation,
$$
\HR_i(X,G)_{\varrho} =0\ \, \mbox{for all}\ \, i>\dim(G)\cdot n_L\ ,
$$
where $ n_L $ is the number of components of the link $L$.
\ethm
We note that the maximal dimension $M$  featuring in Part $(1)$ of
Theorem~\ref{T3} is a familiar geometric invariant which can be
computed explicitly for some simple links (for example, for $(p,q)$-torus knots in $ \R^3 $,  we have $\, M \leq 2 \dim G \,$, see Corollary~\ref{pqtorus}). However, we do not know any general formula for this invariant.

We conclude the introduction by stating one {\it non}-vanishing result. Let $ p,\, q_1, \,q_2,\, \ldots,\, q_m $ be integers  such that $ p>1$ and $\,q_1,\, q_2,\, \ldots,\, q_m\,$ are coprime to $p$. Let $ L_p(q_1, q_2, \ldots, q_m) $ denote the $(2m-1)$-dimensional lens space of type $(p; q_1, q_2, \ldots, q_m)$. Recall that
$ L_p(q_1, q_2, \ldots, q_m) $ is defined as the orbit space $ {\mathbb S}^{2m-1}\!/\Z_p $, where $ \Z_p $ acts on $ {\mathbb S}^{2m-1} \subset \c^m $ by  complex rotations: $\, (z_1, z_2, \ldots, z_m) \mapsto (\lambda_1 z_1, \lambda_2 z_2, \ldots, \lambda_m z_m)\,$ with $ \lambda_j = e^{2 \pi i q_j/p}\,$ ($ j=1,2,\ldots,m$). We have

\bthm[see Theorem \ref{rephomlens}]
\la{T4}
Let $ G $ be $\, \GL_n(\c) \,$ or $\, \SL_n(\c) \,$ for $\, n \ge 1 \,$. Then
$$
\HR_i( L_p(q_1, q_2, \ldots, q_m),\,G) \not= 0
\ \  \mbox{for}\ \ i \not= 2(m-1)k \ ,\ k = 1,\,2,\,3,\,\ldots\  ,
$$
and $\, \HR_i( L_p(q_1, q_2, \ldots, q_m),\,G) = 0\, $  otherwise.
\ethm
The result of Theorem~\ref{T4} should be compared to that of Theorem~\ref{T1}: since $ \pi_1(L_p) \cong \Z_p $, the latter theorem implies $\, \HR_i(\pi_1(L_p),\,G) = 0 $ for all $\, i > 0 $.

\vspace*{1ex}

The paper is organized as follows. In Section \ref{S2}, we introduce notation and review some
basic facts about simplicial sets, simplicial groups and simplicial commutative algebras.
This section contains no new results, except possibly for Theorem~\ref{bounded_homology},
which we could not find in the literature in this form and generality ({\it cf.}, however,
Remark after this theorem). In Section~\ref{S3}, we compute the homology of the derived
cotangent complex $ \T_\varrho \DRep_G(X) \,$ (Theorem~\ref{LinRep}) and its Euler
characteristic (Proposition~\ref{Euler}). We also give a natural interpretation of
$ \T_\varrho\DRep_G(X) $  in terms of functor homology (see Section~\ref{S3.2}). One interesting consequence of this interpretation is a spectral sequence relating the cotangent homology of $ \DRep_G(X) $ to the Pontryagin ring of $X$ in case when the space $X$ is simply connected 
(Corollary~\ref{corfhint1}). 
In Section~\ref{S4}, we introduce the derived character scheme $ \DChar_G(X) $ for
a reductive group $G$ and prove our first main result (Theorem~\ref{T0}).
In Section~\ref{S5}, we prove our vanishing theorems for virtually free groups,
closed surfaces, link complements and lens spaces. An interesting consequence of
these theorems is the existence of a well-defined $K$-theoretic virtual fundamental class $\,[\DRep_G(X)]_K^{\rm vir} \in K_0(\Rep_G[\pi_1(X)])\,$  in the sense of \cite{CK}. We give a new ``Tor formula'' for this class in terms of
functor homology (Corollary~\ref{vfclass0}) and show, in particular, that $\, [\DRep_G(X)]_K^{\rm vir} \not=0
\,$ when $ X = \Sigma_g $ is a surface of genus $ g \ge 2 $ and $ G $
semisimple (Corollary~\ref{vfclass}). 

The paper ends with an Appendix, where we clarify the relation of our work to some earlier work in derived algebraic geometry. Specifically, following a suggestion of the referee, we compare our construction of $ \DRep_G(X) $ to
the To\"{e}n-Vezzosi construction of the derived mapping stack $ \textbf{Map}(X, BG) $ of flat $G$-bundles on $X$ (see \cite{PTVV13}) as well as Kapranov's original construction of the derived moduli space of $G$-local systems,  $ \boldsymbol{R} \Loc_G(X) $. Strictly speaking, the material discussed in the Appendix is not used and not needed for understanding the main results of the paper but it puts these results in a proper geometric context
and provides a dictionary that allows one to translate our theorems into the language of derived algebraic geometry.

\vspace{1ex}

%
%

%

\subsection*{Acknowledgments}
We thank Damien Calaque, Giovanni Felder, Michael Larsen and Tony Pantev for interesting discussions, questions and suggestions. We also thank the anonymous referee for providing insightful comments on the earlier version of this paper: in particular, for explaining to us the connection
between the derived representation scheme $ \DRep_G(X) $ and the To\"{e}n-Vezzosi derived mapping stack of flat $G$-bundles on $X$, which we discuss in the Appendix.
The second author is grateful to Forschungsinstitut f\"{u}r Mathematik (ETH, Z\"{u}rich) for 
excellent working conditions during his stay through the fall of 2017. He also thanks the Department of Mathematics, Indian Institute of Science, Bengaluru, for its hospitality during his visit in the summer of 2017. Research of the first two authors was partially supported by the Simons Foundation Collaboration Grant 066274-00002B as well as NSF grants DMS 1702323 and DMS 1702372.

\section{Preliminaries}
\la{S2}
In this section, we fix notation and review some basic facts about simplicial sets, simplicial groups and simplicial commutative algebras.

\subsection{Simplicial sets}
\la{S2.1}
Let $\Delta$ denote the simplicial category whose objects are the finite ordered sets $[n]\,=\,\{0<1< \ldots <n\}$,
$n\geq 0$, and morphisms are the (weakly) order preserving maps.
A contravariant functor from $\Delta$ to a given category $\mathscr{C}$
is called a {\it simplicial object} in $\mathscr{C}$. The simplicial objects in $\mathscr{C}$ form a category which we denote by
$s\mathscr{C}$. In particular, we write $\sset$ for the category of simplicial sets, and denote by $\sset_0$ its full subcategory consisting of reduced simplicial sets. We recall that a simplicial set $X$ is {\it reduced} if $ X_0 := X([0]) $ is a singleton.

The geometric realization of simplicial sets defines a functor $\, |\,\mbox{--}\,|:\, \sset \to \Top \,$ which provides the category $ \sset $
with the notion of weak equivalence (to wit, a map $ X \to Y $ of simplicial sets is called a {\it weak equivalence} in $ \sset $ if the corresponding map of spaces $|X| \to |Y|$ is a weak homotopy equivalence in the usual
topological sense). By inverting weak equivalences, one defines the homotopy category of simplicial sets, $ \Ho(\sset) $,
and a classical result of simplicial homotopy theory asserts that $\, |\,\mbox{--}\,|\,$ induces an equivalence of
categories  $\,\Ho(\sset)\cong \Ho(\Top) $  (see, e.g., \cite[Chap.~III, \S 16]{M}). When restricted to reduced simplicial sets this equivalence becomes
\begin{equation}
\la{settop}
\Ho(\sset_0)\cong \Ho(\Top_{0,*})
\end{equation}
where  $ \Top_{0,*} $  is the full subcategory of $\Top $ comprising the pointed connected spaces.

Next, we recall that if $ \mathscr{C} $ is an abelian category (for example, the category $ \Vect_k $ of
$k$-vector spaces), the category $ s\mathscr{C}$ is equivalent
to the category $ \mathtt{Ch}_{\ge 0}(\mathscr{C}) $ of non-negatively graded chain complexes over $ \mathscr{C} $.
This classical equivalence (called the Dold-Kan correspondence) is given by a functor $\,N: s\mathscr{C} \to \mathtt{Ch}_{\ge 0}(\mathscr{C}) \,$, which is usually called the {\it normalization functor} (see \cite[8.4]{W}).

Throughout this paper, we will adopt the following standard notational convention.

\vspace{1ex}

\noindent
\textit{Convention.} A (reduced) simplicial set $X$ and its geometric realization $|X|$ will be denoted by the same
symbol (i.e., $|X| \equiv X $). We shall also not distinguish notationally between a simplicial object $V$ in an
abelian category and its normalized chain complex $ N(V) $ (i.e., $ N(V) \equiv V $). In particular, if
$V \in \mathtt{s}\Vect_k $ is a simplicial $k$-vector space, we denote the homology of the chain complex $N(V)$
by $\H_{\ast}(V)$.

\vspace{1ex}

The category $\sset $ of simplicial sets has a natural model structure, with weak equivalences being the
weak homotopy equivalences, the cofibrations the degreewise injective maps, and the fibrations the Kan
fibrations of simplicial sets \cite{Q1}. If $ s\mathscr{C} $ is a simplicial category, where
the objects have underlying simplicial sets (i.e., there is a forgetful functor $ U:\,s\mathscr{C} \to \sset $),
then $ s\mathscr{C} $ can be equipped with an induced model structure by declaring a morphism $ f $ in $ s\mathscr{C} $
to be a weak equivalence (resp., fibration) if $ U(f) $ is a weak equivalence (resp., fibration) in $ \sset $
(see \cite[Part~II, \S 4]{Q1}). In this way, the categories of simplicial objects in basic algebraic categories (such as
groups, modules, associative algebras, Lie algebras, commutative algebras, etc.) acquire natural model structures,
which allows one to do ``homotopy theory" in these categories. In the present paper, we will be mostly concerned with
the category of simplicial groups,  $ \sGr $, and the category of simplicial commutative algebras, $ \scAlg_k $.
In the next two sections we will describe basic features of these model categories.

\subsection{Simplicial groups}
\la{simpgr}
Unlike $\sset$, the model category $\sGr$ is fibrant: by a classical theorem of Moore, every simplicial group is a Kan complex (see \cite[Theorem~17.1]{M}). The cofibrant objects in $ \sGr $ are semi-free simplicial groups and their retracts. We recall
that a simplicial group $ \Gamma \in \sGr $ is called {\it semi-free} if all its components $ \Gamma_n $, $ n \ge 0 $,
are freely generated by some subsets $ B_n \subset \Gamma_n $, and the set of all these generators $\{B_n\}_{n \geq 0}$ is closed under the degeneracy maps of $\Gamma $. The elements of $ B_n $ which are not the images of degeneracies
$ s_i: B_{n-1} \to B_n $ are called nondegenerate generators of $ \Gamma $; we denote these elements by $ \bar{B}_n $.

Semi-free simplicial grous arise naturally from reduced simplicial sets via the Kan construction \cite{Kan1}. Recall
that this construction gives a pair of adjoint functors
\begin{equation}
\la{kanl}
\lgr :\, \sset_0\, \rightleftarrows \,\sGr\,:\overline{W}
\end{equation}
called the {\it Kan loop group functor} and the {\it classifying space functor}, respectively.
For an explicit definition and basic properties of these functors we refer to \cite[Chap.V]{GJ}. Here, we only
recall that the pair \eqref{kanl} is a Quillen equivalence: the functor $ \lgr $ preserves weak equivalences
and cofibrations, the functor $\overline{W}$ preserves weak equivalences and fibrations, and after inverting the weak equivalences, $ \lgr $ and $ \overline{W} $ become inverse to each other (see \cite[Corollary~V.6.4]{GJ}). Thus, in combination with \eqref{settop} we have equivalences of homotopy categories
\[
\,\Ho(\mathtt{Top}_{0,\ast}) \,\cong\, \Ho(\sset_0)\,\cong\, \Ho(\sGr)\, .\]
A simplicial group whose homotopy type coincides (under the above equivalence) with the homotopy type of
a space $X$ is called a {\it simplicial group model} of $X$.

The topological meaning of the loop group construction is clarified by the fundamental theorem of Kan which asserts
that for any reduced simplicial set $X$, there is a (weak) homotopy equivalence
%
$$
|\lgr{X}| \simeq \Omega|X|
$$
where $\Omega|X|$ is the $($Moore$)$ based loop space of $|X|$ (see \cite[Chap. V, Corollary~5.11]{GJ}).
In particular, we have $\pi_n(X) \,\cong\,\pi_{n-1}(\lgr{X})$ for all $n \geq 1$.

Now, for a {\it reduced} CW-complex $X$, there is another construction (again due to Kan \cite{Kan})
that provides a `small' semi-free simplicial group model $\mathbf{\Gamma}X$ of $X$,  such that for all $n \geq 0$,
the set $ \bar{B}_n $
of nondegenerate generators of $\mathbf{\Gamma}X$ in degree $ n $ coincides with the set of $(n+1)$-cells of $X$. The simplicial group $\mathbf{\Gamma}X$ can be constructed inductively as follows: put
a (partial) order on the set $ \{\sigma\} $ of cells of $X$ so that $ \sigma < \sigma' $ whenever
$ \dim \sigma < \dim \sigma' $ in $X$. Let $X_{< \sigma}$ (resp., $X_{\leq \sigma}$) denote the CW-subcomplex of $X$ formed by the cells less than $\sigma$ (resp., at most equal to $\sigma$) under the chosen order. Given an $(n+1)$-cell $\sigma$ of $X$ and given $\mathbf{\Gamma}X_{<\sigma}$, let $x \,\in\,(\mathbf{\Gamma}X_{<\sigma})_{n-1}$ denote any representative of $\lambda\,\in\,\pi_{n}(X_{<\sigma})\,\cong\,\pi_{n-1}(\mathbf{\Gamma}X_{<\sigma})$, where $\lambda$ denotes the homotopy class of the attaching map of $\sigma$. Then, a basis of $\mathbf{\Gamma}(X_{\leq \sigma})$ is given by adjoining $[\sigma]$ (in degree $n$) and all its degeneracies to a basis of $\mathbf{\Gamma}X_{<\sigma}$. The face maps on the extra generators are determined by
$$d_0([\sigma])=x\,,\,\,\,\,d_i([\sigma])=\id\,,\,\,1 \leq i\leq n \ . $$
This construction shows, in particular, that any reduced CW-complex with finitely many cells in each dimension has a semi-free simplicial group model with finitely many generators in each simplicial degree.

Next, following \cite{Q1}, we describe the homology of a local system on $X$ in terms of its simplicial group model. Recall  that if $R$ is a simplcial ring, the categories of left and right simplicial $R$-modules (to be denoted by  $R\text{-}\sMod$ and $\sMod\text{-}R$ resectively) carry standard model structures inherited from $\sset$. A morphism of simplicial modules $f\,:\,M \rar N$ is a free extension if there is a subset $C_n \subset N_n$ for each $n$ with $C:= \cup_n C_n$ is closed under degeneracies such that $f_n$ induces an isomorphism $M_n \oplus C_n.R_n \cong N_n$ for each $n$. A morphism $f\,:\,M \rar N$ in
$\sMod\text{-}R$ is a cofibration iff it is a retract of a free extension.

There is a natural bifunctor
\[ \mbox{--} \otimes_R \mbox{--}\,:\,  \sMod\text{-}R \times R\text{-}\sMod\, \rar\, \Z\mbox{-}\sMod\,\]
defined by $(M \otimes_R N)_n:= M_n \otimes_{R_n} N_n$, $n \geq 0$. This bifunctor has a total left derived functor $\otimes^{\bm L}: \Ho(\sMod\text{-}R) \times \Ho(R\text{-}\sMod) \rar \Ho(\Z\mbox{-}\sMod)$, which can be computed using resolutions: explicitly, for $M \in \sMod\text{-}R$ and $N \in R\text{-}\sMod$, $M \otimes^{\bm L} N \cong P \otimes N \cong N \otimes Q$ in $\Ho(\sMod_{\Z})$ for any cofibrant resolution $P \stackrel{\sim}{\rar} M$ in $\sMod\text{-}R$ and for any cofibrant resolution $Q \stackrel{\sim}{\rar} N$ in $R\text{-}\sMod$. Clearly, in this construction we can replace $\Z$ by any field $k$ and a simplicial ring $R$ by a simplicial $k$-algebra.

Let $X$ be a space and let $M$ be a discrete simplicial left $\Z[\lgr{X}]$-module. Giving such $M$ is equivalent to giving a representation of $\pi_0(\lgr{X})\cong \pi_1(X)$ in abelian groups, which corresponds to a local system on $X$. Viewing
$\Z$ as a constant simplicial right module with trivial $\lgr{X}$-action we can form the simplicial abelian group
$\,\Z \otimes^{\bm L}_{\Z[\lgr{X}]} M \in  \Ho(\Z\mbox{-}\sMod)$. Then, we have
\blemma
\la{localsys}
There is an isomorphism of graded abelian groups
$$
\pi_{\ast}(\Z \otimes^{\bm L}_{\Z[\lgr{X}]} M) \,\cong\, \H_{\ast}(X,M)
$$
where $\H_\ast(X,M)$ is the homology of $X$ with coefficients in the local system corresponding to $M$.
\elemma

\subsection{Simplicial commutative algebras}
\la{S2.3}
Recall that $\scAlg_k$ denotes the model category of simplicial commutative algebras over a field $k$.
When $k$ has characteristic zero (as we always assume in this paper), this model category is known to
be Quillen equivalent to the model category of non-negatively graded commutative DG $k$-algebras,
$\cDGA_k^+ $, where the weak equivalences are the DG algebra maps inducing isomorphisms on homology.
Specifically, by \cite[Part I, Sect. 4]{Q}, the Quillen equivalence between these
model categories is given by
\begin{equation} \la{NNstar_equiv} N^{\ast} \,:\, \cDGA^+_k \rightleftarrows \scAlg_k\,:\,N \,,\end{equation}
where $N$ is the Dold-Kan normalization functor and $N^{\ast}$ its left adjoint taking semi-free commutative DG algebras to semi-free simplicial commutative algebras ({\it cf.} \cite[Prop. A.2]{BRY}). Using \eqref{NNstar_equiv}, one can prove (see \cite[Prop.~3.1.5]{Lur04} or \cite[Section~11.1]{WY} for details):
\bprop  
\label{finite_type_equiv}
For a simplicial commutative $k$-algebra $A$, the following conditions are equivalent:
\begin{enumerate}
\item[$(1)$] $\pi_0(A)$ is a finitely generated $k$-algebra, and each $\pi_n(A)$ is a finitely generated module over $\pi_0(A)$.
\item[$(2)$] The homotopy type of $N(A) \in \cDGA_k^+$ has a representative
$B$ which is a semi-free DG algebra with finitely many generators in each homological degree.
\item[$(3)$] The homotopy type of $A$ has a representative in $\sCommAlg_k$ such that each $A_n$ is a finitely generated algebra over $k$.
\end{enumerate}

We say that $A$ is {\rm of quasi-finite type} $($or {\rm almost of finite presentation$)$ over} $k$ if it satisfies
one of the $($equivalent$)$ conditions $(1)$-$(3)$.
\eprop

Next, we recall that for an algebra $A \,\in\,\sCommAlg_k$, we can define an {\it augmentation} over $k$ in three equivalent
ways: either a morphism $A \rar k$ in $\sCommAlg_k$ or a morphism $A \rar k$ in $\Ho(\sCommAlg_k)$, or a morphism
$\pi_0(A) \rar k$ in $\cAlg_k$. We denote such an augmentation by $\varepsilon$, thinking of it as a $k$-point of the affine scheme $\Spec[\pi_0(A)]$.

Choosing an augmentation of $A$ allows one to ``linearize" $A$. Namely, given $\varepsilon : A \rar k$,
consider the simplicial $k$-module $\cQ(A) := \mathfrak{m}_{\varepsilon}/\mathfrak{m}_{\varepsilon}^2 $ where $\mathfrak{m}_{\varepsilon} :=\Ker(\varepsilon)$. This defines a functor on the category of augmented simplicial commutative algebras:
\begin{equation}
\la{linfun}
\cQ:\ \sCommAlg_{k/k} \rar \sMod_k\ ,
\end{equation}
called the {\it linearization $($or abelianization$)$ functor}.

The functor $\cQ$ is left Quillen: it has a right adjoint  given by $ M \mapsto k \ltimes M $, where
$k \ltimes M $ is the semidirect product (square-zero extension) of the discrete simplicial algebra $k$
with a simplicial $k$-module $M$. It follows that the linearization functor has a total left derived functor
\begin{equation}
\la{dlinfun}
 \L \cQ : \ \Ho(\sCommAlg_{k/k}) \rar \Ho(\sMod_k)\ ,
\end{equation}
which assigns to $ (A,\varepsilon) $ the linearization $ \cQ(P) = \Omega^1(P) \otimes_P k $ of a cofibrant replacement
$ P $ of $A$ in $ \sCommAlg_k $ with respect to the augmentation $ P \sonto A \xrightarrow{\varepsilon} k $.
The object $ {\bm L} \cQ(A, \varepsilon) $ is isomorphic to $\,{\bm L}_{A|k} \otimes_A k\,$ in $\,\Ho(\sMod_k)\,$, 
where $\,{\bm L}_{A|k} :=  \Omega^1(P) \otimes_P A \,$ is the {\it classical cotangent complex} of $A$ relative
to $k$ in the sense of \cite{A74, Q70}. When the augmentation on $A$ is fixed, we will often omit it
from the notation, writing $ {\bm L} \cQ(A) $  instead of $ {\bm L} \cQ(A, \varepsilon) $. 

The next proposition is a standard homological result, well known to experts; as we could not find a suitable reference, we briefly outline its proof for the benefit of the reader.
\bprop \la{aqhomspecseq}
For any augmented simplicial commutative $k$-algebra $A$, there is a convergent homology spectral sequence
$$
{\it E}^1_{p,q} = \H_q[\L\cQ(A_p)] \,\implies\, \H_{p+q}[\L\cQ(A)]\ .
$$
\eprop
\bproof
The linearization functor is also defined for augmented DG commutative algebras. This gives
$$\cQ\,:\,\cDGA_{k/k}^+ \rar \Com_k\,,\,\, B \mapsto \mathfrak{m}/\mathfrak{m}^2\,, $$
where $\mathfrak{m}$ denotes the DG augmentation ideal of $B$. There is an isomorphism of functors from $\cDGA_{k/k}^+$ to $\sMod_k$:
$$\cQ \circ N^{\ast} \,\cong\, N^{-1} \circ \cQ \ .$$
Indeed, the right adjoints of the functors on both sides are easily seen to be isomorphic. Hence the functors themselves are isomorphic.

Now, if $B \stackrel{\sim}{\rar} N(A)$ is any semi-free resolution, $N^{\ast}(B)$ gives a semi-free simplicial resolution of $A$. It follows that the homology of $\L\cQ(A)$ is isomorphic to the homology of $\cQ(B)$. Taking $B$ to be the canonical resolution\footnote{Namely, the ``cobar-bar" resolution obtained by the Quillen adjunction between the category of (conilpotent) DG Lie coalgebras and the category of commutative DG algebras (see, e.g. \cite[Theorem 6.1]{BFPRW} and the subsequent remark therein).} of $N(A)$, we see that $\H_{\ast}[\L\cQ(A)]$ is isomorphic to the homology of the shifted (reduced) Harrison chain complex $\mathrm{CHarr}_{\ast}(N(A),k)[-1]$ (see \cite[Sect. 4.2.10]{L}). It is easy to see that this complex is quasi-isomorphic to the total complex of the (first quadrant) double complex associated to the simplicial complex which assigns $\mathrm{CHarr}_{\ast}(A_n,k)[-1]$ to the simplex $[n]$. The spectral sequence associated with the filtration of the above double complex by rows is precisely
the desired spectral sequence.
\eproof
Let $A$ be a simplicial commutative $k$-algebra of quasi-finite type. By Proposition \ref{finite_type_equiv}, the homotopy type of $N(A)$ has a representative $B\,\in\,\cDGA_k^+$ that is semi-free and has finitely many generators in each degree. It follows that the homotopy type $A$ has the representative $N^\ast(B)\,\in\,\scAlg_k$ that is semi-free and has finitely many generators in each simplicial degree. Hence, for every $\varepsilon\,:\,A \rar {k}$, the vector spaces $\H_i[{\bm L}{\cQ}(A,\varepsilon)]$ are finite-dimensional over ${k}$. Moreover, $\dim_{{k}}\,\H_i[{\bm L}{\cQ}(A)]$ does not exceed the number of generators of $N^\ast(B)$ in simplicial degree $i$ for any $i \geq 0$.

Fix an augmentation $\varepsilon: A \rar k$. For $q \geq 0$, let $\hat{\pi}_q(A)$ denote the completion of $\pi_q(A)$ with respect to the $\mathfrak{m}$-adic filtration given by powers of the augmentation ideal of $\pi_0(A)$ defined by $\varepsilon$.
\bprop
\la{Quillen_specseq}
Let $A$ be a simplicial commutative $k$-algebra of quasi-finite type. There is a natural spectral sequence concentrated in the region $\,q \geq 0,\, p \geq -q\,$:
\begin{equation}
\la{Qsp}
{\it E}^2_{p,q}\,:=\,\H_{p+q}[\Sym^q({\bm L}{\cQ}(A))]\,\implies\,\hat{\pi}_{p+q}(A)
\end{equation}
converging to the completed homotopy groups of $A$.
\eprop
\bproof
Pick a representative $B\,\in\,\cDGA_k^+$ of the homotopy type $N(A)$ that is semi-free and that has finitely many generators in each degree (such a choice exists by Proposition \ref{finite_type_equiv}). The augmentation $\varepsilon$ on $A$ defines an augmentation $\varepsilon\,:\,B \rar k$. Let $\mathfrak{m}:=\ker(\varepsilon)\,\subset\,B$. Let $\widehat{B}$ (resp., $\widehat{\mathfrak{m}}$) denote the completion of $B$ (resp., $\mathfrak{m}$) with respect to the $\mathfrak{m}$-adic filtration. Consider the increasing filtration $F_{\bullet}\widehat{B}$, where $F_{-p}\widehat{B}\,=\,\widehat{\mathfrak{m}}^p$ for $p \geq 0$. Note that $\widehat{B}$ is complete with respect to this exhaustive filtration. Since $\widehat{\mathfrak{m}}^p/\widehat{\mathfrak{m}}^{p+1}\,\cong\,\mathfrak{m}^p/\mathfrak{m}^{p+1}\,\cong\, \Sym^p(\mathfrak{m}/\mathfrak{m}^2)$ for all $p$,  the spectral sequence (concentrated in the region $\{(-p,q)\,|\,p \geq 0, q \geq p\}$) associated with this filtration is
\begin{equation}
\la{specseqfil}
 {\it E}^1_{-p,q}\,=\,\H_{q-p}[\Sym^p(\mathfrak{m}/\mathfrak{m}^2)]\, {\implies} \,\pi_{q-p}[\widehat{B}] \ .
\end{equation}
Since $B$ has finitely many generators in each degree, the spectral sequence \eqref{specseqfil} is regular and therefore, converges by the complete convergence theorem (see \cite[Theorem 5.5.10]{W}). A standard d\'{e}calage argument (see \cite[Ex. 5.4.3]{W}) then yields the convergent spectral sequence
\begin{equation} \la{decalage}
 {\it E}^2_{p,q}\,=\,\H_{p+q}[\Sym^q(\mathfrak{m}/\mathfrak{m}^2)]\,{\implies} \,\pi_{p+q}[\widehat{B}]
\end{equation}
associated to the shifted filtration $\tilde{F}_r\widehat{B}_n\,:=\,F_{r-n}\widehat{B}_n$.

Since $\cQ \circ N^{\ast}\,\cong\,N^{-1} \circ \cQ$, we have $\mathfrak{m}/\mathfrak{m}^2\,\cong\,\bm L \cQ(A)$. It therefore, suffices to verify that $\pi{\ast}[\widehat{B}]\,\cong\,\hat{\pi}_{\ast}(A)$. Let $\mathfrak{n} \subset B_0$ denote the kernel of the augmentation on $B_0$. Since the $\mathfrak{n}$-adic and $\mathfrak{m}$-adic filtrations on $B$ define the same topology, since each $B_n$ is finitely generated as a $B_0$-module and since $B_0$ is Noetherian $\widehat{B}\,\cong\,B \otimes_{B_0} \widehat{B}_0$, where $\widehat{B}_0$ denotes the completion of $B_0$ with respect to the $\mathfrak{n}$-adic filtration. Since the functor $(\,\mbox{--}\,) \otimes_{B_0} \widehat{B}_0$ is exact, we have
$$ \pi_{\ast}[\widehat{B}]\,\cong\, \pi_{\ast}(B) \otimes_{B_0} \widehat{B}_0\,\cong \,\pi_{\ast}(B) \otimes_{\pi_0(B)} \widehat{\pi_0(B)}\,,$$
where $\widehat{\pi_0(B)}$ is the completion of $\pi_0(B)$ with respect to the adic filtration defined by the augmentation ideal (the augmentation being induced by $\varepsilon$). The desired proposition now follows from the fact that $\pi_{\ast}(B)\,\cong\,\pi_{\ast}(A)$.
\eproof

\begin{remark} The spectral sequence of Proposition~\ref{Quillen_specseq} may be viewed as `desuspension' of Quillen's fundamental spectral sequence for the augmentation map  $ \varepsilon: A \rar k $. Indeed, when applied to the suspension $ \Sigma A $ of an augmented simplicial  $k$-algebra $A$, the spectral sequence \eqref{Qsp} becomes 
({\it cf.} \cite[Theorem 6.3]{Q70})
\begin{equation}
\la{ssQuillen} E^2_{p,q}=\H_{p+q}[\Sym^q({\bm L}_{k|A})]\,\implies\,\Tor^{A}_{p+q}(k,k) \, , 
\end{equation}
where $ {\bm L}_{k|A}$ is the classical cotangent complex associated to the morphism $ \varepsilon: A \rar k$. 
To see this, we first observe that, by \cite[I.4.4, Prop.~2]{Q1}, there is a natural isomorphism in $ \Ho(\sMod_k) $
$ $
$$ 
\L \cQ(\Sigma A)\,\cong\,  \L \cQ(A)[1]\ .
$$
Now, since $\pi_0(\Sigma A)\,\cong\,k$, the completion of $ \pi_\ast(\Sigma A)  $ in the
limit term of the spectral sequence \eqref{Qsp} for $ \Sigma A $ is trivial, and the limit itself is 
therefore given by $ \pi_\ast(\Sigma A) \cong \Tor^A_\ast(k,k) $ (see, e.g., \cite[Sect. 11.3]{DS}).

We note that when $k$ has characteristic zero, the Quillen spectral sequence \eqref{ssQuillen} degenerates for
any augmented $k$-algebra $A$ (see \cite[Theorem~7.3]{Q70}). By constrast, the spectral sequence \eqref{Qsp} of Proposition \ref{Quillen_specseq} is far from being degenerate in general, even when ${\rm char}(k) = 0 $.
\end{remark}
\vspace*{1ex}

The following theorem provides sufficient conditions for the vanishing of higher homotopy groups of a simplicial commutative algebra in terms of its cotangent homology. We will use this result repeatedly in our calculations in Section~\ref{S5}.
\begin{theorem} 
\la{bounded_homology}
Let $A$ be a simplicial commutative algebra defined over an algebraically closed field $k$. Assume that
$A$ is of quasi-finite type over $k$ in the sense of Proposition~\ref{finite_type_equiv}. 

$(a)$ 
If there is an integer $ N \ge 0 $ such that, for every augmentation $ \varepsilon:\, A  \rar {k} $,
$$ 
\dim_{{k}} \H_1[{\bm L}{\cQ}(A)] \leq N \quad \mbox{and} \quad 
\H_i[\bm L\cQ(A)] = 0\ ,\ \forall\,i > 1 \ ,
$$
then $\,\pi_n(A)=0\,$ for $\,n > N\,$. 

$(b)$ If for every augmentation $ \varepsilon:\, A  \rar {k} $, 
$\,\H_i[\bm L\cQ(A)] = 0\,$ for all $ i > 0 $, 
then $\,\pi_n(A)=0\,$ for  $\,n > 0\,$. Moreover, $\,{\rm Spec}[\pi_0(A)] $ is a smooth $k$-scheme in this case.
\end{theorem}
\bproof
$(a)$
Let $\mathfrak{m}$ be any maximal ideal of $\pi_0(A)$. Since $k$ is algebraically closed, $\mathfrak{m}$ defines an augmentation $\varepsilon\,:\,A \rar k$. Since $\H_i[\bm L\cQ(A)]=0$ for $i>1$ and $\dim_{{k}} \H_1[{\bm L}{\cQ}(A)] \leq N$,
$\H_n[\Sym({\bm L}{\cQ}(A))]$ vanishes for $n>N$. It follows from Proposition \ref{Quillen_specseq} that $\hat{\pi}_n(A)=0$ for $n > N$.
Since $A$ is of quasi-finite type, $\pi_0(A)$ is a Noetherian $k$-algebra and $\pi_i(A)$ is a finitely generated $\pi_0(A)$ module for $i \geq 0$. It follows that $\hat{\pi}_n(A)\,\cong\, \widehat{\pi_n(A)_{\mathfrak{m}}}$ for all $n$, where $\widehat{(\mbox{--})}$ denotes completion with respect to the $\mathfrak{m}_{\mathfrak{m}}$-adic filtration. Since $\hat{\pi}_n(A)=0$ for $n>N$, $\pi_n(A)_{\mathfrak{m}}=0$ for $n >N$ by Krull's Theorem. This proves part $(a)$.

$(b)$  The first claim of part $(b)$ follows immediately from part $(a)$ (take $ N = 0 $). Thus, if the condition of $ (b) $ holds, then $A$ is weakly equivalent to  $\pi_0(A)$ in $ \scAlg_k $. In that case,
for any augmentation $\varepsilon\,:\,A \rar k$, there are natural isomorphisms 
\begin{equation} 
\la{harrison} 
\H^1[\bm L\cQ(A)^\ast] \,\cong\, D^1(\pi_0(A)|k,k) \,\cong\,\mathrm{Harr}^2(\pi_0(A),k)\,,\end{equation}
where $ \bm L\cQ(A)^\ast $ denotes the graded $k$-linear dual of ${\bm L}{\cQ}(A)$, $ D^1$ and $\mathrm{Harr}^2$  
are the 1-st Andr\'{e}-Quillen and the 2-nd Harrison cohomology of $ \pi_0(A) $ with coefficients
in $k$ (viewed as a $ \pi_0(A)$-module via $\varepsilon$). The first isomorphism in 
\eqref{harrison} is induced by the canonical projection $ A \sonto \pi_0(A) $, while the second is the standard 
isomorphism relating Andr\'{e}-Quillen and Harrison cohomology (see, e.g., \cite[Corollary~8.8.9]{W}). 
Since $\,\H_1[{\bm L}{\cQ}(A)]=0\,$, we conclude from \eqref{harrison} that $\mathrm{Harr}^2(\pi_0(A),k)\,=\,0$. By  \cite[Corollary 20]{Har}, the morphism $\varepsilon$ corresponds then to a smooth $k$-point of 
$\mathrm{Spec}[\pi_0(A)]$. Since $\varepsilon$ is arbitrary, ${\rm Spec}[\pi_0(A)]$ is smooth.
\eproof

\begin{remark} 
Theorem~\ref{bounded_homology} can be found in the literature in various forms and generality. 
For example, the result of part $(a)$ appears in (the proof of) \cite[Theorem~2.2.2]{CK}. Part $(b)$ is 
closely related to \cite[Proposition~2.2.5.1 and Corollary~2.2.5.3]{TV08}, although these results are stated 
in \cite{TV08} in a different language and under different assumptions than our Theorem~\ref{bounded_homology}. In particular, under the assumptions of Theorem~\ref{bounded_homology}$(b)$, it follows 
that $ k \rar A $  is a strongly smooth morphism in the sense of \cite[Definition~2.2.2.3 and Theorem~2.2.2.6]{TV08}. 
We thank the referee for kindly pointing out to us the latter reference.
\end{remark}

\section{Linearization of representation homology}
\la{S3}
\subsection{Derived cotangent complex}
\la{S3.1}
Let $X$ be a space and let $G$ be an affine algebraic group over a field $k$ with Lie algebra $\g$.
As in the Introduction, we denote by $\mathcal{A}_G(X)$ the derived representation algebra of $X$. Recall
that $\mathcal{A}_G(X) = {\bm L}(\Gamma X)_G $ is an object of $ \Ho(\scAlg_k) $ defined by applying
the derived representation functor in $G$ to a simplicial group model of $X$.
Giving an augmentation of $\mathcal{A}_G(X)$ over $k$ is equivalent to giving an augmentation of $\HR_0(X,G)$, which, in turn, is equivalent to fixing a representation $\varrho\,:\,\pi_1(X) \rar G(k)$. Composing such a representation with the adjoint action of $G(k)$ in $ \g $, we get the action $\mathrm{Ad} \circ \varrho:\,\pi_1(X) \rar \GL_k(\g)$. In this way, $\g$ acquires the structure of a left $\pi_1(X)$-module. We denote the corresponding local system on $X$ by $\mathrm{Ad}\,\varrho$, and write  $\mathrm{Ad}^{\ast}\varrho$ for the dual local system with fiber $ \g^* $.

Note that $\mathcal{A}_G(X)$ given together with augmentation $\varrho$ is an object of $\Ho(\scAlg_{k/k})$. The derived linearization functor $\L\cQ$  applied to the pair $(\mathcal{A}_G(X), \varrho)$  gives a simplicial $k$-module
$\mathrm{T}^{\ast}_{\varrho}\DRep_G(X):=\L\cQ_{\rho}[\mathcal{A}_G(X)] \in \Ho(\sMod_k)$, which may be interpreted geometrically as the derived cotangent space to $\DRep_G(X)$ at $\varrho$. We have
\bthm
\la{LinRep}
There are natural isomorphisms of vector spaces
\begin{equation}
\la{LinRepGr}
\H_i[\mathrm{T}^{\ast}_{\varrho}\DRep_G(X)]\,\cong\, \left\{
                                                 \begin{array}{ll}
                                                    Z^1(\pi_1(X), \mathrm{Ad}\,\varrho)^{\ast} & \text{if }\ i=0\\*[1ex]

                                                    \H_{i+1}(X,\mathrm{Ad}^{\ast}\varrho) & \text{if }\ i>0 \ \\
                                                 \end{array}
                                                 \right.
\end{equation}
\ethm
%
%
%
%
To prove Theorem~\ref{LinRep} we need some notation and a few technical results. First, we recall
one simple fact about extensions of algebraic groups (see \cite[Lemma 2.1]{LuM}).

%
\begin{lemma}
\la{dualnos}
Let $ k \ltimes V $ be the square-zero extension of $k$ by a vector space $V$. Then $\,G(k \ltimes V) \,\cong\, G(k) \ltimes (\g \otimes V)$, where the left $G(k)$-module structure $\g \otimes V$ arises from the adjoint action of $G(k)$ on $\g$.
\elemma

Next, we consider the slice category $\mathtt{Gr}/\pi$ of $\mathtt{Gr}$ over a fixed group $\pi$. By definition, the objects of $\mathtt{Gr}/\pi$ are the group homomorphisms $ \Gamma \to \pi $ and
the morphisms are the obvious commutative triangles  over $ \pi $.
By the universal property of the representation functor, any representation $\varrho\,:\,\pi \rar G(k)$
defines an augmentation of the algebra $\pi_G$, which induces an augmentation on $ \Gamma_G $ for every
$\Gamma \in \mathtt{Gr}/\pi$. This last augmentation is given by the composite map $\Gamma_G \rar \pi_G \stackrel{\varrho}{\rar} k$, which corresponds to the pullback (via the structure map) of $\varrho$ to a representation of $\Gamma$. Thus, the representation functor $(\mbox{--})_G$ may be viewed as a functor
$(\mbox{--})_G\,:\,\mathtt{Gr}/\pi \rar \cAlg_{k/k}$, and the associated adjunction becomes
\begin{equation} \la{augadj} (\,\mbox{--}\,)_G :\,\mathtt{Gr}/\pi \rightleftarrows \cAlg_{k/k}\,: \pi \times_{G(k)} G(\,\mbox{--}\,) \ . \end{equation}
Now, given $\Gamma \,\in\,\mathtt{Gr}/\pi$, equip $\g$ with the left $\Gamma$-action via the composite map
$$
\mathrm{Ad}\,\varrho\,:\, \Gamma \rar \pi \xrightarrow{\varrho} G(k) \xrightarrow{\mathrm{Ad}} \GL_k(\g)\ .
$$
By duality, this defines a right $\Gamma$-action on $\g^{\ast}$ and hence makes $\g^{\ast}$ a right $\Z[\Gamma]$-module.
\bprop \la{linrepprop}
For $\Gamma \in \mathtt{Gr}/\pi$, there is a natural isomorphism of $k$-vector spaces
\[
\cQ(\Gamma_G) \,\cong\,  \g^{\ast} \otimes_{\Z[\Gamma]} \Omega(\Gamma)\,,
\]
where $\Omega(\Gamma)$ is the augmentation ideal of $\Z[\Gamma]$.
\eprop
%
%
\bproof
For any $V \in \ve_k$, there are natural isomorphisms
\begin{eqnarray*}
\Hom_k(\cQ(\Gamma_G), V) & \cong & \Hom_{\cAlg_{k/k}}(\Gamma_G, k \ltimes V)\\
                         & \cong & \Hom_{\mathtt{Gr}/\pi}(\Gamma, \pi \times_{G(k)} G(k \ltimes V))\,\,\,\,\,[\text{ by }\eqref{augadj}]\\
                         & \cong & \Hom_{\mathtt{Gr}/\pi}(\Gamma,  \pi \ltimes (\g \otimes V)) \,\,\,\,\,[\text{ by Lemma \ref{dualnos}}]\\
                         & \cong & \Der(\Gamma, \g \otimes V)\\
                         & \cong & \Hom_{\Z[\Gamma]}(\Omega(\Gamma), \g \otimes V)\\
                         & \cong & \Hom_k(\g^{\ast} \otimes_{\Z[\Gamma]} \Omega(\Gamma), V)\ .
\end{eqnarray*}
The desired proposition then follows from the Yoneda lemma.
\eproof
\bproof[Proof of Theorem \ref{LinRep}]
Recall (see \cite[Definition B.2]{BRY}) that a map $f\,:\,A \rar B$ of simplicial commutative algebras is a {\it smooth extension}\footnote{We warn the reader that our notion of smoothness differs from the standard geometric one used, e.g., in 
\cite{TV08}.} if $f$ can be written as an (infinite) composition
$$ A={sk}_{-1}(A) \rar {sk}_0(A) \rar {sk}_{1}(A) \rar \ldots \rar \varinjlim {sk}_{\ast}(A)=B$$
such that each map ${sk}_{n-1}(A) \rar {sk}_{n}(A)$ is a pushout of the form
$$
\begin{diagram}[small]
    \partial \Delta[n] \otimes C_n & \rTo & \Delta[n] \otimes C_n\\
      \dTo & & \dTo\\
       {sk}_{n-1}(A) & \rTo & {sk}_n(A)
       \end{diagram}
$$
for some smooth commutative $k$-algebra $C_n$. A simplicial commutative algebra $A$ is called {\it smooth} if the structure map $k \rar A$ is a smooth extension.  Since $\lgr{X}$ is a semi-free simplicial group and since $\O(G)$ is a smooth commutative algebra, the simplicial commutative algebra $(\lgr{X})_G$ is smooth. It follows from \cite[Prop. B.3]{BRY} that
the canonical map $\L \cQ_{\varrho}[\mathcal{A}_G(X)] \rar \cQ_{\varrho}[\mathcal{A}_G(X)]$ is an isomorphism in $\Ho(\sMod_k)$. In other words, $\L\cQ_{\varrho}[\mathcal{A}_G(X)]$ is represented in $\Ho(\sMod_k)$ by $\cQ_{\varrho}[(\lgr{X})_G]$. By Proposition \ref{linrepprop}, there is an isomorphism of simplicial $k$-vector spaces
\[
\cQ_{\varrho}(\lgr{X}_G) \,\cong\, \g^{\ast} \otimes_{\Z[\lgr{X}]} \Omega(\lgr{X}) \,,
\]
where the discrete simplicial $k$-vector space $\g^{\ast}$ acquires its right $\lgr{X}$-action via the homomorphisms
$\lgr{X} \rar \pi_1(X) \xrightarrow{\mathrm{Ad}} \GL_k(\g)$. Here, $\pi_1(X)$ and $\GL_k(\g)$ are viewed as constant simplicial groups. Since the categories of (simplicial) left and right $\Z[\lgr{X}]$-modules are equivalent, we have
\[
\cQ_{\varrho}(\lgr{X}_G) \,\cong\, \g^{\ast} \otimes_{\Z[\lgr{X}]} \Omega(\lgr{X}) \,\cong\,  \Omega(\lgr{X}) \otimes_{\Z[\lgr{X}]} \g^{\ast}\ .
\]
Since $\lgr{X}$ is semi-free, $\Omega(\lgr{X})$ is a free $\lgr{X}$-module at every level (see, e.g. \cite[Chap. IV.2, Ex. 3]{KB}). Hence, $\Omega(\lgr{X}) \otimes_{\Z[\lgr{X}]} \g^{\ast} \,\cong\, \Omega(\lgr{X}) \otimes^{\L}_{\Z[\lgr{X}]} \g^{\ast}$ in $\Ho(\sMod_k)$. It follows that in $\Ho(\sMod_k)$,
\begin{equation} \la{ethm1}
\mathrm{T}^{\ast}_{\varrho}\DRep_G(X) \,\cong\, \Omega(\lgr{X}) \otimes^{\L}_{\Z[\lgr{X}]} \g^{\ast}\ .
\end{equation}
The exact sequence of simplicial $\Z[\lgr{X}]$-modules
$$ 0 \rar \Omega(\lgr{X}) \rar \Z[\lgr{X}] \rar \Z \rar 0$$
 gives a distinguished triangle in the derived category $\mathscr{D}(k)$ of complexes of $k$-vector spaces
\[
\Omega(\lgr{X}) \otimes^{\L}_{\Z[\lgr{X}]} \g^{\ast} \rar \Z[\lgr{X}] \otimes^{\L}_{\Z[\lgr{X}]} \g^{\ast} \cong \g^{\ast} \rar \Z \otimes^{\L}_{\Z[\lgr{X}]} \g^{\ast} \rar \Omega(\lgr{X}) \otimes^{\L}_{\Z[\lgr{X}]} \g^{\ast}[1] \ .
\]
The associated long exact sequence together with \eqref{ethm1} and Lemma~\ref{localsys} give us the isomorphisms
\[\H_i[\T^{\ast}_{\varrho}\DRep_G(X)]\,\cong\,
\H_{i+1}(X,\mathrm{Ad}^{\ast}\varrho) \ , \quad i > 0\ ,
\]
as well as the exact sequence
\[
0 \rar \H_1(X, \mathrm{Ad}^{\ast}\varrho) \rar \H_0[\T^{\ast}_{\varrho}\DRep_G(X)] \rar \g^{\ast} \rar \H_0(X, \mathrm{Ad}^{\ast}\varrho) \rar 0 \ .
\]
This last exact sequence is dual to
\[
0 \rar \H^0(X,\mathrm{Ad}\,\varrho) \rar \g \rar \Der(\pi_1(X),\mathrm{Ad}\,\varrho) \rar \H^1(X,\mathrm{Ad}\,\varrho) \rar 0\ .
\]
Whence $ \H_0[\T^{\ast}_{\varrho}\DRep_G(X)] \cong \Der(\pi_1(X),\mathrm{Ad}\,\varrho)^* $,
which proves the desired result for  $ i = 0 $.
\eproof
The next proposition gives an explicit formula for the Euler characteristic of the cotangent complex of $\DRep_G(X)$ in terms of the topological Euler characteristic  $\chi_{\rm top}(X)$  of $X$.
\bprop \la{Euler}
Let $X$ be a finite reduced CW-complex. Then, for every $\varrho\,:\,\pi_1(X) \rar G(k)$, $\T^{\ast}_{\varrho}\DRep_G(X)$ is a bounded complex with finite-dimensional homology. Its Euler characteristic is given by
$$
\chi[\T^{\ast}_{\varrho}\DRep_G(X)]\,=\,[1-\chi_{\rm top}(X)]\, \dim G\,.  
$$
In particular, $\chi[\T^{\ast}_{\varrho}\DRep_G(X)]$ is independent of $\varrho$.
\eprop

\bproof
Recall (see Section \ref{simpgr}) that $X$ has a semi-free simplicial group model $\mathbf{\Gamma}X$ such that $n$-dimensional cells of $X$ correspond to non-degenerate generators of $\mathbf{\Gamma}X$ in degree $n-1$, for $n \geq 1$. The representation $\varrho$ extends to representations of the groups $(\mathbf{\Gamma}X)_n$ for all $n \geq 0$. This in turn, gives augmentations on the commutative algebras $\cO[\Rep_G((\mathbf{\Gamma}X)_n)]$ (each of which is isomorphic to a tensor product of copies of $\cO(G)$) for each $n$. Identifying the cotangent space of $\Rep_G[(\mathbf{\Gamma}X)_n]$ at the representation fixed by $\varrho$ to the cotangent space of $\Rep_G[(\mathbf{\Gamma}X)_n]$ at the trivial representation via right translation, we see that $\T^{\ast}_{\varrho}\DRep_G(X)$ is a simplicial vector space with nondegenerate basis elements in degree $n$ corresponding to basis elements of $\g^{\ast}$ labelled with each $(n+1)$-cell of $X$. Since $X$ is a finite cell complex, the first statement follows. Further, from this description, we see that
$$ \chi[\T^{\ast}_{\varrho}\DRep_G(X)]\,=\,\dim_k(\g^{\ast})(n_1-n_2+n_3- \ldots)\,,
$$
where $n_i$ is the number of $i$-cells in $X$. This proves the second assertion.
\eproof

\subsection{Functor homology interpretation}
\la{S3.2}
In \cite[Section~4.2]{BRY}, we gave a natural interpretation of representation homology  in terms of abelian homological algebra: specifically, as a derived functor tensor product over the (small) category of  f.g. free groups. In this section, we extend this interpretation to the cotangent homology of $ \DRep_G(X) $. We begin by reviewing the basic construction of \cite{BRY}.

Let $\mathfrak{G}$ denote the full subcategory of $\mathtt{Gr}$ whose objects are the free groups $\langle n \rangle$ based on the sets $\underline{n}:=\{x_1,x_2,\ldots,x_n\}$, $n \geq 0$ ($\langle 0 \rangle$ is the identity group by convention). The category $\mathfrak{G}\text{-}\Mod$ (resp., $\Mod\text{-}\mathfrak{G}$) of covariant (resp., contravariant) functors from $\mathfrak{G}$ to the category of $k$-vector spaces is an abelian category with sufficiently many injective and projective objects. We view the objects of $\mathfrak{G}\text{-}\Mod$ (resp., $\Mod\text{-}\mathfrak{G}$) as left (resp., right) $\mathfrak{G}$-modules. There is a natural bifunctor
\[
\mbox{--} \otimes_{\mathfrak{G}} \mbox{--} \,:\, \Mod\text{-}\mathfrak{G} \times \mathfrak{G}\text{-}\Mod \rar \ve_k
\]
which is right exact (with respect to each argument), preserves sums, and is left balanced (in the sense of \cite{CE}). The derived functors of $\mbox{--} \otimes_{\mathfrak{G}} \mbox{--}$ with respect to each argument are thus isomorphic, and we denote their common value by $\Tor^{\mathfrak{G}}_{\ast}(\mbox{--},\mbox{--})$.

It is known (see, e.g., \cite{Ha16}) that any commutative Hopf $k$-algebra $ {\mathcal H} $ defines a left $\mathfrak{G}$-module assigning $\, \langle n \rangle \mapsto {\mathcal H}^{\otimes n} $. In particular, if $G$ is an affine
algebraic group over $k$, the left $\mathfrak{G}$-module associated to $ \O(G) $ can be written in the form: $ \langle n \rangle \mapsto \O[\Rep_G(\langle n \rangle)]\,$, which makes the functoriality clear.
 Dually, a cocommutative Hopf algebra defines a right $\mathfrak{G}$-module. In particular, the simplicial cocommutative Hopf algebra $k[\lgr{X}]$ defines a simplicial right $\mathfrak{G}$-module, and hence a chain complex of right $\mathfrak{G}$-modules. It was proven in \cite{BRY} (see \cite[Corollary 4.3]{BRY}) that
\[
\HR_{\ast}(\Gamma,G) \,\cong\, \Tor^{\mathfrak{G}}_{\ast}(k[\Gamma],\O(G))\ .
\]
This formula is a consequence of a more general fact (see {\it loc. cit}, Theorem 4.2):
\[
\HR_\ast(X,G) \,\cong\, \H_\ast[k[\lgr{X}] \otimes^{\L}_{\mathfrak{G}} \O(G)] \ .
\]
The above isomorphisms provide a natural interpretation of representation homology in terms of classical abelian derived functors.

Now, let $\mathfrak{G}/\pi$ denote the slice category of $ \mathfrak{G} $ in $ \Gr $ over a fixed group $ \pi \in {\rm Ob}(\Gr) $. Thus, the objects in $\mathfrak{G}/\pi$ are the pairs $(\langle n\rangle, \phi)$, where $\langle n \rangle  $ is the free group based on the set $\underline{n}:=\{x_1,x_2,\ldots,x_n\}$ and $\phi:\,\langle n \rangle \rar \pi$ is a group homomorphism, and the morphisms $\,(\langle n\rangle, \phi) \to (\langle m\rangle, \psi)\,$ in $\mathfrak{G}/\pi$ are the commutative triangles in $ \Gr \,$:
$$
\begin{diagram}[small] 
\langle n \rangle & \rTo(3,0)^{f} & \langle m \rangle \\
                      & \rdTo(1,2)_{\phi}   \ldTo(1,2)_{\psi}&   \\
                       &  \pi & \\
\end{diagram}
$$
Let $G$ be an affine algebraic group scheme and let $\varrho:\,\pi \rar G$ be a representation of $\pi$ into $G$. The commutative Hopf algebra $\cO(G)$ defines a functor $\underline{\cO}(G)\,:\,\mathfrak{G} \rar \cAlg_k\,,\,\langle n \rangle \mapsto \cO(G)^{\otimes n}$. Further, for $(\langle n \rangle ,\phi)\,\in\,\mathfrak{G}/\pi$, one has the augmentation maps (depending on $\varrho$)
$$
\cO(G)^{\otimes n} \cong \langle n \rangle_G \stackrel{\phi_G}{\longrightarrow} \pi_G \stackrel{\varrho}{\rar} k \ .
$$
This enriches  $\underline{\cO}(G)$ to a functor $\underline{\cO}(G)_{\varrho}\,:\,\mathfrak{G}/\pi \rar \cAlg_{k/k}$. Combining
this last functor with linearization \eqref{linfun}, we define the
functor
\begin{equation} \la{uAd}
 \underline{\Ad}^*\varrho :=\,\cQ \circ \underline{\cO}(G)_{\varrho}:\ \mathfrak{G}/\pi \rar \ve_k\ ,\quad  (\langle n \rangle,\phi) \mapsto
 (\Ad^* \varrho)^{\oplus n}\ ,
\end{equation}
which assigns to $\, (\langle n \rangle,\phi) \in \mbox{Ob}(\mathfrak{G}/\pi)\,$ the vector space $\,\Omega(\langle n \rangle) \otimes_{\Z\langle n \rangle} \g^* \cong (\g^*)^{\oplus n}\,$ 
determined by the coadjoint representation
$$
\langle n \rangle \xrightarrow{\phi} \pi \xrightarrow{\varrho} G \xrightarrow{\Ad^*} \GL_k(\g^*) \ .
$$
In particular, if $ \varrho: \pi \to G $ is the trivial representation, the functor \eqref{uAd} is isomorphic to (the restriction to $ \mathfrak{G}/\pi $ of) the abelianization functor on $\mathfrak{G} $ tensored with $ \g^* $: we denote this functor by
\begin{equation} 
\la{natmod}
\underline{\g}^*:\,\mathfrak{G} \to \ve_k\ ,\quad \langle n\rangle \mapsto 
\langle n\rangle_{\rm ab} \otimes_{\Z} \g^*  = (\g^*)^{\oplus n}\ .
\end{equation}

Next, an object $\,(\langle m \rangle, \psi) \in \mbox{Ob}(\mathfrak{G}/\pi)\,$ defines a right $\mathfrak{G}/\pi$-module
$$
k[\langle m \rangle ,\psi]:\ (\langle n \rangle,\phi) \mapsto k[\Hom_{\mathfrak{G}/\pi}[(\langle n \rangle,\phi),(\langle m \rangle,\psi)]] \,=\,k[\psi^{-1}(\phi(x_1)) \times \ldots \times \psi^{-1}(\phi(x_n))]
\ .
$$
By Yoneda's Lemma ({\it cf.} \cite[Lemma 4.3]{BRY}), the module $k[\langle m \rangle ,\psi]$ is projective, and $k[\langle m \rangle ,\psi] \otimes_{\mathfrak{G}/\pi} {\mathcal N}\,\cong\, {\mathcal N}[(\langle m \rangle ,\psi)]$ for any left $\mathfrak{G}/\pi$-module $ {\mathcal N} $. 
More generally, any $ (\Gamma, \psi) \in {\rm Ob}(\mathtt{Gr}/\pi) $ defines a right $\mathfrak{G}/\pi$-module
\begin{equation} 
\la{gmodpimod} 
k[\Gamma, \psi]:\ (\langle n \rangle,\phi) \mapsto k[\Hom_{\mathfrak{G}/\pi}((\langle n \rangle,\phi),(\Gamma, \psi)) ]
\,=\,k[\psi^{-1}(\phi(x_1)) \times \ldots \times \psi^{-1}(\phi(x_n))]\ .
\end{equation}
If $\Gamma$ is a based free group, the module $k[\Gamma, \psi]$ is flat and $k[\Gamma, \psi] \otimes_{\mathfrak{G}/\pi} {\mathcal N}\,\cong\, \tilde{{\mathcal N}}[(\Gamma, \psi)]$, where the functor on the right-hand side is the left Kan extension of $ {\mathcal N} $ along the inclusion of 
$ \mathfrak{G}/\pi $ into the category of all based free groups equipped with a homomorphism to $\pi$.

Now, for a space $X$,  set $ \pi := \pi_1(X) $. Then, $\lgr{X}$ is naturally  a simplicial object in $\mathtt{Gr}/\pi$, with the structure map $ p: \lgr{X} \to \pi_0(\lgr{X}) $ being the canonical projection onto $ \pi_0(\lgr{X}) = \pi \,$.
Thus we can associate to $X$ the simplicial right $\mathtt{Gr}/\pi$-module $k[\lgr{X},p]$ as in \eqref{gmodpimod}. The proof of the following theorem is a trivial modification of that of \cite[Theorem 4.2]{BRY}; we leave its details to the interested reader.
\bthm \la{fhint}
For any representation $\varrho\,:\,\pi_1(X) \rar G(k)$, there is a natural isomorphism
$$ \H_{\ast}[\mathrm{T}^{\ast}_{\varrho}\DRep_G(X)] \,\cong\,\H_{\ast}[k[\lgr{X}, p]\, \otimes^{\L}_{\mathfrak{G}/\pi}\, \underline{\Ad}^*\varrho]\,,$$
where $\underline{\Ad}^*\varrho $ is the left $\mathfrak{G}/\pi$-module
defined in \eqref{uAd}.
\ethm

By standard homological algebra (see, e.g., \cite[Appl.~5.7.8]{W}), Theorem~\ref{fhint} implies 
the existence of the spectral sequence
\begin{equation}
\la{specseq} 
{E}^2_{p,q}\,=\,\Tor^{\mathfrak{G}/\pi}_p(\pi_q(k[\lgr{X},p]), \,\underline{\Ad}^*\varrho) \,\implies\, \H_{p+q}[\mathrm{T}^{\ast}_{\varrho}\DRep_G(X)]\ ,
\end{equation}
where $ \pi_q(k[\lgr{X}, p]) $ denotes the $q$-th homotopy group of the simplicial $\mathfrak{G}/\pi$-module $ k[\lgr{X}, p] $. When $X$ is simply connected, the spectral sequence \eqref{specseq} can be written in a 
more geometric form. In this case, the module $ k[\lgr{X}, p] $  may be identified with the right $\mathfrak{G}$-module $\underline{\H}_q(\Omega X;k)$, which assigns to $\langle m \rangle$ the degree $q$ component of the graded vector space $\H_\ast(\Omega X;k)^{\otimes m}$, where $\Omega X$ denotes the (Moore) based loop space of $X$ ({\it cf.} \cite[Theorem 4.3]{BRY}), and we have
\bcor
\la{corfhint1}
For simply connected $X$, there is a first quadrant spectral sequence
$$ {E}^2_{p,q}\,=\,\Tor^{\mathfrak{G}}_p(\underline{\H}_q(\Omega X;k),\, \underline{\g}^*) \,\implies\, \H_{p+q}[\mathrm{T}^{\ast}_{\varrho}\DRep_G(X)]\ ,
$$
where $ \underline{\g}^* $ is the left $ \mathfrak{G}$-module defined in \eqref{natmod}.
\ecor

On the other extreme, in the case of aspherical spaces, the spectral sequence  \eqref{specseq} degenerates giving the following
\bcor
\la{corfhint2}
For any discrete group $\Gamma$, and for any representation
$\varrho\,:\,\Gamma \rar G(k)$,
$$ \H_{\ast}[\mathrm{T}^{\ast}_{\varrho}\DRep_G(\Gamma)] \,\cong\,\Tor^{\mathfrak{G}/\Gamma}_{\ast}(k , \,\underline{\Ad}^*\varrho)\ .$$
\ecor
\bproof
We need to check that for $X=\B\Gamma$ and $\pi=\Gamma$, the module $ \pi_q(k[\lgr{X}, p]) $ is isomorphic to $k$ for $\,q=0\,$ and vanishes for $\,q>0$. Recall that $ p:\,\lgr{X} \rar \Gamma$ denotes 
the canonical projection to the discrete simplicial group $\Gamma = \pi_0(\lgr{X})$. For any 
$\, \gamma \in \Gamma\,$, $ p^{-1}(\gamma)$ is then a simplicial set. The $\mathfrak{G}/\Gamma$-module $\lgr{X}$ assigns to $(\langle n \rangle,\phi)$ the simplicial $k$-vector space $k[ p^{-1}(\phi(x_1)) \times \ldots \times p^{-1}(\phi(x_n))]\,\cong\, \otimes_{i=1}^n k[p^{-1}(\phi(x_i))]$. By Kunneth's Theorem, there is an isomorphism of graded $k$-vector spaces
$$ \pi_{\ast}(k[p^{-1}(\phi(x_1)) \times \ldots \times p^{-1}(\phi(x_n))])\,\cong\, \otimes_{i=1}^n \pi_{\ast}(k[p^{-1}(\phi(x_i))])\ .
$$
It therefore suffices to verify that for any $ \gamma \in \Gamma$, $\pi_0(k[p^{-1}(\gamma)])\,\cong\, k$, while $\pi_i (k[p^{-1}(\gamma)])\,\cong\,0$ for $i>0$. Indeed, for any space $X$, the map $\,p: \lgr{X} \rar \Gamma \,$ is a Kan fibration, with $ p^{-1}(\gamma) $ being its homotopy fiber.
Hence, we have a long exact sequence of homotopy groups
$$ 
\ldots \rar \pi_n[p^{-1}(\gamma)] \rar \pi_n(\lgr{X}) \rar \pi_n(\Gamma) \rar \ldots \rar \pi_0[p^{-1}(\gamma)] \rar \pi_0(\lgr{X}) \rar \pi_0(\Gamma) \ . $$
For $X = \B \Gamma $, the fibration $\lgr{X} \rar \Gamma$ is acyclic, i.e. $\pi_i(\lgr{X}) \,\cong\,\pi_i(\Gamma)$ for all $i \geq 0$. In this case, it follows from the above exact sequence that $p^{-1}(\gamma)$ is contractible. 
\eproof
\noindent
Comparing the result of Corollary~\ref{corfhint2} to that of Theorem~\ref{LinRep} for $ X = \B \Gamma $, we get
$$
\Tor^{\mathfrak{G}/\Gamma}_{\ast}(k , \,\underline{\Ad}^*\varrho)
\,\cong\, \left\{
                                                 \begin{array}{ll}
                                                    Z^1(\Gamma,\, \mathrm{Ad}\,\varrho)^{\ast} & \text{if }\ i=0\\*[1ex]

                                                    \H_{i+1}(\Gamma,\,\mathrm{Ad}^{\ast}\varrho) & \text{if }\ i>0 \ \\
                                                 \end{array}
                                                 \right.
$$
This last isomorphism can be obtained directly by using the canonical resolution
of the trivial module $ k $ in the category of right
$\mathfrak{G}/\Gamma$-modules.

\section{Character homology}
\la{S4}
Throughout this section we assume, for simplicity, that $k=\c$ and $G$ is a complex connected {\it reductive} affine algebraic group with Lie algebra $\g$. Recall that for any group $\Gamma$, the algebraic group $G$ acts naturally on the affine scheme $\Rep_G(\Gamma)$ by conjugation. This action induces an action of $G$ on the commutative algebra $ \Gamma_G $ representing $\Rep_G(\Gamma) $,
and we write $\Gamma_G^G$ for the subalgebra of $G$-invariants in $ \Gamma_G $. The affine scheme $\mathrm{Spec}(\Gamma_G^G)$ is called the {\it $G$-character scheme} of $\Gamma$ and denoted  $\Char_G(\Gamma)$. The assignment $\, \Gamma \mapsto \Gamma_G^G\,$ defines the subfunctor $\,(\,\mbox{--}\,)_G^G\,:\,\Gr \rar \cAlg_k\,$ of the representation functor which \mbox{---}
when extended to simplicial groups \mbox{---} has a total left derived functor
\begin{equation}
\la{Lchar}
\L(\,\mbox{--}\,)_{G}^G :\,\Ho(\mathtt{sGr}) \rar \Ho(\mathtt{s}\cAlg_k)\,\text{.}
\end{equation}

Now, for a space $ X \in \Top_{0,\ast} $, we fix a simplicial group model
$ \mathbf{\Gamma}X $ and set
$$
\mathcal{B}_G(X) :=
\mathcal{B}_G(\mathbf{\Gamma}X) :=
 \L(\mathbf{\Gamma}X)_{G}^G\ .
$$
Then, we define the {\it derived character scheme} of $ X $ in $G$
formally by $\,\DChar_G(X) := \Spec[\mathcal{B}_{G}(X)] \,$, i.e.
$\,\mathcal{B}_{G}(X)$ viewed as an object of the opposite category $\Ho(\scAlg_k)^{\mathrm{op}}$, and refer to the  homotopy groups $\,\pi_*[\mathcal{B}_{G}(X)] = \L_\ast(\mathbf{\Gamma}X)_{G}^G\,$ as the {\it character homology} of $X$ in $G\,$.

Since $G$ is a reductive algebraic group, the canonical morphism of functors $\,\L(\,\mbox{--}\,)_{G}^G \to \L(\,\mbox{--}\,)_{G}\,$ induces an isomorphism $\,\L_\ast(\,\mbox{--}\,)_{G}^G \cong
[\L_\ast(\,\mbox{--}\,)_{G}]^G $ which applied to a simplicial group $ \mathbf{\Gamma}X $, reads
$$
\pi_*[\mathcal{B}_G(X)]\, \cong\, \HR_{\ast}(X,G)^G
$$
where $ \HR_{\ast}(X,G)^G $ denotes the $G$-invariant part of representation homology of $X$ in $G$.  We will use this last isomorphism to identify the character homology of $X$ in $G$ with $\, \HR_{\ast}(X,G)^G $. In particular, we have
$$
\pi_0[\mathcal{B}_G(X)] \cong \HR_0(X,G)^G\,\cong\,\cO[\Char_G(\pi_1(X))]\ .
$$

Now, given a space $X$, fixing a representation $\varrho:\pi_1(X) \rar G$ is equivalent to giving an augmentation of the commutative algebra $\HR_0(X,G)$. This in turn, yields (by restriction) an augmentation $\varrho\,:\,\HR_0(X,G)^G \rar k$, which is equivalent to giving a homomorphism of simplicial commutative algebras $\varrho\,:\,\mathcal{B}_G(X) \rar k$. Applying the derived linearization functor $\L\cQ$ to the pair $(\mathcal{B}_G(X),\varrho)$, we obtain a simplicial vector space
  $$\T^{\ast}_{\varrho}\DChar_G(X)\,:=\,\L\cQ_{\varrho}[\mathcal{B}_G(X)]\,,$$
which may be interpreted geometrically as the derived cotangent complex to  $\DChar_G(X)$ at the character $\varrho$. Our goal is to compute the homology of the derived cotangent complex of the derived character scheme in terms of the homology of the local system determined by a representation.

To state our main theorem we recall some standard terminology from representation theory
(see, e.g., \cite{LuM, S}). If $\Gamma$ is a (discrete) group, a representation $\varrho\,:\,\Gamma \rar G $ is called {\it irreducible} if $\varrho(\Gamma)$ is not contained in any parabolic subgroup of $G $. Further, $\varrho$ is called {\it completely reducible} if for every parabolic subgroup $P \subset G$ containing $\varrho(\Gamma)$, there exists a Levi subgroup $ L \subseteq P $ such that
$ \varrho(\Gamma) \subseteq L \subseteq P\,$. Thus, a fortiori, any irreducible representation is completely reducible. In the case $G=\GL_n(k)$, a representation $\varrho\,:\,\Gamma \rar \GL_n(k)$ is irreducible if $k^n$ is a simple $\Gamma$-module (via $\varrho$) and completely reducible if $k^n$ is a semisimple $\Gamma$-module. It is known that the preimage in $\Rep_G(\Gamma)$ of any (closed) point in $\Char_G(\Gamma)$ contains at least one completely reducible representation (see, e.g. \cite[Sect. 11]{S}).
Note that for any representation $\varrho\,:\,\Gamma \rar G$, the centralizer of $\varrho(\Gamma)$ in $G$ contains the center $\mathcal{Z}(G)$ of $G$. Following \cite{JM} (see also \cite{S}), we call
a representation $\varrho$ {\it good} if its $G$-orbit in $ \Rep_G(\Gamma) $ is closed and its
stabilizer  in $G$ (i.e., the centralizer of $ \varrho(\Gamma) $) coincides with the center
$\mathcal{Z}(G)$. One can show that every good representation is irreducible
({\it cf.} \cite[Theorem~30]{S}), but the converse is not always true\footnote{For $G =\GL_n(k)$ or
$\mathrm{SL}_n(k)$, every irreducible representation in $G$ of a discrete group $ \Gamma $ is good. However, for $\,G=\mathrm{PSL}_2(k),\,\mathrm{SO}_n(k),\,\mathrm{Sp}_{2n}(k)\,$ (or any quotient of these groups by a finite subgroup), there exist discrete groups having `bad' irreducible representations in $G$  (see \cite[Sect. 4]{S}).}.

The following theorem, which is one of the main results of the present paper, is a natural generalization of \cite[Theorem 53]{S}.

\bthm
\la{lincharhom}
Let $X$ be a space such that $\pi_1(X)$ is finitely generated. Then  \\
$(i)$ For every completely irreducible representation $\varrho\,:\,\pi_1(X) \rar G(k) $, there are natural maps
 \begin{equation} \la{ecothom} \H_i[\T^{\ast}_{\varrho}\DChar_G(X)] \rar \H_{i+1}(X, \mathrm{Ad}^{\ast}\varrho)\,,\,\,\, \forall \, i \geq 0 \ .\end{equation}
 $(ii)$ If a representation $\varrho\,:\,\pi_1(X) \rar G(k) $ is good, then the maps \eqref{ecothom} are isomorphisms.
 %
 %
 \ethm
To prove Theorem~\ref{lincharhom} we need the following lemma.
 \blemma \la{freecothom}
 Let $F$ be a free group and let $\varrho\,:\,F \rar G$ be a good representation that factors through some finitely generated quotient of $F$. Then,
 $$\H_i[\L\cQ_{\varrho}(F_G^G)]=0 \ ,\quad \forall\, i>0\ . $$
 \elemma
 \bproof
Note that $F$ can be expressed as a direct limit of finitely generated free subgroups that surject onto $\varrho(F)$. We may therefore, assume without loss of generality that $F$ is a finitely generated free group. In this case, $\cO[\Rep_G(F)]$ is reduced, which implies that the representation scheme $\Rep_G(F)$ coincides with the associated variety. By \cite[Prop. 27]{S}, the set of irreducible representations of $F$ is Zariski open in $\Rep_G(F)$. Since $G$ is reductive, the projection $\Rep_G(F) \rar \Char_G(F)$ is a good categorical quotient, whence the set $\Char^i_G(F)$ of irreducible characters of $F$ is Zariski open in $\Char_G(F)$. By \cite[Corollary50]{S}, the set of characters of good representations is a smooth (Zariski) open subset of $\Char^i_G(F)$. It follows that the point $\varrho$ of $\Char_G(F)$ is smooth. It therefore, suffices to verify that if $B \in \cAlg_{k/k}$ is such that the augmentation on $B$ defines a smooth point on $\mathrm{Spec}(B)$, then $\H_i[\L\cQ(B)]=0$ for $i>0$. To see this, note that by the definition of Andr\'{e}-Quillen homology,
 $$ \H_i[\L\cQ(B)]\,\cong\, {D}_i(B|k, k)\,,$$
 where ${D}$ stands for Andr\'{e}-Quillen homology, and where the second argument on the right hand side is $k$ equipped with $B$-module structure through the given augmentation. Now, if $\mathfrak{m}$ is any maximal ideal of $B$, by \cite[Sect. 6.5]{Iyer} we have:
 $$ {D}_i(B|k, k)_{\mathfrak{m}}\,\cong\, {D}_i(B_{\mathfrak{m}}|k, k_{\mathfrak{m}}) \text{ for all } i \geq 0\ .$$
 If $\mathfrak{m}$ is not the augmentation ideal of $B$, $k_{\mathfrak{m}}=0$. Hence  ${D}_i(B|k,k)_{\mathfrak{m}}=0$. On the other hand, if $\mathfrak{m}$ is the augmentation ideal of $B$, then $B_{\mathfrak{m}}$ is smooth and $k_{\mathfrak{m}}=k$, which implies that ${D}_i(B_{\mathfrak{m}}|k, k_{\mathfrak{m}})$ vanishes for $i>0$ (see. e.g. \cite[Theorem 9.6]{Iyer}). Thus,
${D}_i(B|k, k)=0$ for $i>0$, as desired.
 \eproof

\bproof[Proof of Theorem \ref{lincharhom}]
It follows from \cite[Theorem 53]{S} (by taking linear duals) that if $\Gamma$ is a finitely generated group, and $\varrho\,:\,\Gamma \rar G$ is a completely reducible representation, then the natural projection of schemes $\Rep_G(\Gamma) \rar \Char_G(\Gamma)$ induces a map of $\c$-vector spaces
\begin{equation} \la{cotspacechar} \T^{\ast}_{\varrho}\Char_G(\Gamma)\,:=\,\cQ_{\varrho}(\Gamma_G^G) \rar \H_1(\Gamma,\mathrm{Ad}^{\ast}\varrho) \end{equation}
that is natural in $\Gamma$. If, moreover, $\varrho$ is good, then the map \eqref{cotspacechar} is an isomorphism. Since the functor $\cQ$ commutes with direct limits (as does the functor $(\,\mbox{--}\,)_G$) and since homology commutes with direct limits, the map \eqref{cotspacechar} is defined even if $\Gamma$ is not finitely generated provided that the completely reducible representation $\varrho$ factors through a finitely generated quotient of $\Gamma$. This can be seen by representing $\Gamma$ as a direct limit of finitely generated subgroups that surject onto $\varrho(\Gamma)$. If in addition $\varrho$ is good, the map \eqref{cotspacechar} is an isomorphism.

Next, we note that if $\Gamma$ is a free group, the complex $0 \rar \Omega(\Gamma) \stackrel{\alpha}{\hookrightarrow} \Z[\Gamma] \rar 0$ is a free resolution of the trivial module $\Z$ in the category $\Z[\Gamma]$-modules (see, e.g. \cite[Chap. IV.2, Ex.3]{KB}). Hence, in this case there is the exact sequence
$$0 \rar \H_1(\Gamma, \mathrm{Ad}^{\ast}{\varrho}) \rar \Omega(\Gamma) \otimes_{\Z[\Gamma]} \g^{\ast} \rar \g^{\ast} \rar \H_0(\Gamma,\mathrm{Ad}^{\ast}{\varrho}) \rar 0 \,, $$
where $\g^{\ast}$ acquires a $\Gamma$-module structure via the representation $\mathrm{Ad}^{\ast}\varrho$. As a result, we have
$$ \H_1(\Gamma,\mathrm{Ad}^{\ast}\,{\varrho}) \,\cong\, \Ker(\alpha \otimes_{\Z[\Gamma]}\id_\g^{\ast})\ .$$
Since the subspace of $\Gamma$-invariants $(\g^{\ast})^{\Gamma}$ is precisely the Lie algebra $\mathfrak{z}$ of the centralizer $C(\varrho)$ of $\varrho(\Gamma)$, the projection $\g^{\ast} \rar \H_0(\Gamma,\mathrm{Ad}^{\ast}\varrho)$ is the linear dual of the inclusion $\mathfrak{z} \hookrightarrow \g$. The image of the map $\alpha \otimes_{\Z[\Gamma]}\id_\g^{\ast}$ is therefore, the subspace $(\g/\mathfrak{z})^{\ast}$ of $\g^{\ast}$. Thus, there is an exact sequence of vector spaces
\begin{equation} \la{exact1} 0 \rar K(\Gamma) \rar \Omega(\Gamma) \otimes_{\Z[\Gamma]}\g^{\ast} \rar (\g/\mathfrak{z})^{\ast} \rar 0 \,, \end{equation}
where $  K(\Gamma):= \Ker(\alpha \otimes_{\Z[\Gamma]}\id_\g^{\ast})$. Next, observe that since the natural map $\lgr{X} \rar \pi_0(\lgr{X})\,\cong\,\pi_1(X)$ is surjective in every simplicial degree, any completely reducible (resp., good) representation of $\pi_1(X)$ restricts to a completely reducible (resp.,good) representation of the (free) group $(\lgr{X})_n$ for every $n$. It follows that for a completely reducible representation $\varrho\,:\,\pi_1(X) \rar G$, there is a map of simplicial vector spaces
\begin{equation} \la{qrho} \cQ_{\varrho}[(\lgr{X})_G^G] \rar K(\lgr{X})\,,\end{equation}
which is an isomorphism by \cite[Theorem 53 (2)]{S} when $\varrho$ is good. Since $\lgr{X}$ is semi-free, $K(\lgr{X})$ fits into the long exact sequence of simplicial vector spaces
\begin{equation} \la{exact} 0 \rar K(\lgr{X}) \rar \Omega(\lgr{X}) \otimes_{\Z[\lgr{X}]} \g^{\ast} \rar (\g/\mathfrak{z})^{\ast} \rar 0\,,\end{equation}
where $(\g/\mathfrak{z})^{\ast}$ is viewed as a discrete simplicial vector space. By (the proof of) Theorem \ref{LinRep},
$$ \H_i[\Omega(\lgr{X}) \otimes_{\Z[\lgr{X}]} \g^{\ast}]\,\cong \,
\left\{
\begin{array}{ll}
\mathrm{Z}^1(\pi_1(X), \mathrm{Ad}\,\varrho)^{\ast} & \text{ if } i=0\\*[1ex]
\H_{i+1}(X,\mathrm{Ad}^{\ast}\varrho) & \text{ if } i>0 \ .\\
\end{array}
\right.
$$
It follows from the long exact sequence of homologies associated with the sequence \eqref{exact} that
 \begin{equation} \la{kerhom} \H_i[K(\lgr{X})]\,\cong\,\H_{i+1}(X,\mathrm{Ad}^{\ast}\varrho)\ .\end{equation}
Therefore, the composition in $\Ho(\sMod_{k})$
$$ \T^{\ast}_{\varrho}\DChar_G(X) \,\cong\, \L\cQ_{\varrho}[(\lgr{X})_G^G] \rar \cQ_{\varrho}[(\lgr{X})_G^G] \rar K(\lgr{X}) $$
induces the map
$$ \H_i[\T^{\ast}_{\varrho}\DChar_G(X)] \rar \H_{i+1}(X,\mathrm{Ad}^{\ast}\varrho)\,,\,\,\,\forall\, i \geq 0\ .$$
This proves part $(i)$. If $\varrho$ is good, then the map \eqref{qrho} of simplicial vector spaces is an isomorphism. Part $(ii)$ therefore, follows from \eqref{kerhom} once we verify that
\begin{equation} \label{lqtoq}  \L\cQ_{\varrho}[(\lgr{X})_G^G] \rar \cQ_{\varrho}[(\lgr{X})_G^G] \end{equation}
is an isomorphism in $\Ho(\sMod_{\c})$. Since the group $(\lgr{X})_n$ is free for each $n$, and since $\varrho$ induces a good representation on $(\lgr{X})_n$ for each $n$, we have
$$\H_i[\L\cQ_{\varrho}((\lgr{X})_n)]=0\,,\,\,\,\forall\,i \geq 0 $$
by Lemma \ref{freecothom}. It follows that for the pair $((\lgr{X})_G^G, \varrho)\,\in\,\scAlg_{k/k}$, the spectral sequence in Proposition \ref{aqhomspecseq} collapses on the ${\it E}^1$ page. The map \eqref{lqtoq} therefore induces an isomorphism on homologies, as desired.
\eproof

\noindent
We conclude this section with two examples.

\begin{example}
Consider the Heisenberg manifold $X:=H_3(\R)/H_3(\Z)$, where the Heisenberg group $H_3(R)$ of a ring $R$ is the group of unipotent upper-triangular $3 \times 3$ matrices with entries in $R$. Since $H_3(\R)$ is homeomorphic to the contractible space $\R^3$, $X$ is an aspherical manifold with fundamental group $H_3(\Z)$. It is well known that $H_3(\Z)_{\ab}\,\cong\,\Z^2$, where $(\,\mbox{--}\,)_{\ab}$ stands for abelianization. Moreover,
$$\Char_{\mathrm{SL}_2}(X)\,\cong\,\Char_{\mathrm{SL}_2}(\Z^2) \sqcup \{\varrho\}\,,$$
where the first component consists of the characters of representations $H_3(\Z) \rar \mathrm{SL}_2(\c)$ that factor through the abelianization. It is easy to verify that no representation of $\Z^2$ in $\mathrm{SL}_2(\c)$ is good. Hence, none of the characters in the first component are characters of good representations. On the other hand, $H_3(\Z)$ has the presentation
$$H_3(\Z)\,\cong\, \langle \alpha,\beta,\gamma\,|\, \alpha \beta =\beta \alpha \,,\,\alpha \gamma =\gamma \alpha\,,\,
[\gamma, \beta] =\alpha\,\rangle \ . $$
In terms of the above presentation, the unique character $\varrho$ in the second component of $\Char_{\mathrm{SL}_2}(X)$ is the character of the representation
$$ \varrho\,:\,H_3(\Z) \rar \mathrm{SL}_2(\c)\,,\,\,\alpha \mapsto -\id\,,\,\beta \mapsto \left( \begin{array}{cc}
                                 i & 0 \\
                                 0 & -i \\
                                 \end{array} \right)\,,\,
                                 \gamma \mapsto \left( \begin{array}{cc}
                                                         0 & -1\\
                                                         1 & 0\\
                                                         \end{array} \right)\ . $$
It can be checked by a direct computation that the representation $\varrho$ is good, and that the cohomologies $\H^0(X,\mathrm{Ad}\,\varrho)$ and $\H^1(X,\mathrm{Ad}\,\varrho)$ vanish.
By Poincar\'{e} duality the homologies
$ \H_2(X,\mathrm{Ad}^{\ast}\varrho)$ and $\H_3(X,\mathrm{Ad}^{\ast}\varrho)$ vanish as well.
Since the local systems $\mathrm{Ad}\,\varrho$ and $\mathrm{Ad}^{\ast}\varrho$ are isomorphic via the Killing form,
$\H^i(X,\mathrm{Ad}\,\varrho)=0$ for $i=0,1$. By universal coefficients, $\H_i(X,\mathrm{Ad}^{\ast}\varrho)=0$ for $i=0,1$ as well. It follows from Theorem \ref{lincharhom} that
$$
\H_i[\T^{\ast}_{\varrho}\DChar_{\mathrm{SL}_2}(X)]\,=\,0 \ ,\ \, \forall\,i \geq 0\ . 
$$
It follows from Theorem~\ref{bounded_homology}$(a)$ that the localization
$\HR_{\ast}(X,\mathrm{SL}_2)^{\mathrm{SL}_2}_{\varrho}$ of the $\mathrm{SL}_2$-character homology of $X$ vanishes in positive degree. Since $\cO[\Char_{\mathrm{SL}_2}(X)]_{\varrho}=k$, we conclude that
$$\HR_{\ast}(X,\mathrm{SL}_2)^{\mathrm{SL}_2}_{\varrho}\,\cong\,k \ . $$
\end{example}

\begin{example}
We now illustrate that the map \eqref{ecothom} in Theorem \ref{lincharhom} is far from being an isomorphism in general, even when $\varrho$ is completely reducible. Let $X=\c\mathbb{P}^r,\, r \geq 1$. Since $X$ is simply connected, $\pi_1(X)$ has the unique (trivial) representation $\varrho$ into any complex reductive group $G$. Clearly, $\varrho$ is completely reducible, and the local system $\mathrm{Ad}^{\ast}{\varrho}$ is the constant local system $\g^{\ast}$.  Let $m_1,\ldots,m_l$ denote the exponents of the Lie algebra $\g$. Then there is an isomorphism of graded vector spaces
\begin{equation}
\la{eiso1}
\H_{\ast}[\T^{\ast}_{\varrho}\DChar_G(\c\mathbb{P}^r)] \,\cong\, \bigoplus_{i=1}^l \big(k.\xi_1^{(i)} \oplus k.\xi_3^{(i)} \oplus \ldots \oplus k.\xi_{2r-1}^{(i)} \big)\,,
\end{equation}
where the basis element $\xi_{2s-1}^i$ has homological degree $2rm_i+2s-1$. On the other hand,
\begin{equation}\la{eiso2}
\H_{i+1}(\c\mathbb{P}^r,\g^{\ast})\,\cong\, \left\{
\begin{array}{ll}
\g^{\ast} & \text{ if } i=1,3,5,\ldots, 2r-1 \\*[1ex]
0 & \text{ otherwise.}\\
\end{array}
\right.
\end{equation}
Indeed, by \cite[Theorem~6.4, Corollary 6.3]{BRY}, there is an isomorphism of graded commutative algebras
\begin{equation} \la{charhomcpr} \Sym_{k}\big(k.\xi_1^{(i)} \oplus \ldots \oplus k.\xi_{2r-1}^{(i)}\,|\, i=1.2.\ldots,l \big) \,\cong\, \HR_{\ast}(\c\mathbb{P}^r, G)^G  \,, \end{equation}
where the generator $\xi_{2s-1}^i$ has homological degree $2rm_i+2s-1$. From  {\it loc. cit.} (see also Sect. 6.1 therein as well as \cite[Sect. 7]{BFPRW}), it is clear that the isomorphism \eqref{charhomcpr} is induced by a quasi-isomorphism of commutative DG algebras between the symmetric algebra of a complex whose homology is isomorphic to
$$\bigoplus_{i=1}^l \big(k.\xi_1^{(i)} \oplus k.\xi_3^{(i)} \oplus \ldots \oplus k.\xi_{2r-1}^{(i)} \big) $$
and a commutative DG algebra representing $\DChar_G(\c\mathbb{P}^r)$ in $\Ho(\cDGA_{\c}^+)$. The 
isomorphism \eqref{eiso1} follows immediately from the above observation. The verification of
\eqref{eiso2} is an easy  exercise which we leave to the reader.
\end{example}


%
\begin{remark}
In order to recover the groups
$\, \H_{\ast}(X, \mathrm{Ad}^{\ast}\varrho) \,$ as a cotangent homology for {\it every}\, representation $ \varrho $, one should replace the derived character scheme $ \DChar_G(X) $ in Theorem~\ref{lincharhom} by the derived mapping stack $ \textbf{Map}(X, BG) $  (see Appendix~A, Section~\ref{TVcon}). The derived stack $ \textbf{Map}(X, BG) $ may be viewed as a homotopy quotient of the derived representation scheme $ \DRep_G(X) $, and thus -- from the derived geometric point of view -- it is a more natural object to consider. 
\end{remark}

\section{Vanishing theorems}
\la{S5}
In this section, we prove vanishing theorems for representation homology of some classical spaces. 
One important consequence of these results is the existence of a well-defined virtual fundamental class of the corresponding derived representation schemes. The proofs are based on Proposition~\ref{Quillen_specseq} and Theorem~\ref{bounded_homology} from Section~\ref{S2}. We begin by reformulating these results in the form that we will use in this section.

Let $X$ be a CW-space having a cellular model with finitely many cells in each dimension. Given an algebraic group $G$ over $k$, a representation $\varrho\,:\,\pi_1(X) \rar G(k)$ corresponds to a $k$-point in $\Rep_G[\pi_1(X)]$. One may therefore, consider the localization of the representation homology at $\varrho$ which we denote by $\HR_{\ast}(X,G)_{\varrho}$. Applying Proposition \ref{Quillen_specseq} to the derived representation algebra $\mathcal{A}_G(X)$, we obtain the spectral sequence ({\it cf.} \eqref{CKsp})
\begin{equation*}
E_{p,q}^2 = \pi_{p+q}[\Sym^p(\T_{\varrho}^*\DRep_G(X))]\ \ \Longrightarrow\ \ \widehat{\HR}_{p+q}(X,G)_{\varrho} \,,
\end{equation*}
where $\widehat{\HR}_{\ast}(X,G)_{\varrho}$ stands for the completion of the local graded ring ${\HR}_{\ast}(X,G)_{\varrho}$ with respect to the filtration by powers of the maximal ideal of $\HR_0(X,G)_{\varrho}$ defined by $\varrho$. In this situation, Theorem~\ref{bounded_homology}$(a)$ asserts:  if 
for every representation $\varrho\,:\,\pi_1(X) \rar G(\bar{k})$ over the algebraic closure $\bar{k}$ of $k$,$\,\H_i[\T_{\varrho}^*\DRep_G(X)]=0\,$ for all $i \geq 2$ and $\,\dim_{\bar{k}} \H_1[\T_{\varrho}^*\DRep_G(X)] \leq N\,$, then $\HR_n(X,G)=0$ for all $\,n>N\,$.

\subsection{Virtually free groups}
Recall that a group is called {\it virtually free} if it contains a free group as a subgroup of finite index. The following vanishing theorem is our second main result.
\bthm \la{vfgroups}
Let $\Gamma$ be a f.g. virtually free group. Then, $ \Rep_G(\Gamma) $ is a smooth $k$-scheme, and the natural map $\mathcal{A}_G(\Gamma) \rar \cO[\Rep_G(\Gamma)]$ is an isomorphism in $\Ho(\scAlg_k)$. In particular, 
$$
\HR_i(\Gamma, G)=0  \quad \text{for all}\quad i>0\ . 
$$
\ethm
\bproof
By Lemma \ref{cspacevf} below, ${\rm B}\Gamma$ has a semi-free simplicial group model $\mathbf{\Gamma}$ with finitely many generators in each degree (see Section \ref{simpgr}). Hence, for each $ n \ge 0 $, 
the degree $n$ component of the simplicial commutative algebra $\mathcal{A}_G(\Gamma) \cong \mathbf{\Gamma}_G $ is isomorphic to the tensor product of finitely many copies of (the finitely generated algebra) $\O(G)$. It follows then from Proposition~\ref{finite_type_equiv} that $\mathcal{A}_G(\Gamma) $ is of quasi-finite type.

Let $\g$ denote the Lie algebra of $G({k})$. Any augmentation $\varrho\,:\,\mathcal{A}_G(\Gamma) \rar {k}$ corresponds to a representation $\varrho\,:\,\Gamma \rar G({k})$. Composing this representation with the coadjoint action of $G$,
we get a group homomorphism $\mathrm{Ad}^*\,\varrho\,:\,\Gamma \rar 
\GL_k(\g^*) $, which equips $\g^*$ with the structure of a $\Gamma$-representation.
By Theorem~\ref{LinRep}, there are natural isomorphisms 
\begin{equation}
\la{macL}
\H_i[\L \cQ_\varrho(\mathcal{A}_G(\Gamma))]\,\cong\, \H_{i+1}(\Gamma,\,\mathrm{Ad}^{\ast}\varrho)\,\cong\, \HH_{i+1}({k}[\Gamma],\, \widetilde{{\rm Ad}^{\ast}}\varrho) \ ,\quad  i >0 \ ,
\end{equation}
where $\widetilde{{\rm Ad}^{\ast}}\varrho$ denotes the $ k[\Gamma]$-bimodule $ \g^* $ with the above $\Gamma$-action on the left and the trivial $\Gamma$-action on the right. The last isomorphism in \eqref{macL} is the classical Mac Lane Isomorphism (see, e.g., \cite[Prop. 7.4.2]{L}). Now, since $\Gamma$ is finitely generated virtually free, by \cite[Theorem 2]{LeB}, its group algebra ${k}[\Gamma]$ is quasi-free in the sense of 
Cuntz-Quillen \cite{CQ}.  It follows then from \cite[Prop. 3.3]{CQ} that ${k}[\Gamma]$ has a projective bimodule resolution of length at most two. Hence $\HH_p({k}[\Gamma],\, \widetilde{{\rm Ad}^{\ast}}\varrho)=0$ for $p>1$, and therefore  
$\H_i[\L\cQ_\varrho(\mathcal{A}_G(\Gamma))]=0$ for $i>0$. The desired result follows now from Theorem~\ref{bounded_homology}$(b)$.
\eproof
\blemma \la{cspacevf}
The classifying space $\, {\mathrm B}\Gamma $ of a finitely generated virtually free group has a reduced CW-model with finitely many cells in each dimension.
\elemma
\bproof
First, we note that the classifying space $BA$ of a finite group $A$ has a reduced simplicial model with finitely many simplices in each degree. It therefore, has a reduced CW model with finitely many cells in each dimension.

Next, we recall from \cite[Section 1.B]{Hat} that if
$A \stackrel{\varphi}{\hookleftarrow} C \stackrel{\psi}{\hookrightarrow} D$ is a diagram of groups, with both arrows injective, then the classifying space of the amalgamated free product $A \ast_C D$ is constructed by taking $BA, BC, BD$ and the mapping cylinders of $B\varphi$ and $ B\psi$, identifying the two ends of the mapping cylinder of $B\varphi$ with $BC$ and $BA$ and the two ends of the mapping cylinder of $B\psi$ with $BC$ and $BD$. If $BA, BC$ and $BD$ are taken to be reduced CW-models with finitely many cells in each dimension, this construction gives a CW model of the classifying space of $A \ast_C D$ with finitely many cells in each dimension. This CW complex, however contains three $0$-cells that are part of a contractible $1$-dimensional CW-subcomplex: if $x$, $y$ and $z$ denote the $0$-cells in $BA, BC$ and $BD$ respectively, this subcomplex has an edge between $x$ and $y$ and an edge between $y$ and $z$. Identifying this subcomplex with a single point gives us the required reduced CW-model for $A \ast_C D$.

We also note that given two monomorphisms $\varphi,\psi\,:\, C \hookrightarrow A$, the classifying space of the corresponding HNN extension is constructed from reduced CW-models $BC$ and $BA$ by gluing cells of dimension two or more to $BA \vee S^1$, the glued cells being the cells of $BC \times I$ that are of dimension two or more (see \cite[Example 1.B.13]{Hat}). It follows that if $BA, BC$ have finitely many cells in each dimension, the classifying space of the corresponding HNN extension also has a reduced CW model with finitely many cells in each dimension.
Finally, we recall that any finitely generated virtually free group is a fundamental group of a graph of finite groups (see, e.g., \cite[Section 2]{LeB}), i.e. can be constructed from a finite collection of finite groups in finitely many steps by taking (injective) amalgamated free products and (injective) HNN extensions as above. This proves the desired lemma.
\eproof

Since any finite group is virtually free, Theorem~\ref{vfgroups} implies
\bcor \la{finitegps} 
If $ \Gamma $ be a finite group, then $ \Rep_G(\Gamma) $ is smooth and $\,\HR_i(\Gamma,G)=0$ for all $i>0$.
\ecor
We remark that the smoothness of $ \Rep_G(\Gamma) $ for a finite group $ \Gamma $ is a well-known
result, which follows immediately from A. Weil's classical Rigidity Theorem \cite{Weil} (see, e.g., \cite[Corollary 
45]{S}). Weil's Theorem, however, doesn't apply to virtually free groups, since for such groups 
$ \H^1(\Gamma,\, {\rm Ad}\,\varrho) \not= 0 $ in general.

The following result shows that the higher representation homology also vanishes if we take $G$ to be a 
finite group.
\bprop \la{rephomintofgp}
Let $X$ be an arbitrary space.
If  $G$ is a finite algebraic group, then  $\,\HR_i(X,G)\,=\,0 \,$ for all $\, i > 0 $.
\eprop
\bproof
Let $Y$ be a CW-complex with finitely many cells in each dimension. Then $Y$ has a semi-free simplicial group model $\mathbf{\Gamma}$ with finitely many generators in each degree. It follows (as in the proof of Theorem \ref{vfgroups}) that $\mathcal{A}_G(Y)$ is of quasi-finite type.  Note that the Lie algebra $\g$ of $G({k})=G$ is zero. It follows that the local system $\mathrm{Ad}\,\varrho$ on $Y$ is zero for any augmentation $\varrho\,:\,\mathcal{A}_G(Y) \rar {k}$ (which corresponds to a representation $\varrho\,:\,\pi_1(Y) \rar G({k})=G$). Hence, by Theorem \ref{LinRep}, 
$$\H_i[\T^{\ast}_{\varrho}\DRep_G(Y)] = 0 \ , \ \,\forall\,i \ge 0 \ . $$
Theorem~\ref{bounded_homology}$(b)$ then implies that $\pi_i[\mathcal{A}_G(Y)]=0$ for all $i>0$. Thus, $\HR_i(Y,G)=0$ for all $i>0$. Since any space $X$ can be realized as a direct limit of CW-complexes with finitely many cells in each dimension, and since representation homology commutes with direct limits, the desired result holds.
\eproof
Finally, combining Theorem \ref{vfgroups} with \cite[Corollary 4.3]{BRY}, we have
\bcor \la{cvfgrps}
For any f.g.  virtually free group $\,\Gamma$, we have
 \begin{equation} \la{vftor} \Tor^{\mathfrak{G}}_i(k[\Gamma],\O(G))=0 \ \text{ for all }\ i>0\ . 
\end{equation}
\ecor

\subsection{Surfaces}
\la{S5.2} Let $G$ be a connected reductive affine algebraic group over $\c$. Recall that a representation $\, \varrho: \Gamma \to G \,$ is called (scheme theoretically) smooth if the corresponding (closed) point of the representation scheme $ \Rep_G(\Gamma) $ is simple, i.e. belongs to a unique irreducible component of $ \Rep_G(\Gamma) $ the dimension of which is equal to
the dimension of $ \T_\varrho \Rep_G(\Gamma)\,$. We note that every smooth representation is reduced
(i.e. corresponds to a reduced point of $ \Rep_G(\Gamma)$), and the set of all smooth representations
is Zariski open in $ \Rep_G(\Gamma)_{\rm red}\,$. Since the set of irreducible representations is also
Zariski open in $ \Rep_G(\Gamma)_{\rm red}\,$, ``most'' irreducible representations are smooth.

For surface groups $\,\Gamma = \pi_1(\Sigma_g) $, $ g \ge 2 \,$, 
a theorem of Goldman \cite{Go1} implies that $ \varrho \in  \Rep_G(\Gamma) $ is smooth if and only if $\, \dim\, [C(\varrho)/{\mathcal Z}(G)] = 0 \,$, where $ C(\varrho) $ is the centralizer of $ \varrho(\Gamma) $ in $G$. In particular, the good representations of surface groups
(i.e., the ones for which $ C(\varrho) = {\mathcal Z}(G) $) are all smooth.
More generally, it is known  (see \cite[Prop. 37]{S}) that {\it every}\, irreducible representation of 
$ \pi_1(\Sigma_g) $ is smooth. 

In this section, we prove
\bthm \la{vrephomsurfaces}
Let $\Sigma_g$ be a closed connected orientable surface of genus $g \geq 1$. Then

\vspace{1ex}

$(i)$ $\HR_i(\Sigma_g,G)=0\ $ for all $\ i >\dim G \,$.

\vspace{1ex}

$(ii)$ If $\varrho\,\in\,\Rep_G[\pi_1(\Sigma_g)]$ is smooth, then
$$\HR_i(\Sigma_g,G)_{\varrho}=0 \ \,\text{for }\, i>\dim_{\varrho}\Rep_G[\pi_1(\Sigma_g)] - (2g-1)\dim G\,,$$
 where $\dim_{\varrho} $ denotes local dimension of the irreducible component of $\Rep_G[\pi_1(\Sigma_g)]$ containing $\varrho$.
\ethm

\remark \, Part $(i)$ of Theorem \ref{vrephomsurfaces} can also be proven using the description of representation homology of surfaces in terms of Hochschild homology (see \cite[Sect. 7.1]{BRY}) along with the fact that $\cO(G)$ (viewed as a $\cO(G \times G)$-module via multiplication) {\it locally} has a projective resolution of length $\dim G $. Part $(ii)$, however, does not appear to be accessible from these considerations.

\vspace{1ex}

In general, Part $(ii)$ implies sharper vanishing statements than Part $(i)$: for example,
\bcor \la{vrephomtorus}
If $G$ is a simply-connected complex reductive group, then, for every smooth representation $\varrho\,:\,\Z^2 \rar G(\c)$, we have $\,\HR_i(\mathbb{T}^2, G)_{\varrho} = 0 \,$ if $\,i>\mathrm{rank}\, G\,$.
\ecor
\bproof
By \cite[Theorem C]{Ri} (which also holds for $\GL_n$ as the centralizer of any semisimple element of $\GL_n$ is connected), $\Rep_G(\Z^2)= \overline{G \cdot (T \times T)}$, where $T$ is a maximal torus of $G$ and $T \times T$ is viewed as a subscheme of $\Rep_G(\Z^2)$. It follows that $\Rep_G(\Z^2)$ is irreducible and $\dim\,\Rep_G(\Z^2) = \dim[G \cdot (T \times T)]$. The dimension of the generic fibre of the projection to the first factor $G \cdot (T \times T) \subset G \times G \rar G$ is equal to the dimension of the fibre over a regular point in $T$. This dimension is easily verified to be $\mathrm{rank}\,G\,$. Thus, $\dim \Rep_G(\Z^2) =\dim G +\mathrm{rank}\, G\,$. The desired result is now immediate from Theorem \ref{vrephomsurfaces}$(ii)$\,.
\eproof
\remark\ As noted in the above proof, Corollary \ref{vrephomtorus} holds for the groups $\GL_n $,$\, n \geq 1\,$ as well. In this case, however, we have a stronger statement: by computations of \cite[Sect. 7.1]{BRY}, $\HR_{\ast}(\mathbb{T}^2, \GL_n)$ is isomorphic to the homology of the Koszul complex
$$\big(k[X, Y, T][\det(X)^{-1}, \,\det(Y)^{-1}],\, dT=[X,Y]\big)\ . $$
By \cite[Theorem 27]{BFR}, this last homology vanishes in degree $> n$.
\bcor
\label{vrephomhighergenus}
Let $\Sigma_g$ be a closed connected orientable surface of genus $g \geq 2$.
Then, for every smooth representation $\varrho\,:\,\pi_1(\Sigma_g) \rar G(\c)$, we have
$$
\HR_i(\Sigma_g,G)_{\varrho}=0\ \ \mbox{if} \ \, i>\dim \mathcal{Z}(G)\ .
$$
In particular, if $G$ is semisimple, then $\HR_i(X,G)_{\varrho}$ vanishes for all $i>0$.
\ecor
\bproof
Let $\pi:=\pi_1(\Sigma_g)$. As pointed out in the proof of \cite[Theorem 3]{Go}, the natural map $\Rep_G(\pi) \rar \Rep_{G/\mathcal{Z}(G)}(\pi)$ is a principal $\Rep_{\mathcal{Z}(G)}(\pi)$-bundle. By \cite[Lemma 1]{Go}, the dimension of every connected component of $\Rep_{G/\mathcal{Z}(G)}(\pi)$ is $(2g-1)\dim[G/\mathcal{Z}(G)]$. Clearly, $\Rep_{\mathcal{Z}(G)}(\pi)\,\cong\, \mathcal{Z}(G)^{2g}$. Hence, the dimension of every connected component of $\Rep_G(\pi)$ is $\,(2g-1)\dim G \,+\,\dim \mathcal{Z}(G) =
(1-\chi_{\rm top}(\Sigma_g))\dim G \,+\, \dim \mathcal{Z}(G)\,$.
The desired result is now immediate from Theorem \ref{vrephomsurfaces}$(ii)$.
\eproof
An important consequence of Theorem \ref{vrephomsurfaces} is the existence\footnote{We remark that the 
$K$-theoretic virtual fundamental class exists for any quasi-compact derived scheme, which is quasi-smooth in the sense that its cotangent complex is perfect with homology vanishing in degree $ i > 1 $.} 
of a {\it $K$-theoretic virtual fundamental class} of the derived  scheme $ \DRep_G(\Sigma_g) $, which is formally defined by 
(see \cite[Sect. 3.2]{CK})
\begin{equation}
\la{vircl}
[\DRep_G(\Sigma_g)]^{\mathrm{vir}}_K\ :=\ \sum_{i \ge 0}\, (-1)^i\, [\HR_i(\Sigma_g,G)] \ .
\end{equation}
where $\, [\HR_i(\Sigma_g,G)]\,$ is the class of the $ \cO[\Rep_G(\pi_1(\Sigma_g))]$-module
$ \HR_i(\Sigma_g,G) $ in the Grothendieck group $ K_0(\Rep_G[\pi_1(\Sigma_g)])$.
Our interpretation of representation homology in terms of functor homology
leads to a nice Tor-formula for  \eqref{vircl} similar to the well-known Serre formula in classical intersection theory.
\bcor
\label{vfclass0}
The class \eqref{vircl} is well-defined in $ K_0(\Rep_G[\pi_1(\Sigma_g)]) $ and given by the  formula
$$
[\DRep_G(\Sigma_g)]^{\mathrm{vir}}_K \,=\, \sum_{i=0}^{\dim G} (-1)^i\, [\Tor^{\mathfrak{G}}_i(k[\pi_1(\Sigma_g)],\,\cO(G))]\ .
$$
\ecor
\bproof
We have
\begin{eqnarray*}
[\DRep_G(\Sigma_g)]^{\mathrm{vir}}_K &=& \sum_{i=0}^{\dim G}(-1)^i [\HR_i(\Sigma_g,G)] \qquad \text{ [by Theorem \ref{vrephomsurfaces} $(i)$]}\\
	&=& \sum_{i=0}^{\dim G}(-1)^i [\HR_i(\pi_1(\Sigma_g),G)] \quad \text{ [since $\Sigma_g$ is aspherical]}\\
	&=& \sum_{i=0}^{\dim G}(-1)^i \Tor^{\mathfrak{G}}_i(k[\pi_1(\Sigma_g)],\cO(G))\quad \text{ [by \cite[Corollary 4.3]{BRY}].}
\end{eqnarray*}
This proves the assertion.
\eproof
\begin{remark}
The claim of Corollary~\ref{vfclass0} holds for an arbitrary affine algebraic group $G$ defined
over a field $k$ of characteristic $0$.
\end{remark}

\vspace{1ex}

The next observation shows that there is an interesting dichotomy between $ \GL_n $ and semisimple algebraic groups with regard to representations of the surface groups.
\bcor
\label{vfclass} $(a)$ If $\,G=\GL_n$, $\,n\ge 1$, then for all $\, g\ge 1\,$,
$$
[\DRep_{G}(\Sigma_g)]^{\mathrm{vir}}_K = 0\ \ \mbox{in} \ \ K_0(\Rep_G[\pi_1(\Sigma_g)])_{\Q} \ ,
$$
$(b)$ If $G$ is a semisimple  group, then for every $\, g\ge 2\,$,
$$
[\DRep_{G}(\Sigma_g)]^{\mathrm{vir}}_K \not= 0\ \ \mbox{in} \ \ K_0(\Rep_G[\pi_1(\Sigma_g)])_{\Q} \ .
$$
\ecor
\bproof $(a)$
As shown in \cite[Section 7.1]{BRY}, for $\,G=\GL_n$, $n\ge 1$,
the derived representation scheme $\DRep_G(\Sigma_g)$ can be represented by the the global Koszul complex
$$\big(k[X_1,\dots,X_g,Y_1,\ldots,Y_g,\Theta;\det(X_1)^{-1},\ldots,\det(X_g)^{-1},\det(Y_1)^{-1},\ldots,\det(Y_g)^{-1}]\,,\,d\big)\ .$$ 
Here, each of the generic square matrices $X_i,1\leq i \leq g$, $Y_j,1 \leq j\leq g$ and $\Theta$ carries $n^2$ commuting variables. The differential $d$ is given by 
$$d\Theta\,=\,\prod_{i=1}^gX_iY_iX_i^{-1}Y_i^{-1}-\id_n\ . $$
It therefore, suffices to show that the $K$-theoretic virtual fundamental class of an affine
DG-scheme represented by a global Koszul complex $(A \otimes_k \Lambda^{\ast}V, d)$ vanishes in rationalized $K$-theory, where $A$ is any smooth $k$-algebra and where $V$ is a finite-dimensional $k$-vector space.

For any homologically graded sheaf $\mathcal F:=\{\mathcal F_i\}_{0 \leq i \leq n}$ on a scheme $Y$, let $[\mathcal F]$ denote the class 
$$\sum_{i \in \Z} (-1)^i\,[\mathcal F_i]\,\in\,K_0(Y)\ .$$ 
For the rest of this proof, we follow conventions from \cite{CK}: in particular, we switch from homological to cohomological grading by inverting degrees. Let $X=\mathrm{Spec}(A \otimes_k \Lambda^{\ast}V, d)$. Then, $X^0=\mathrm{Spec}(A)$ and $\pi_0(X)=\mathrm{Spec}(A/I)$, where $I$ denotes the ideal generated by $d(V)$. Consider the sheaf of DG algebras $\mathbb{O}:=\cO_X \otimes_{\cO_{X^0}} \cO_{\pi_0(X)}$ on $\pi_0(X)$. In our case, $\mathbb{O}=A/I \otimes_k \Lambda^{\ast}V$ (with trivial differential).

Let $N$ denote the normal cone of $\pi_0(X)$ in $X^0$, with $i\,:\,\pi_0(X) \rar N$ and $p\,:\,N \rar \pi_0(X)$ denoting the inclusion of the zero section and the canonical projection respectively. There is a natural action of $\mathbb{G}_m$ on $N$, whose locus of fixed points is $\pi_0(X)$. The maps $i$ and $p$ are $\mathbb{G}_m$-equivariant. By \cite[Lemma 3.3.7]{CK},
$$i_\ast[\H^{\ast}(\cO_X)]\,=\,[\H^{\ast}(p^{\ast}\mathbb{O})] \,\in\, K_0^{\mathbb{G}_m}(N)\ .$$
Since $p^{\ast}\mathbb{O}=\cO_N \otimes_k \Lambda^{\ast}V$,
$$[\H^{\ast}(p^{\ast}\mathbb{O},\delta)]\,=\,[p^{\ast}\mathbb{O}] \,=\,[\cO_N \otimes_k \Lambda^{\ast}V]=0 \ .$$
Thus, $i_\ast[\H^{\ast}(\cO_X)]=0$ in $K^{\mathbb{G}_m}_0(N)$. As in \cite[Sec. 3.3]{CK}, we argue (using the localization result \cite[Thm. 2.1]{Tho} of Thomason) that $i_\ast\,:\,K^{\mathbb{G}_m}_0(\pi_0(X)) \rar K^{\mathbb{G}_m}_0(N)$ is injective. It follows that $[\H^{\ast}(\cO_X)]=0$ in $K^{\mathbb{G}_m}_0(\pi_0(X))$. Since $K^{\mathbb{G}_m}_0(\pi_0(X))\,\cong\,K_0(\pi_0(X)) \otimes \c[\mu,\mu^{-1}]$, $[\H^{\ast}(\cO_X)]=0$ in $K_0(\pi_0(X)) \otimes \mathbb{Q}$ as well.

$(b)$ Let $j\,:\,U \hookrightarrow Y$ denote the natural inclusion, where $U$ denotes the smooth locus of $Y:=\Rep_G[\pi_1(\Sigma_g)]$. It is well known that the restriction map $j^{\ast}\,:\,K_0(Y) \otimes \mathbb{Q} \rar K_0(U) \otimes \mathbb{Q}$ is surjective. By Corollary \ref{vrephomhighergenus}, we then conclude
$\, j^{\ast}[\DRep_G(\Sigma_g)]^{\mathrm{vir}}_K\,=\, [\cO(U)] \neq 0\,$.
\eproof

\begin{example}
\la{conj1tori}
 We now show that Conjecture \ref{Conj1} holds for $k=\c$ and $G= (\c^*)^r$. In this case, $\cO(G)=\c[t_1^{\pm 1},\ldots,t_r^{\pm 1}]$. By \cite[Sect. 7.1.2]{BRY}, we have
$$\HR_{\ast}(\Sigma_g, G)\,\cong\,\Tor^{\cO(G)}_{\ast}(k,\cO(G^{2g}))\,, $$
where $k$ is equipped with the $\cO(G)$-module structure coming from the canonical augmentation $\cO(G) \rar k$, and $\cO(G^{2g})$ is equipped with the $\cO(G)$-module structure coming from the compositire map
$$ \cO(G) \rar k \rar \cO(G^{2g})\ .$$
Here, the first arrow is the canonical augmentation on $\cO(G)$ and the second arrow is the unit. Since $k$ has a free $\cO(G)$-module resolution given by the Koszul complex
$$\cO(G) \otimes_k \Lambda^{\ast}(\tau_1,\ldots,\tau_r)\,,\qquad d\tau_i\,=\,t_i-1\,, $$
we have
$$ \HR_{\ast}(\Sigma_g,G)\,\cong\,\cO(G^{2g}) \otimes_k \Lambda^{\ast}(\bar{\tau}_1,\ldots,\bar{\tau}_r)\ . $$
It follows that in this case,
$$ \HR_i(\Sigma_g,G)=0 \text{ for } i>r=\dim \mathcal{Z}(G)                                                                                                                                                                                                                \ .$$

We now turn to proving Theorem \ref{vrephomsurfaces}. We begin by studying the corresponding linearized representation homologies case by case, starting with the torus.
\end{example}

\subsubsection{The 2-torus} 
\la{2-torus} 
The $2$-dimensional torus has a simplicial group model given by (see \cite[Section 7.1.1]{BRY})
$$ \lgr{(\mathbb T^2)} \,\simeq\, \mathrm{hocolim}(1 \leftarrow \mathbb F_1 \stackrel{\alpha}{\rar} \mathbb F_2)\,,$$
the map $\alpha$ is defined on generators by $c \mapsto [a,b]:=aba^{-1}b^{-1}$ and the homotopy colimit is taken in the category  of simplicial groups.
It follows from \cite[Lemma 4.2]{BRY}  that
\begin{equation} 
\label{torusdrep} 
\DRep_G(\mathbb{T}^2) \,\simeq\, \mathrm{hocolim}[k \leftarrow \cO(G) \stackrel{\alpha_\ast}{\rar} \cO(G \times G)]\,,\end{equation}
where the homotopy colimit is taken in the category of simplicial commutative algebras, and the map $\alpha_\ast$ is induced by $\alpha$ (explicitly, $\alpha_\ast(f)(x,y) =f([x,y])$ for all $f \in \cO(G)$ and $x,y \in G$).
A representation $\varrho\,\in\,\Rep_G(\mathbb{T}^2)$ is equivalent to a pair $(x,y)\,\in\,G \times G$ such that $[x,y]=\id_G$. Equipping $\cO(G \times G)$ with the augmentation determined by $(x,y)$, we obtain an augmentation on $\DRep_G(\mathbb{T}^2)$, making \eqref{torusdrep} an isomorphism in the homotopy category $\Ho(\scAlg_{k/k})$ of {\it augmented} simplicial commutative $k$-algebras.
Since the derived linearization functor $\L\cQ$ commutes with homotopy pushouts, and by \cite[Prop. B.3]{BRY} (which applies because $\cO(G)$ and $\cO(G \times G)$ are smooth commutative algebras),
\begin{equation} 
\la{linreptorus} \T^{\ast}_{\varrho}\DRep_G(\mathbb{T}^2) \,\simeq\, \cn[ \begin{diagram}[small]\g^{\ast} & \rTo^{\cQ(\alpha_\ast)}& \g^{\ast} \oplus \g^{\ast} 
\end{diagram}]\,, 
\end{equation}
where $\cn$ is the mapping is in the category $\Com_k$ of complexes of  $k$-vector spaces and the cotangent space to $G \times G$ at $(x,y)$ is identified with the cotangent space to $G \times G$ at the identity (i.e, $\g^{\ast} \oplus \g^{\ast}$) via right translation.
The map $Q(\alpha_\ast)$ is dual to the map
$$\g \oplus \g \rar \g\,,\,\,\,\, (u,v) \mapsto (\id-\mathrm{Ad}(y))\cdot u - (\id-\mathrm{Ad}(x))\cdot v \ .$$
It follows from \eqref{linreptorus} that $\H_i[\T^{\ast}_{\varrho}(\DRep_G(\mathbb{T}^2))] =0$ for $i>1$, while
$$
\H_0[\T^{\ast}_{\varrho}\DRep_G(\mathbb{T}^2)] \cong  \Coker[\cQ(\alpha_\ast)]
\quad \mbox{and} \quad
\H_1[\T^{\ast}_{\varrho}\DRep_G(\mathbb{T}^2)] \cong  \Ker[\cQ(\alpha_\ast)] \ .
$$
Hence, we have
\blemma \la{h1ttorus}
For any representation $\varrho\,:\,\pi_1(\mathbb{T}^2) \rar G({k})$,
$\dim_{{k}} \H_1[\T^{\ast}_{\varrho}(\DRep_G(\mathbb{T}^2))] \leq \dim(G)$, with equality corresponding to the trivial representation.
\elemma
\subsubsection{Riemann surfaces} Let $\Sigma_g$ be a closed connected orientable surface of genus $g \geq 1$. $\Sigma_g$ has a simplicial group model (see \cite[Section 7.1.1]{BRY})
$$\lgr{(\Sigma_g)} \,\simeq\,\mathrm{hocolim}(1 \leftarrow \mathbb{F}_1 \stackrel{\alpha^g}{\rar} \mathbb{F}_{2g} )\,,  $$
where $\mathbb{F}_1$ and $\mathbb{F}_{2g}$ denote free groups generated by $c$ and $\{a_1,b_1,\ldots,a_g,b_g\}$ respectively, and $\alpha^g$ is given by
$$
\alpha^g(c) = [a_1,b_1][a_2,b_2]\ldots [a_g,b_g]\ .$$
This implies that
\begin{equation} \la{drepsigmag}\DRep_G(\Sigma_g)\,\simeq\, \mathrm{hocolim}[k \leftarrow \cO(G) \stackrel{\alpha^g_{\ast}}{\longrightarrow} \cO(G^{2g})] \ .\end{equation}
It then follows from \cite[Lemma 4.2]{BRY} that at any representation $\varrho$ of $\pi_{1}(\Sigma_g)$,
$$
\T^{\ast}_{\varrho}\DRep_G(\Sigma_g) \,\cong\, \cn[\begin{diagram}[small] \g^{\ast} & \rTo^{\cQ(\alpha^g_{\ast})}&(\g^{\ast})^{2g} \end{diagram}]\,,$$
where the cotangent space to $G^{2g}$ at the tuple $(x_1,y_1,\ldots,x_g,y_g)$ corresponding to the representation $\varrho$ is identified with the cotangent space to $G^{2g}$ at the identity via right translation. As in the case of the torus,
$\H_i[\T^{\ast}_{\varrho}(\DRep_G(\Sigma)_g))]
=0$ for $i>1$ and
$$
\H_0[\T^{\ast}_{\varrho}(\DRep_G(\Sigma_g))]\,\cong\,\Coker[\cQ(\alpha^g_{\ast})]\,,\quad
\H_1[\T^{\ast}_{\varrho}(\DRep_G(\Sigma_g))]\,\cong\,\Ker[\cQ(\alpha^g_{\ast})]\ . 
$$
The argument in the torus case works with trivial modificatios to give:
\blemma \la{h1tsigmag}
For any $\varrho\,:\,\pi_1(\Sigma_g) \rar G({k})$,
 $\dim_{{k}}\H_1[\T^{\ast}_{\varrho}(\DRep_G(\Sigma_g))] \leq \dim G $, with equality corresponding to the trivial representation.
 \elemma
By \cite[Section 7.1.3]{BRY}, a similar result holds for non-orientable surfaces.

\subsubsection{Proof of Theorem \ref{vrephomsurfaces}}
Part $(i)$ of Theorem \ref{vrephomsurfaces} follows immediately from Theorem~\ref{bounded_homology} and Lemmas \ref{h1ttorus} and \ref{h1tsigmag}. For Part $(ii)$, note that for $\Sigma_g$ as above, $\H_i[\T^{\ast}_{\varrho}\DRep_G(\Sigma_g)]=0$ for $i>1$, and $\chi[\T^{\ast}_{\varrho}\DRep_G(\Sigma_g)]=(1-
\chi_{\rm top}(\Sigma_g))\dim G\,$ by Proposition \ref{Euler}. Since $\dim\H_0[\T^{\ast}_{\varrho}\DRep_G(\Sigma_g)]=\dim\Rep_G[\pi_1(\Sigma_g)]$ for smooth $\varrho$,
$$
\dim\H_0[\T^{\ast}_{\varrho}\DRep_G(\Sigma_g)]=\dim\Rep_G[\pi_1(\Sigma_g)]\,-\,(1-\chi_{\rm top}(\Sigma_g)) \dim G\ .
$$
The desired statement then follows from Theorem~\ref{bounded_homology}.

\subsection{Link complements}
By a link $L$ in $\R^3$ we mean a smooth oriented embedding of a disjoint union $\mathbb{S}^1 \sqcup \ldots \sqcup \mathbb{S}^1$ of (a finite number of)
copies of $\mathbb{S}^1$ in $\R^3$. The link complement $\R^3 \setminus L$ is defined to be the complement of (an open tubular neighborhood of) the image of $L$. By a well-known theorem of J. W. Alexander, every link in $\R^3$ can be obtained geometrically as the closure of some (not unique) braid $\beta$ in $\R^3$ (we write $L=\hat{\beta}$ to indicate this relation). Let $B_n$ denote the group of braids on $n$ strands in $\R^3$. The group $B_n$ is generated by $n-1$ ``flips" $\sigma_1,\ldots,\sigma_{n-1}$ subject to the relations
$$ \sigma_i \sigma_j=\sigma_j\sigma_i \,\,\,\text{(if $|i-j|>1$)} \ ,\quad \sigma_i\sigma_j\sigma_i =\sigma_j\sigma_i\sigma_j \,\,\, \text{(if $|i-j|=1$)} \ .$$
There is a faithful representation (the Artin representation) $B_n$ in $\mathrm{Aut}(\mathbb{F}_n)$. Explicitly, the generator $\sigma_i$ acts of $\mathbb{F}_n$ via
$$x_i \mapsto x_ix_{i+1}x_i^{-1}\,,\,\, x_{i+1} \mapsto x_i \,,\,\, x_j \mapsto x_j \,\,(\text{for}\,j \neq i,i+1)\ . $$
The group $B_n$ naturally surjects onto $S_n$ (the ``flip" $\sigma_i$ is mapped to the transposition $(i, i+1)$ under this surjection). We denote the image of $\beta$ under this surjection by $\bar{\beta}$. A refinement of a classical theorem of Artin and Birman (see \cite[Prop. 7.1]{BRY}) shows that if $L=\hat{\beta}$ for some $\beta \,\in\,B_n$, then the link complement $\R^3 \setminus L$ has the simplicial group model
\begin{equation} \la{sgmodlc} \lgr{(\R^3 \setminus L)} \, \simeq\, \mathrm{hocolim}(\mathbb{F}_n \xleftarrow{(\id,\beta)} \mathbb{F}_n \sqcup \mathbb{F}_n \xrightarrow{(\id,\id)} \mathbb{F}_n )\,,\end{equation}
where $\beta$ acts on $\mathbb{F}_n$ via the Artin representation.  Note that the Artin representation induces an (contravariant) action of $B_n$ on $\Rep_G(\mathbb{F}_n)=G^n$. Let $\varrho\,:\,\pi_1(\R^3\setminus L) \rar G(k)$ be any representation. Such a representation corresponds to an element of $G^n$ fixed by $\beta$. In other words, the map $\beta_\ast\,:\,\cO(G^n) \rar \cO(G^n)$ preserves the corresponding augmentation of $\cO(G^n)$. The same argument 
as in Section~\ref{2-torus} (see the discussion before formula \eqref{linreptorus}) yields
\blemma
\la{cotcplx}
$\ \T^{\ast}_{\varrho}\DRep_G(\R^3 \setminus L)\,\simeq\,\cn[(\g^{\ast})^{\oplus n} \xrightarrow{\id-\cQ(\beta_\ast)}
(\g^{\ast})^{\oplus n}]\,$.
\elemma
%
%
%
%
%
%
%

Now, let $L=\hat{\beta}$ be a link in $\R^3$, where $\beta\,\in\,B_n$ and let $X:=\R^3 \setminus L$. Denote by $M$ the maximum dimension of $\T_{\varrho}\Rep_{G({k})}[\pi_1(X)]$ as $\varrho$ runs over all representations $\pi_1(X) \rar G({k})$. Denote the number of cycles in a cycle decomposition of a permutation $\sigma$ by $n(\sigma)$. If $ L=\hat{\beta} $, then the number $ n(\beta) $
equals the number of components of the link $ L $.
\bthm \la{linkcompl} With the above notation, we have

\vspace*{1ex}

$(i)$ $\HR_i(X,G)=0$ for all $i>M$. In particular, $\HR_i(X,G)=0$ for all $i>n \cdot \dim(G)$.

\vspace*{1ex}

$(ii)$ For smooth $\varrho\,:\,\pi_1(X) \rar G({k})$,
$$\HR_i(X,G)_{\varrho}=0 \text{ for } i>\dim_{\varrho}\Rep_G[\pi_1(X)]\,,$$
where $\dim_{\varrho}\Rep_G[\pi_1(X)]$ denotes the local dimension of the irreducible component of $\Rep_G(\pi_1(X))$ containing $\varrho$.

\vspace*{1ex}

$(iii)$ For smooth $\varrho$ in the irreducible component of $\Rep_G[\pi_1(X)]$ containing the trivial representation, $\HR_i(X,G)_{\varrho} =0$ for $i>\dim(G)\cdot n(\beta)$.
\ethm

\remark\ The second statement in Theorem \ref{linkcompl}, part $(i)$ may also be shown using the description of $\HR_{\ast}(X,G)$ in terms of Hochschild homology (see \cite[Theorem 7.1]{BRY}). The remaining results in this section do not appear to be easily recovered from that point of view.

\bproof
By Lemma \ref{cotcplx}, $\H_i[\T^{\ast}_{\varrho}\DRep_G(X)]=0$ for $i>1$. Further, $\chi[\T^{\ast}_{\varrho}\DRep_G(X)]=0$. Thus,
$$\dim_{{k}}\H_1[\T^{\ast}_{\varrho}\DRep_G(X)] = \dim_{{k}}\H_0[\T^{\ast}_{\varrho}\DRep_G(X)]=\dim_{{k}} \T_{\varrho}\Rep_{G({k})}(\pi_1(X))\ .$$
 Again, by Lemma \ref{cotcplx}, this dimension is bounded above by $n \cdot \dim(G)$ for all $\varrho$. Part $(i)$ of the desired theorem thus follows from Theorem~\ref{bounded_homology}. For smooth $\varrho$, we have
 $$\dim_{{k}}\H_1[\T^{\ast}_{\varrho}\DRep_G(X)] = \dim_{{k}} \T_{\varrho}\Rep_{G({k})}[\pi_1(X)] = \dim[\Rep_G(\pi_1(X))]\ .$$
 Part $(ii)$ thus follows from Theorem~\ref{bounded_homology}.
Part $(iii)$ follows from the fact that at the trivial representation, $\cQ(\beta_{\ast})$ is indeed the map given by permuting the $n$ copies of $\g^{\ast}$ by $\bar{\beta}$. $\pi_1[\T^{\ast}_{\varrho}\DRep_G(X)]$ is therefore, the fixed point space of $\bar{\beta}$. The dimension of this space is easily seen to be $\dim(G) \cdot n(\beta)$.
\eproof
The following example illustrates that the number $M$ in Theorem \ref{linkcompl}, Part $(i)$ may be considerably smaller than $n \cdot \dim(G)$, where $n$ is the minimum number of strands such that there exists a $\beta \in B_n$ with $L=\hat{\beta}$. Indeed, let $L$ be a $(p,q)$-torus knot, where $2<p<q$ and $p,q$ relatively prime. Then, it is well known that the minimal presentation of $L$ as a braid closure is given by $L=\hat{\beta}$, where $\beta = (\sigma_1\ldots\sigma_{p-1})^q \,\in\,B_p$. Nevertheless,
\bcor \la{pqtorus}
Let $X$ be the complement of a $(p,q)$-torus knot in $\R^3$. Then

\vspace*{1ex}
$(i)$ $\HR_i(X,G)=0$ for $i>2 \cdot \dim(G)$.

\vspace{1ex}

$(ii)$ For smooth $\varrho$ in the irreducible component of $\Rep_G(\pi_1(X))$ containing the trivial representation, $\HR_i(X,G)=0$ for $i>\dim(G)$.
\ecor
\bproof
By Theorem \ref{linkcompl}, part $(i)$ follows once we verify that for any representation $\varrho\,:\,\pi_1(X) \rar G({k})$, $\dim_{{k}} \T_{\varrho}\Rep_G(\pi_1(X)) \leq 2 \cdot \dim(G)$. Since $L$ is a knot, $\pi_1(X)\,\cong\,\pi_1(\mathbb{S}^3 \setminus L) \,\cong\, \langle x, y\,|\, x^p=y^q \rangle$. Thus, any representation $\varrho\,:\,\pi_1(X) \rar G({k})$ is determined by a pair $(X_0,Y_0)$ in $G({k})^{\times 2}$ such that $X_0^p=Y_0^q$. The tangent space to $\Rep_G(\pi_1(X))$ at such a representation is easily seen to be the space
\begin{equation} \la{tsptorus} \{(u,v)\,\in\,\g \oplus \g\,|\, \sum_{i=0}^{p-1} \mathrm{Ad}(X_0^i)(u)\,=\,\sum_{i=0}^{q-1} \mathrm{Ad}(Y_0^i)(v) \}\ . \end{equation}
Clearly, for all $\varrho\,:\,\pi_1(X) \rar G({k})$, $\dim_{{k}} \T_{\varrho}\Rep_G(\pi_1(X)) \leq 2 \cdot \dim(G)$. This completes the required verification for part $(i)$. For part $(ii)$, we note that the tangent space to $\Rep_G(\pi_1(X))$ at the trivial representation is $\{(u,v)\,\in\,\g \oplus \gamma\,|\,p \cdot u =q \cdot v\}$. This is isomorphic to $\g$ as a vector space. It follows that the corresponding irreducible component of $\Rep_G(\pi_1(X))$ has dimension at most $\dim(G)$. Part $(ii)$ then follows from Theorem \ref{linkcompl}, part $(ii)$.
\eproof

\subsection{Lens spaces}
Given an integer $p>1$ and $ m>1 $ integers $\,q_1, q_2,\ldots,q_m\,$, which are relatively prime to $p$, the lens space 
$ L=L_p(q_1,\ldots,q_m)\,$ is defined to be the quotient $\mathbb{S}^{2m-1}/\Z_p$ of the unit sphere $\mathbb{S}^{2m-1} \subset \c^m$ modulo the (free) action of $\Z_p$ generated by the rotation $\gamma(z_1,\ldots,z_m)=(e^{2\pi i q_1/p}z_1,\ldots,e^{2 \pi i q_m/p}z_m)$. By definition, $ L $ is a compact connected manifold with universal cover
$ \mathbb{S}^{2m-1} $ and fundamental group $\Z_p$.

Assume, for simplicity, that $ k = \c $ and $G$ is either $\GL_n $ or $ \SL_n $ with $n \geq 1$. Then
\blemma \label{cothomlens}
For any representation $\varrho\,:\,\Z_p \rar G(\c)$, $\H_i[\T^{\ast}_{\varrho}\DRep_G(L)]=0$ for all $\,i \neq 0,\, 2m-2\,$. Moreover, $\H_{2m-2}[\T^{\ast}_{\varrho}\DRep_G(L)] \neq 0$.
\elemma
\bproof
Following \cite[Example 2.43]{Hat}, we construct a $\Z_p$-equivariant cell structure on $\tilde{L}:=\mathbb{S}^{2m-1}$ that descends to a cell structure on $L$ as follows. Subdivide the unit circle $C$ in the $m$-th $\c$-factor of $\c^m$ by taking the points $e^{2 \pi i j /p}\,\in\,C$ as vertices. Joining the $j$-th vertex of $C$ to the unit sphere $\mathbb{S}^{2m-3} \subset \c^{m-1}$ by arcs of great circles in $\mathbb{S}^{2m-1}$ yields a $(2m-2)$-dimensional cell $B_j^{2m-2}$ bounded by $\mathbb{S}^{2m-3}$. Specifically, $B_j^{2m-2}$ consists of the points $\cos \theta(0,\ldots,0,e^{2 \pi i j/p})+\sin\theta(z_1,\ldots,z_{m-1},0)$ for $0 \leq \theta \leq \pi/2$. Similarly, joining the $j$-th edge of $C$ to $\mathbb{S}^{2m-3}$ gives a $(2m-1)$-dimensional cell $B_j^{2m-1}$ bounded by $B_j^{2m-2}$ and $B_{j+1}^{2m-2}$, the subscripts being taken mod $p$. The rotation $\gamma$ takes $\mathbb{S}^{2m-3}$ to itself, while permuting the $B_j^{2m-2}$'s and $B_j^{2m-1}$'s. It is easy to check that $\gamma^{l_m}$ takes each $B_j^{2m-2}$ and $B_j^{2m-1}$ to the next one, where $q_m\,l_m \equiv 1 \,(\text{mod}\, p)$. The lower dimensional cells in the desired structure are obtained by repeating the above construction on $\mathbb{S}^{2m-3} \subset \c^{m-1}$.

Thus, as $\Z[\Z_p]$-modules, we have $ C_k(\tilde{L}) \,\cong\,\Z[\Z_p]$ for $0 \leq k \leq 2m-1$, with differential given by
$$
d x\,=\,\begin{cases}
 (\sum_{r=0}^{p-1} \gamma^r) \cdot x & \text{if } x \in C_{2k}\,,\ 1 \leq k \leq m-1\\*[1ex]
 (\gamma^{l_k}-1)\cdot x & \text{if } x \in {C}_{2k-1}\,,\ 1 \leq k \leq m
\end{cases} 
$$
By Theorem \ref{LinRep}, $\H_\ast[\T^{\ast}_{\varrho}\DRep_G(L)]$ is the homology of the complex
$$
\ldots \,\to\, {C}_k(\tilde{L}) \otimes_{\Z[\Z_p]} \g^{\ast}\, \to \, {C}_{k-1}(\tilde{L}) \otimes_{\Z[\Z_p]} \g^{\ast} \,\to\,\ldots \,\to\, {C}_1(\tilde{L}) \otimes_{\Z[\Z_p]} \g^{\ast} \,\to\, 0 \ , 
$$
where $\g^{\ast}$ is the $\Z_p$-representation dual to $\mathrm{Ad} \, \varrho$ and ${C}_k(\tilde{L}) \otimes_{\Z[\Z_p]} \g^{\ast}$ has degree $k-1$ for $ k \geq 1 $. The graded linear dual of this complex is the {\it cohomologically}\, graded complex
\begin{equation}\la{cxg}
0 \to \g \xrightarrow{\delta^0}  \g \xrightarrow{\delta^1} \g \xrightarrow{\delta^2} \ldots 
\xrightarrow{\delta^{2m-4}}  \g \xrightarrow{\delta^{2m-3}} \g \to 0 
\end{equation}
with nonzero terms concentrated in cohomological degrees $\,0, 1,\ldots, 2m-2\,$ and the differentials given by
$$
\delta^{2k}=\sum_{r=0}^{p-1} X^r \quad \mbox{for}\quad 0 \leq k \leq m-2\ ,
$$
$$
\delta^{2k-1}= X^{l_{k+1}}-1 \quad \mbox{for}\quad 1 \leq k \leq m-1\,, 
$$
where $X:=\mathrm{Ad}\, \varrho(\gamma)\,$. Up to isomorphism, the cohomology of \eqref{cxg} depends  only on the conjugacy class of $\varrho(\gamma)$. Since $\varrho(\gamma)$ is an element of order $ p $ in $ G(\c)$, its minimal polynomial divides $ x^p - 1 $ and hence factors into distinct linear factors over $ \c $. It follows that $\varrho(\gamma)$ is conjugate to a diagonal matrix, with eigenvalues being $p$-th roots of unity (possibly non-primitive).
From the structure of differentials $ \delta^j $ it is easy to see that the complex \eqref{cxg} is acyclic in all degrees,
except possibly in degree $ 0 $ and $ 2m-2 $. The kernel of $ \delta^{2m-3} $ is nontrivial (it contains, 
in particular, the diagonal matrices in $ \g $), hence $ \im(\delta^{2m-3}) \not= \g $, which means that the 
cohomology  in degree $ 2m-2 $ is always nonzero. This proves the desired lemma.
\eproof
\bthm \la{rephomlens} Let $G$ be either $\GL_n(\c) $ or $ \SL_n(\c) $ for $n \geq 1$. Then
$$\HR_i(L_p(q_1,\ldots,q_m), G) = \begin{cases}
0 &  \text{iff} \ \, i \not\equiv 0\ \mathrm{mod}\,(2m-2)\\*[1ex]
\neq 0 & \text{iff}\ \, i \equiv 0\ \mathrm{mod}\,(2m-2)
\end{cases}  
$$
\ethm
%
%
\bproof
To simplify the notation we set  $L:=L_p(q_1,\ldots,q_m)$. Choose a representative $B$ of $\mathcal{A}_{G}(L)$ in $\cDGA_{\c}^+$ that has finitely many generators in each degree. Such a representative exists since $L$ has a finite cellular model, whence $\mathcal{A}_{G}(L)$ is of quasi-finite type. Let $\varrho\,:\,\Z_p \rar G(\c)$ be a representation.
By Proposition \ref{Quillen_specseq}, one has a convergent spectral sequence
$$
E^2_{p,q}:=\H_{p+q}[\Sym^p(\mathrm{T}^{\ast}_{\varrho}\DRep_G(L))] \implies \widehat{\HR}_{p+q}(L,G)_{\varrho}\,, $$
where $\widehat{\HR}_{\ast}(L,G)_{\varrho}$ stands for the completion of the local graded ring $\HR_{\ast}(L,G)_{\varrho}$ with respect to the filtration by powers of the maximal ideal of $\HR_0(L,G)$ defined by $\varrho$. By Lemma \ref{cothomlens}, the $ E^2_{p,q}$ terms of the above spectral sequence are nonzero iff $2m-2$ divides $p+q$. It follows that $ \widehat{\HR}_{i}(L,G)_{\varrho} \neq 0$ iff $2m-2 \,|\,i$.
By Krull's Intersection Theorem, $\HR_i(L,G)_{\varrho} \neq 0$ iff $2m-2\,|\,i$. This proves the desired result.
\eproof
\remark\ Let $\{q_i\,,\,i \in \mathbb{N}\}$ be an infinite sequence of integers coprime to $p$. Then, one has the infinite dimensional lens space $L_{\infty}:=L(q_1,q_2,\ldots)$. Note that $L_{\infty}$ is the (homotopy) direct limit of the sequence of lens spaces $L(q_1,\ldots,q_m)\,,\,m \in \mathbb{N}$. It follows from this observation and Theorem \ref{rephomlens} that $\HR_i(L_{\infty},G)$ vaishes for $i>0$. This agrees with Theorem \ref{vfgroups}, given that $L_{\infty}$ is a $K(\Z_p,1)$ space.

\appendix 

\section{Relation to derived algebraic geometry}

In this Appendix, we discuss the relation of the derived representation scheme $\DRep_G(X)$ to two basic examples of moduli spaces in derived algebraic geometry: the derived moduli space  $  \boldsymbol{R} \Loc_G(X,*) $ of $G$-local systems on a \emph{pointed} space $(X,*)$ studied in  \cite{K}, and the derived stack of flat $G$-bundles $ \textbf{Map}(X,BG) $ on an \emph{unpointed} space $X$ introduced in \cite[Corollary 2.6]{PTVV13} (see also \cite[Definition 2.2.6.2]{TV08} for $G = \GL_n$). Specifically, we will construct below a pointed mapping stack $ \textbf{Map}((X,*),(BG,*)) $, which is equivalent to our derived representation scheme $\DRep_G(X)$ and 
fits in the homotopy cartesian diagram\footnote{We thank the referee for suggesting us 
this relation.}
\begin{equation} \label{mapping_spaces_pullback}
\begin{diagram}[small]
\textbf{Map}((X,*),(BG,*))  & \rTo &  \textbf{Map}(X,BG)  \\
 \dTo & & \dTo\\
\ast = {\bm R}\Spec(k) & \rTo &  BG = \textbf{Map}(*,BG) 
\end{diagram}
\end{equation}
In this way, we realize $\DRep_G(X) \simeq \textbf{Map}((X,*),(BG,*))$
as a homotopy fiber of the To\"{e}n-Vezzosi derived stack $\,\textbf{Map}(X,BG)\,$ over $ BG $, 
so that $\textbf{Map}(X,BG)$ is, in fact, the homotopy quotient of $\DRep_G(X) $ by the natural $G$-action:
\begin{equation}  
\label{DRep_stack_MapXBG}
\textbf{Map}(X,BG) \, \simeq \, [\,\DRep_G(X)/G \,]\ .
\end{equation}

In a different direction, we will also consider the construction of the derived moduli space $\boldsymbol{R} \Loc_G(X,*) $ of $G$-local systems on $ (X,*) $ introduced in \cite{K} and recall from \cite{BRY} that it is 
equivalent to the derived representation scheme $\DRep_G(X)$. Combined together, these two comparison results 
establish a connection between Kapranov's and To\"{e}n-Vezzosi's constructions which are given a priori in different languages (see Corollary~\ref{LocMap} below). This connection is probably well known to experts in
derived algebraic geometry but we could not find an explicit statement in the literature.

\subsection{To\"{e}n-Vezzosi construction}
\la{TVcon}
Recall that the notion of derived stacks \cite{TV05} can be formalized in terms of simplicial presheaves on the (model) category $\dAff_k$, which is defined as the opposite category to the category of simplicial commutative algebras (i.e., $\dAff_k := (\sCommAlg_k)^{\op}$).
More precisely, a \emph{derived stack} is a simplicial presheaf on $\dAff_k$, i.e. a functor $F : \sCommAlg_k \rightarrow \sset$, that is objectwise fibrant, preserves equivalences and satisfies a certain hyperdescent property (see \cite[Definition 4.6.6]{TV05} and 
\cite[Definition 2.2.2.14]{TV08}). 

In this Appendix, we will only be concerned with derived stacks defined as a global quotient
$[ \, \Spec(A)/G  \, ]$ for some $A \in \sCommAlg_k$ acted on by an affine algebraic group $G$. 
To give a precise meaning to $[ \, \Spec(A) / G  \, ]$, we recall a more elementary description of derived stacks in terms of simplicial derived affine schemes, as in \cite{Pri13}. 
Namely, we start with the derived Yoneda functor, which associates a simplicial presheaf ${\bm R}h_X $ to each derived affine scheme $X \in \dAff_k$, defined by
\begin{equation*}
{\bm R}h_X \, = {\rm Map}_{\dAff_k}(-,X) \, : \, \dAff_k^{\op} \, \rightarrow \, \sset
\end{equation*}
where ${\rm Map}_{\dAff_k}(-,-)$ is a functorial homotopy mapping space in the model category $\dAff_k$.
For $X = \Spec(A)$, i.e. the derived affine scheme opposite to  $A \in \sCommAlg_k$, we will simply denote ${\bm R}h_X$ by ${\bm R}\Spec(A)$.
Now, given a simplicial derived affine scheme $X_{\bullet} \in s(\dAff_k)$, one can take the derived Yoneda functor levelwise, the result is a \emph{bisimplicial} presheaf 
$ {\bm R}h_{X_{\bullet}} \, : \, \dAff^{\op} \, \rightarrow \, s(\sset) $.
This allows us to make the following

\bdf
The \emph{derived stack associated to a simplicial derived affine scheme} $X_{\bullet} \in s(\dAff_k)$ is the stackification of the simplicial presheaf 
\[ |X_{\bullet}|^{\rm pre} \, : \, 
\dAff^{\op} \, \xrightarrow{{\bm R}h_{X_{\bullet}}} \, s(\sset)  \, \xrightarrow{\diag} \sset  \]
We denote it by $|X_{\bullet}|^{\st} := (|X_{\bullet}|^{\rm pre} )^{\st}$.
\edf

Now, given any derived affine scheme $X = \Spec(A)$ with an action of an affine algebraic group $G$, or more precisely, a coaction on $A \in \sCommAlg_k$ by the commutative Hopf algebra $\cO(G)$, 
one can define the nerve of this action in the usual way, which gives the simplicial derived affine scheme $\mathcal{N}(G \ltimes \Spec\, A) \in s(\dAff_k)$. This leads us to the following
\bdf
The quotient stack $[ \, \Spec(A)/G  \, ]$ is the derived stack associated to the 
simplicial derived affine scheme $\mathcal{N}(G \ltimes \Spec(A))$.
\edf

In particular, if $A = k$, so that $\Spec(k)$ is the final object $\ast$, then the quotient stack $[ \, \ast \, / \, G  \, ]$ will be simply denoted as $BG$. In other words, it is the derived stack associated to the simplicial derived affine scheme $B_{\bullet}G$ defined in the usual way by $B_nG = G \times \stackrel{n}{\ldots} \times G$.

Notice that the simplicial derived affine scheme $B_{\bullet}G$ satisfies $B_0G = \Spec(k)$.
Therefore, there is a canonical map ${\bm R}\Spec(k) \rightarrow BG$ of simplicial presheaves.
By choosing the homotopy mapping spaces ${\rm Map}_{\dAff_k}(-,-)$ appropriately, we may assume that 
the simplicial presheaf ${\bm R}\Spec(k)(-) := {\rm Map}_{\dAff_k}(-,\Spec(k))$ takes constant value $\ast$ (the singleton)%
\footnote{For example, in \cite{TV05}, one defines the derived Yoneda embedding ${\bm R}h$ on the model category 
$\cC = \dAff_k$ by ${\bm R}h_Y(X) = \Hom_{\dAff_k}(\Gamma(X),R(Y))$ where $\Gamma : \mathcal{C} \rightarrow \mathcal{C}^{\Delta}$ is a cofibrant replacement functor, and $R : \mathcal{C} \rightarrow \mathcal{C}$ is a fibrant replacement functor. If the fibrant replacement functor $R$ is chosen so that it is the identity on the final object $\Spec(k)$ then this criterion is satisfied.},
so that each $BG(A)$ is canonically pointed. In other words, $BG$ is in fact a presheaf
\[
BG \, : \, \dAff_k^{\op} \, \rightarrow \, \sset_*
\]
with values in pointed simplicial sets. This leads to the following definition

\bdf  
\label{mapping_stack_def}
For any simplicial set $X \in \sset$, define the \emph{mapping stack} $\textbf{Map}(X,BG)$ to be the 
simplicial presheaf 
\begin{equation*}
\textbf{Map}(X,BG) \, : \, \dAff_k^{\op} \, \rightarrow \, \sset \, , \qquad \quad A \, \mapsto \, {\rm map}(X, BG(A)) \, ,
\end{equation*}
where ${\rm map}(X,Y)$ is the simplicial set  ${\rm map}(X,Y)_n := \Hom_{\sset}(X \times \Delta^n, Y)$.

Similarly, for a pointed simplicial set $(X,*) \in \sset_*$, define the \emph{pointed mapping stack} $\textbf{Map}((X,*),(BG,*))$ to be the 
simplicial presheaf 
\begin{equation*}
\textbf{Map}((X,*),(BG,*)) \, : \, \dAff_k^{\op} \, \rightarrow \, \sset \, , \qquad \quad A \, \mapsto \, {\rm map}((X,*), (BG(A),*)) \, ,
\end{equation*}
where ${\rm map}((X,*),(Y,*))$ is the simplicial set  ${\rm map}(X,Y)_n := \Hom_{\sset}((X \times \Delta^n, \ast \times \Delta^n), (Y,*))$.
\edf

The simplicial presheaves $\textbf{Map}(X,BG)$ and $\textbf{Map}((X,*),(BG,*))$ are still derived stack%
\footnote{The simplicial presheaf $\textbf{Map}(X,BG)$ is simply the powering of $X \in \sset$ 
on the fibrant object $BG$ in the simplicial model category $(\dAff_k, W)^{\sim, \tau}$ of derived stacks (see \cite{TV05}), and is therefore fibrant. The fact that the pointed mapping space 
$\textbf{Map}((X,*),(BG,*))$ is a derived stack follows by writing it as a homotopy pullback \eqref{mapping_spaces_pullback}.}
and they clearly sit in a (homotopy) pullback diagram \eqref{mapping_spaces_pullback}.

The main result of this Appendix compares the 
derived stack $\textbf{Map}((X,*),(BG,*))$ with the derived representation scheme,
which is a derived affine scheme $\DRep_G(X) \in \dAff_k$ associated to a reduced simplicial set $X \in \sset_0$.
To construct $\DRep_G(X)$, one first take a simplicial group model $\lgr(X) \in \sGr$ of $X$, given by the Kan loop group function $\lgr(-)$
(see Section \ref{simpgr} for details). 
Then applying the derived functor \eqref{Lrep} to the simplicial group $\lgr(X)$, one obtains a simplicial commutative algebra $\lgr(X)_G$.
The derived rerepsentation scheme $\DRep_G(X)$ is then defined to be the derived affine scheme 
$\DRep_G(X) = \Spec \,( \lgr(X)_G )$.
We then have the following

\bprop  \label{DRep_scheme_MapXBG}
The pointed mapping stack $\mbox{\rm \textbf{Map}}((X,*),(BG,*))$ is equivalent to the derived representation scheme $\DRep_G(X)$ in the sense that there is an equivalence of derived stacks
\[
\mbox{\rm \textbf{Map}}((X,*),(BG,*)) \, \simeq \, {\bm R}h_{ \DRep_G(X) }
\]
\eprop

\bproof
Recall from \cite{BRY} that, although the adjunction 
\begin{equation}  
\label{GG_adj}
(\,\mbox{--}\,)_G \,:\, \sGr\, \rightleftarrows\, \sCommAlg_k \,:\, G(\,\mbox{--}\,)
\end{equation}
is not a Quillen adjunction, the total left derived functor ${\bm L}(-)_G : \Ho(\sGr) \rightarrow \Ho(\sCommAlg_k)$ nonetheless
preserves homotopy pushouts and (homotopy) coproducts.
If we denote by $(-)_G^{\infty} : \sGr^{\infty}  \rightarrow  \sCommAlg_k^{\infty}$ the associated 
$\infty$-functor, then it preserves $\infty$-pushouts and $\infty$-coproducts, and hence arbitrary $\infty$-colimits (see \cite[Theorem 4.2.4.1, Proposition 4.4.2.7]{Lur09}). As a result, it has 
an $\infty$-right adjoint (see \cite[Proposition A.3.7.6, Corollary 5.5.2.9]{Lur09}).
Strictify this $\infty$-right adjoint, so that it is represented by a functor $\tilde{G} : \sCommAlg_k \rightarrow \sGr$ that preserves weak equivalences.
Then this functor is characterized by the adjunction, expressed as homotopy equivalences 
\begin{equation}  \label{LG_G_adj}
{\rm Map}_{\sCommAlg_k}({\bm L}(\Gamma)_G , A) \, \simeq \, 
{\rm Map}_{\sGr}(\Gamma , \tilde{G}(A))
\end{equation}
for any simplicial group $\Gamma \in \sGr$.

Now, consider the Quillen equivalence \eqref{kanl}, and take $\Gamma = \lgr(X)$ for some $X \in \sset_0$, then the right hand side of \eqref{LG_G_adj} is equivalent to the homotopy mapping space
${\rm Map}_{\sset_0}(X , \overline{W} (\tilde{G}(A)))$.
Moreover, the term $\overline{W} (\tilde{G}(A))$ in this mapping space has an alternative description by the following
\blemma
The simplicial set $\overline{W} (\tilde{G}(A))$ is equivalent to the simplicial set $|B_{\bullet}G|^{\rm pre}(\Spec (A))$.
\elemma 
\bproof
For any simplicial group $\Gamma \in \sGr$, we have an equivalence $\diag(B\Gamma) \simeq \overline{W}(\Gamma)$ (see, e.g., \cite{Tho08}). Applying this to $\Gamma = \tilde{G}(A)$,
we see that it suffices to show that there is an equivalence $B_{\bullet}(\tilde{G}(A)) \simeq 
{\bm R}h_{B_{\bullet}G} (\Spec (A))$ of bisimplicial sets.
In fact, we will show that there is a levelwise equivalence, i.e. for each $n \geq 0$, there is an equivalence 
$B_n(\tilde{G}(A)) \simeq 
{\bm R}h_{B_n G} (\Spec (A))$ of simplicial sets.
Indeed, take $\Gamma = \mathbb{F}_n$, the free group of $n$ variables, in the adjunction \eqref{LG_G_adj}, then it takes the form
\begin{equation}  \label{RhG_tildeG_equiv}
{\bm R}h_{G\times \stackrel{n}{\ldots} \times G}(\Spec(A)) := {\rm Map}_{\sCommAlg_k}(\mathcal{O}(G)^{\otimes n} , A) \, \simeq \, 
{\rm Map}_{\sGr}(\mathbb{F}_n , \tilde{G}(A)) \simeq  \tilde{G}(A)^n
\end{equation}

In other words, one obtains an equivalence ${\bm R}h_{B_n G} (\Spec (A)) \simeq B_n(\tilde{G}(A)) $
by taking $\Gamma = \mathbb{F}_n$ in the adjunction \eqref{LG_G_adj}.
Moreover, when $n$ varies, the groups $\{ \mathbb{F}_n \}_{n \geq 0}$ form in fact a cosimplicial group in such a way that $B_{\bullet} H = \Hom_{\Gr}(\mathbb{F}_{\bullet} , H)$ as a simplicial set, for any group $H$.
In this way, the equivalences ${\bm R}h_{B_n G} (\Spec (A)) \simeq B_n(\tilde{G}(A)) $ in 
\eqref{RhG_tildeG_equiv} respect the simplicial structure as $n$ varies,
and hence can be strictified into a levelwise equivalence of bisimplicial sets. This completes the proof of the lemma.
\eproof

Thus, combining this lemma with \eqref{LG_G_adj}, and again taking $\Gamma = \lgr(X)$, one has
\begin{equation}  \label{GXG_XBG_adj}
{\rm Map}_{\sCommAlg_k}(\lgr(X)_G , A) \, \simeq \, 
{\rm Map}_{\sset_0}(X , |B_{\bullet}G|^{\rm pre}(\Spec (A)))
\end{equation}

Now, notice that there is a Quillen adjunction
\begin{equation*}
\iota :\, \sset_0\, \rightleftarrows \,\sset_*\,: E
\end{equation*}
where $\iota$ is the inclusion functor, and $E(X)$ is the subsimplicial set consisting of simplices with all vertices lying on the basepoint $\ast \in X$. 
As a result, the right-hand-side of \eqref{GXG_XBG_adj} can be rewritten as
${\rm Map}_{\sset^*}(X , |B_{\bullet}G|^{\rm pre}(\Spec (A)))$,
the homotopy mapping space taken in the model category of pointed simplicial sets.
Thus, if one defines a pointed mapping (pre)stack $\textbf{Map}((X,*),(|B_{\bullet}G|^{\rm pre},*))$
in the same way as in Definition \ref{mapping_stack_def}, 
then this shows that the right hand side of \eqref{GXG_XBG_adj} can be identified with the value of the derived prestack $\textbf{Map}((X,*),(|B_{\bullet}G|^{\rm pre},*))$ on $\Spec(A)$.
Moreover, the left hand side of \eqref{GXG_XBG_adj} is also, by definition, the value of the 
derived stack ${\bm R} \Spec(\lgr(X)_G)$ on $\Spec(A)$.
This shows that there is an equivalence of derived prestacks
\[
{\bm R} \Spec(\lgr(X)_G) \, \simeq \, \textbf{Map}((X,*),(|B_{\bullet}G|^{\rm pre},*))
\]
which in particular shows that $\textbf{Map}((X,*),(|B_{\bullet}G|^{\rm pre},*))$ is in fact a derived stack.

Now, since the stackification functor is left exact (see \cite[Proposition 4.6.7]{TV05}), one can show that the stackification of 
$\textbf{Map}((X,*),(|B_{\bullet}G|^{\rm pre},*))$ is the pointed mapping stack $\textbf{Map}((X,*),(BG,*))$
in Definition \ref{mapping_stack_def}.
Since we have already seen that the former is a derived stack, the stackification is an equivalence, which therefore finishes the proof.
\eproof

Notice that the homotopy pullback \eqref{mapping_spaces_pullback} 
exhibits the pointed mapping stack $\textbf{Map}((X,*),(BG,*))$ as a principal $G$-bundle over
the unpointed mapping stack $\textbf{Map}(X,BG)$. This means that the pointed mapping space $\textbf{Map}((X,*),(BG,*))$ has a $G$-action whose homotopy quotient coincides with the 
unpointed mapping stack $\textbf{Map}(X,BG)$.
One can check that, under the equivalence in Proposition \ref{DRep_scheme_MapXBG} between
the pointed mapping stack $\textbf{Map}((X,*),(BG,*))$ and the derived representation scheme 
$\DRep_G(X)$, this $G$-action coincides%
\footnote{
Indeed, for any simplicial group $\Gamma$, the $G$-action on the simplicial commutative algebra $\Gamma_G$ is induced by the adjunction \eqref{GG_adj}, so that $G(A)$ acts on $\Hom_{\sCommAlg_k}(\Gamma_G,A) \cong \Hom_{\sGr}(\Gamma,G(A))$ by conjugation on $G(A)$.
Retaining the notation in the proof of Proposition \ref{DRep_scheme_MapXBG},
one can likewise consider an action of $\tilde{G}(A)$ on the spaces \eqref{LG_G_adj},
which would then allow one to compare the two actions by showing that the corresponding Borel spaces $EG \times_G -$, suitably defined as derived pre-stacks, are both equivalent to $\textbf{Map}(X,|B_{\bullet}G|^{\rm pre})$.
} 
with the usual one on $\DRep_G(X)$, so that we have
an equivalence \eqref{DRep_stack_MapXBG}.

\vspace*{1ex}

\begin{remark}
If $G$ is a reductive group, then in view of the equivalence \eqref{DRep_stack_MapXBG}, the derived character scheme $ \DChar_G(X) $ (as defined in Section~\ref{S4} of the present paper) is quasi-isomorphic to the (simplicial commutative) ring of global functions on the derived mapping stack $\textbf{Map}(X,BG)$.
\end{remark}

\subsection{Kapranov's construction}
We now turn to the comparison with another derived moduli space: 
the derived moduli space of $G$-local systems $  \boldsymbol{R} \Loc_G(X,*) $ on the \emph{pointed} space $(X,*)$ studied in  \cite{K}.
The construction of this derived moduli space is as follows.

First, one considers the underived case. Given any $X \in \sset$, one considers an affine scheme 
$\Hom_{\sset}(X,BG)$ defined as an end functor between $X : \Delta^{\op} \rightarrow \set$ 
and $BG : \Delta^{\op} \rightarrow \Aff_k$.
Dually, one can take the ring of functions $\mathcal{O}(BG) : \Delta \rightarrow \CommAlg_k$,
then we have
\[
\mathcal{O}(\, \Hom_{\sset}(X,BG) \,) \, = \, 
X \otimes_{\Delta} \mathcal{O}(BG) \, \, \in \CommAlg_k
\]
In the derived case, one can embed the category $\CommAlg_k$ either in $\DGCA_k^+$, as in \cite{K}, or in $\sCommAlg_k$, as we will do in the following, and consider
\begin{equation}  \label{HomXBG_hocoend}
\mathcal{O}(\, {\bm R}\Hom_{\sset}(X,BG) \,) \, = \, 
X \otimes_{\Delta}^{\bm L} \mathcal{O}(BG) \, \, \in \Ho(\sCommAlg_k)
\end{equation}

A word of caution is in order. Here we are taking the homotopy coend functor between 
the objects $X \in (\set)^{\Delta^{\op}}$ and $\mathcal{O}(BG) \in (\sCommAlg_k)^{\Delta}$,
both with the Reedy model structure.
In particular, weak equivalences in $(\set)^{\Delta^{\op}}$ are \emph{isomorphisms} of simplicial sets.
Thus, despite the notation, the homotopy coend \eqref{HomXBG_hocoend} is in general not invariant under weak equivalences of $X$.
Instead, 
one considers a natural action of $G^{X_0} := G \times \stackrel{X_0}{\ldots} \times G$ on the derived affine scheme ${\bm R}\Hom_{\sset}(X,BG)$, 
where $X_0$ is the set of vertices of $X$.
Then, as shown in \cite{K}, the action of $G^{X_0 \setminus \{\ast\}}$ on  ${\bm R}\Hom_{\sset}(X,BG)$ is free, and the quotient
${\bm R}\Hom_{\sset}(X,BG) / G^{X_0 \setminus \{\ast\}}$ is a homotopy invariant of $X$.
This defines the derived moduli space of $G$-local systems $  \boldsymbol{R} \Loc_G(X,*) $
on the pointed space $(X,\ast)$.

To compare with the derived representation scheme, we take a reduced simplicial set $X \in \sset_0$,
so that in particular, $  \boldsymbol{R} \Loc_G(X,*) = {\bm R}\Hom_{\sset}(X,BG)$ in this case.
To compute the homotopy coend \eqref{HomXBG_hocoend} and relate it to 
$\DRep_G(X)$, we take a specific resolution 
$\mathcal{O}(\overline{W}(RG)) \in (\sCommAlg_k)^{\Delta}$ of $\mathcal{O}(BG)$,
and show that the coend 
$X \otimes_{\Delta} \mathcal{O}(\overline{W}(RG))$
is \emph{isomorphic} to $\lgr(X)_G$ as simplicial commutative algebras (see \cite{BRY}).
While this resolution $\mathcal{O}(\overline{W}(RG))$ is not Reedy cofibrant (the latching morphisms are not cofibrations), it is nonetheless Reedy smooth, 
meaning that the latching morphisms are smooth extensions in a suitable sense.
It is shown in \cite{K} that resolutions by this more general class of objects can be used to compute the homotopy coend \eqref{HomXBG_hocoend}. As a consequence, one obtains an equivalence
\[
\DRep_G(X) \, \simeq \, \boldsymbol{R} \Loc_G(X,*)
\]
Combined with the equivalence \eqref{DRep_stack_MapXBG}, this gives a comparison 
between the two derived moduli spaces $\boldsymbol{R} \Loc_G(X,*)$ and $ \textbf{Map}(X,BG) $:
\bcor
\la{LocMap}
The derived mapping stack $\mbox{\rm\textbf{Map}}(X,BG)$ is equivalent to the quotient of the derived moduli space $\boldsymbol{R} \Loc_G(X,*)$ of $G$-local systems by the conjugation action of $G\,$: i.e., 
\[
\mbox{\rm \textbf{Map}}(X,BG) \, \simeq \, 
[ \, \boldsymbol{R} \Loc_G(X,*)/G \, ]
\]
\ecor
\subsection{Relation to Pridham's work}
The referee has suggested to us that the main results of this Appendix (namely, Proposition~\ref{DRep_scheme_MapXBG}
and the equivalence \eqref{DRep_stack_MapXBG}) can be also deduced from the work of J.~P.~Pridham \cite{Pri132}. Specifically, assuming that the constructions of \cite{Pri132} extend to higher Deligne-Mumford stacks, for any (finite) reduced simplicial set $X$, one can consider the derived moduli stack of $G$-torsors on the Deligne-Mumford hypergroupoid $\, \Spec(k) \times X \in \sAff\,$ defined by $\,(\Spec(k) \times X)_n := \Spec(k ^{\times X_n})$.
In this case, the functor $ {\rm C}_{\rm et}^{\bullet}(X,G)$ given in \cite[Definition 4.34]{Pri132} is simply represented by the cosimplicial group scheme $ {\rm C}_{\rm et}^{n}(X,G) := G^{X_n} $,
and the functor $\underline{{\rm MC}}[\underline{ {\rm C}_{\rm et}^{n}(X,G) }]$  in \cite[Definition~4.3]{Pri132} is represented by $\Spec[\lgr(X)_G]$. From Proposition 4.38 in {\it loc. cit.} we may then recover both our Proposition \ref{DRep_scheme_MapXBG} and the equivalence \eqref{DRep_stack_MapXBG}.

\end{document}